\documentclass[11pt]{article}
\usepackage[utf8]{inputenc}
\usepackage[T1]{fontenc}
\usepackage{amsmath,amssymb,amsthm,mathtools,enumerate}

\usepackage[round]{natbib}

\usepackage[margin=2.5cm]{geometry}
\usepackage{graphicx}
\usepackage{hyperref}
\usepackage[section]{placeins}

\newcommand{\EE}{\mathbb{E}}
\newcommand{\PP}{\mathbb{P}}
\newcommand{\NN}{\mathbb{N}}
\newcommand{\RR}{\mathbb{R}}
\newcommand{\x}{{\boldsymbol x}}
\newcommand{\z}{{\boldsymbol z}}
\newcommand{\ub}{{\boldsymbol u}}
\newcommand{\X}{{\textbf{X}}}
\newcommand{\1}{\mathbf{1}}
\newcommand{\dd}{\mathrm{d}}
\DeclareMathOperator{\Var}{Var}

\mathtoolsset{showonlyrefs,showmanualtags}

\makeatletter
\newcounter{condc}
\newcommand{\condlabel}[2]{\refstepcounter{condc}\protected@edef\@currentlabel{#1}\label{#2}}
\makeatother

\theoremstyle{plain}
\newtheorem{theorem}{Theorem}
\newtheorem{lemma}{Lemma}
\newtheorem{proposition}{Proposition}
\newtheorem{corollary}{Corollary}
\newtheorem{definition}{Definition}
\newtheorem{assum}{Assumption}
\theoremstyle{remark}
\newtheorem{remark}{Remark}

\usepackage{xcolor}

\title{Anytime-valid inference for conditional hazard and survival functions with streaming data}

\author{Eni Musta\footnote{Korteweg-de Vries Institute for Mathematics, University of Amsterdam, Amsterdam, Netherlands. E-mail: e.musta@uva.nl  }   \qquad \qquad Valentin Patilea\footnote{Univ. Rennes, Ensai, CREST-UMR 9194, F-35000, France.  Email: valentin.patilea@ensai.fr  }  }

\begin{document}
\maketitle

\begin{abstract}
We consider the nonparametric estimation of the conditional hazard, cumulative hazard and survival functions of a time-to-event variable from streaming data subject to censoring, truncation, or a cure fraction. A time-to-event model is a map from the distribution of the latent variables to that of the observed ones. For a large class of models this map has an explicit inverse, and the conditional hazard is the ratio of the density of the exactly observed durations to an at-risk function. Both are estimated by one-pass learners, in which each observation contributes once and is never revised, such as recursive kernel, series and histogram estimators. Our main result is an asymptotic confidence band sequence for the conditional hazard at a fixed covariate value, uniform in time and in the sample size, which allows monitoring after a burn-in, and stopping at a data-dependent time. The band can shrink at the optimal   uniform   rate up to a logarithmic factor. Bands for the cumulative hazard and the survival function follow. In simulations, a fixed-sample band inspected repeatedly fails in about half of the streams, while the miscoverage of the recommended bands stays below, or close to, the nominal level.
\end{abstract}


\medskip



\addtocontents{toc}{\protect\setcounter{tocdepth}{-1}}

\section{Introduction}\label{sec:intro}

Time-to-event data are increasingly collected in a continuous way: patients are enrolled one after the other in a registry, customers subscribe and leave a service, machines are put into service and fail. The natural estimator of the quantity of interest is then updated each time a new observation arrives, and the analyst may look at it at any moment, with the option to stop the data collection when the current estimate seems precise enough. This streaming setting raises a computational issue, since the estimator has to be updated at a small cost without storing the past data, and an inferential one, which is more delicate: a confidence interval built on a central limit theorem is valid at a fixed sample size, and looking at it repeatedly, stopping according to some rule, generally inflates the error probability, up to one. Group-sequential designs control this error when the number of looks, or the maximal information, is fixed in advance \citep{JT2000,LD1983}. Without such planning, the appropriate notion is that of a \emph{confidence sequence}, a sequence of intervals which contain the parameter simultaneously at all sample sizes with a prescribed probability, and which is therefore valid at any stopping time \citep{R1970,HRMS2021}. The \emph{asymptotic confidence sequences} of \cite{WASKR2024}, the time-uniform analogue of the intervals based on a central limit theorem, have made the approach widely   applicable. They rest upon   a strong Gaussian approximation of the partial sums of the estimator and a time-uniform boundary for the Gaussian process. So far, these ideas have been developed for means, finite-dimensional functionals and functions whose estimation involves no smoothing, such as the distribution function. We make them available for the conditional hazard function, which requires smoothing, with a tuning that evolves with the sample size, and which is estimated from   time-to-event   data, through the observation mechanism.

The duration of interest is not always observed, and sometimes never observed, because of censoring, truncation or long-term survivors. The viewpoint adopted here, described in Section~\ref{sec:map}, is that a model for the observation mechanism is a map from the   distribution of latent variables to that of the observables.   For a large class of mechanisms this map has an explicit inverse, which represents the conditional hazard as   the ratio of the joint density of the exactly observed durations and the covariates to an \emph{at-risk function}, a regression-type function of the observables, possibly affine in a nuisance parameter. Taking this representation as the starting point, we address two questions: \emph{which learners can be used for the numerator and the denominator, with a continuous updating of the estimator}, and \emph{how to build confidence bands for the hazard, the cumulative hazard and the survival function, uniform in time and valid at all sample sizes}. Uniformity in time matters for the questions asked of a hazard curve: where it peaks, whether it is monotone, whether the curves of two covariate values cross. Pointwise confidence sequences cannot answer them. For the survival function, a band gives intervals for the survival probabilities at all horizons and for its quantiles, simultaneously.

The answer to the first question is a structural property that we call the \emph{one-pass} property: the estimator is an average of contributions, one per observation, each computed when the observation arrives, with a tuning parameter rule fixed in advance, and never revised afterwards, or a smooth function of finitely many such averages. The recursive kernel estimators of \cite{WW1969} and \cite{Y1970} have this property, and so do recursive series estimators, with a number of basis functions growing with the index of the observation, and dyadic streaming histograms, with a resolution growing with the index. Besides the computational advantage, the one-pass property has a probabilistic consequence at the heart of our approach: the contributions form a single sequence of independent variables, not a triangular array as for a standard estimator, whose tuning is chosen at the final sample size and applied to all the past observations. This property underlies the laws of the iterated logarithm for recursive density estimators \citep{WD1979}. A strong Gaussian approximation couples one sequence of partial sums with one Gaussian sequence, simultaneously for all sample sizes; a triangular array offers no single sequence to couple. For curve estimation, where the tuning has to change with the sample size, the one-pass structure is thus the   natural route   to anytime-valid inference. Our results are stated under high-level conditions on the one-pass learner, verified for the three learners above.

The answer to the second question is our main contribution. For   an estimator of a function of time with the one-pass property, we construct bands which are uniform in time and in the sample size, which we call asymptotic confidence band sequences, extending the definition of \cite{WASKR2024}. In our application, the bands concern a fixed covariate value $\x$ and a compact interval $[t_0,t_1]$ of durations, with $t_0>0$ in general. The construction combines a strong Gaussian approximation valid for all sample sizes, for sums of independent variables of growing size, with Robbins' normal mixture boundary applied at the points of a nested dyadic grid; no multivariate strong approximation is needed. For the kernel and series learners, the band shrinks at the optimal uniform rate for a twice differentiable function of time and of the continuous covariates, up to a logarithmic factor; the rate-optimal bandwidth is admissible, although a slight undersmoothing is advisable in finite samples. As for all asymptotic confidence sequences, the guarantee for the computable band holds after a burn-in. Bands for the cumulative hazard and the survival function follow by integration, and are narrower.   Simulations illustrate the practical side, show that the estimation of the nuisance parameter has no visible cost, and lead us to recommend   a local-maximum estimator of the variance.

The third contribution is the catalogue of incompleteness mechanisms covered: right censoring, left truncation and right censoring, modified current status data \citep{PR2006b} and cure models \citep{AVK2018}, for which the bands concern the hazard, the cumulative hazard and the survival function of the susceptibles. 

Online and recursive estimation in survival analysis has been considered for the Cox model and parametric models with massive data \citep{WCSY2021,TS2022}, for the hazard rate without covariates under random censoring \citep{DMP2018}, and for conditional densities under censoring \citep{KS2014}; see also \cite{YWHYY2022} for subsampling with massive survival data. For the conditional hazard in survival analysis with a fixed sample, we refer to \cite{E2018,E2024,E2026} and the references therein. Recursive
kernel estimators of the conditional hazard, based on the representation of Section~\ref{sec:map}, were proposed by \cite{AP2026}; here we allow for a general class of one-pass learners and obtain anytime-valid inference. Confidence sequences go back to \cite{R1970}; see \cite{HRMS2021} for the nonasymptotic theory, \cite{WASKR2024} for the asymptotic one, \cite{WLR2025} for strong approximations valid for all sample sizes, and \cite{GK2025} for sharp asymptotic confidence sequences for a mean.   In survival analysis, anytime-valid inference has been developed for testing, with the safe logrank tests of \cite{GLS2021} and \cite{SPLG2024}, and for an average treatment effect on survival at a fixed horizon, in adaptive experiments \citep{WFSSJF2026}.   Time-uniform bands for distribution functions, which involve no smoothing, are obtained by \cite{HR2022}, by \cite{MH2023} for observations which may depend on the past, with a union bound over a dyadic grid, and by \cite{HZR2026}, with an extension to right censoring with a known censoring law. To the best of our knowledge, confidence sequences for a nonparametric conditional hazard or survival function, let alone bands uniform in time, have not been considered before.

Section~\ref{sec:map} presents the map viewpoint of a time-to-event model, and the hazard representations. Section~\ref{sec:learners} introduces the one-pass estimators, the high-level conditions on them, and three learners as examples. Section~\ref{sec:inference} contains the general results for   estimators satisfying these conditions,\ and their consequences for the hazard, the cumulative hazard and the survival function. Section~\ref{sec:examples} verifies the conditions for the three learners and gives the conditions on the tuning parameter rates, and Section~\ref{sec:mechanisms} treats the incompleteness mechanisms, including cure models. Section~\ref{sec:numerics} reports numerical experiments and Section~\ref{sec:discussion} discusses extensions. Proofs and further mechanisms are in the appendices.

\section{When a model is a map}\label{sec:map}

In time-to-event analysis, the variable of interest is a non-negative duration $T$ which, because of censoring, truncation or other incompleteness mechanisms, is not always observed, and sometimes never observed. Naive empirical estimators are then biased, and the bias does not disappear with large data sets or flexible learners. The correction comes from modeling how the observations are generated.

Our viewpoint is that, in many contexts, \emph{a model for the observation mechanism introduces latent variables and defines a map between the (conditional) distributions of the latent variables and the (conditional) distribution of the observables}. The identification assumptions, in general not testable, make this map one-to-one, and the program is then simple: learn the distribution of the observables, with any method, and map it back. In this way, in survival analysis one can learn the distribution of $T$ through appropriate functionals of the distribution of the observed data. This idea was used in \cite{AP2026} to build recursive estimators of the conditional hazard. Here, we allow a general class of one-pass learners, and build  anytime-valid inference. 

\subsection{Observables, the model and the inverse map}\label{subsec:observables}\label{subsec:inverting}

Let $\X$ be a vector of covariates, always observed. Instead of $T$, one observes a non-negative variable $Y$ and an indicator $I$ with values in a finite set of integers, such that $Y=T$ if and only if $I=1$. With left or right censoring, $I\in\{0,1\}$; with left truncation and right censoring, $I\in\{0,1\}$ and the truncation variable is also observed; with the modified current status data of \cite{PR2006b}, $I\in\{0,1,2\}$; with competing risks, $I$ is the cause of the event, and the framework applies after relabeling.   The conditional distribution of the observables $(Y,I)$ given $\X=\x$  is characterized by the sub-distribution functions $H_k(t\mid\x)=\PP(Y\leq t,\ I=k\mid\X=\x)$,   which can be estimated by any method, Bayesian or frequentist, with classical smoothers or machine learning. In all the cases considered below, the sub-distribution $H_1(\cdot\mid\x)$ of the exactly observed durations is absolutely continuous, with a sub-density whose integral is $\PP(I=1\mid\X=\x)$; multiplied by the density $f_\X(\x)$ of the covariates with respect to a dominating measure, it gives the joint density $f_{Y,I,\X}(t,1,\x)$ of $(Y,\X)$ on $\{I=1\}$, the numerator of all the representations below.

The standard example is right censoring, $Y=T\wedge C$ and $I=\1\{T\leq C\}$, with the identification assumption $T\perp C\mid\X$. The latent-to-observed map
 \[
 H_1(\dd t\mid\x)=S_{C\mid\X}(t-\mid\x)F_{T\mid\X}(\dd t\mid\x),\qquad H_0(\dd t\mid\x)=S_{T\mid\X}(t\mid\x)F_{C\mid\X}(\dd t\mid\x),
 \]
 with $F$ and $S=1-F$ the conditional distribution and survival functions, has the explicit inverse
 \[
 \Lambda_{T\mid\X}(\dd t\mid\x)=\frac{H_1(\dd t\mid\x)}{1-H(t-\mid\x)},\qquad S_{T\mid\X}(t\mid\x)=\prod_{s\in[0,t]}\big\{1-\Lambda_{T\mid\X}(\dd s\mid\x)\big\},
 \]
   where $H(t\mid\x)=H_0(t\mid\x)+H_1(t\mid\x)$, and   the last product is a product-integral \citep{GJ1990}.    If the conditional distribution of $T$ given $\X=\x$ admits a density, with $H_1(\dd t\mid\x)= h_1( t\mid\x) \dd t$, the inverse can be expressed using the conditional hazard function $ \lambda_{T\mid\X}$ as  
 \[
 \lambda_{T\mid\X}( t \!\mid\! \x)=\frac{h_1( t\!\mid\!\x)}{1-H(t-\!\mid\!\x)},\quad   \Lambda_{T\mid\X}( t\!\mid\!\x) = \int_0^t \!\! \lambda_{T\mid\X}( s\!\mid\!\x) \dd s, \quad  S_{T\mid\X}(t\!\mid\!\x)=\exp\Big(\!- \int_0^t \!\!\lambda_{T\mid\X}( s\!\mid\!\x) \dd s \!\Big).
 \]
  One then puts a hat on the unknown: empirical frequencies give the Kaplan--Meier estimator \citep{P1977}, kernel smoothers the conditional Kaplan--Meier estimator \citep{B1981,D1987,D1989}, and nearest  neighbors or neural networks the estimators of the recent machine learning literature.   For the nonparametric Bayesian perspective, see \citet[Ch.~8]{R2024}. 

  In summary, the map is given by the incompleteness mechanism (right censoring in the example above), invertibility 
is guaranteed by identification assumptions, such as (conditional) independence, and the learning method of the (conditional) distribution of the observables is a separate choice.

  Let us point out that in some cases, the data leave open the choice of the model and the map. Consider for example   observations where $T$ is either observed ($I=1$), or only known to be larger ($I=0$) or smaller ($I=2$) than $Y$. 
At least four models apply to such data: the left-right censoring model of \cite{T1974}, the two modified current status models of \cite{PR2006b} and the twice censoring models of \cite{PR2006a}. Their maps differ, and so do the distributions of $T$ they recover from the same distribution of $(Y,I)$. The choice between them is not testable,   and should be based on the specific application. See \cite{PR2006b,PR2006a} for some discussions. The choice of the model also matters for the theoretical study and the   computation: Turnbull's model has no explicit inverse map, whereas the others do.

\subsection{Hazard representations and the at-risk function}\label{subsec:representations}

Several types of censoring, truncation, long-term survivors and competing risks have an explicit inverse map, which can be written in a common form. When the conditional distribution of $T$ has a density,
\begin{equation}\label{eq:haz_rep}
	\lambda_{T\mid\X}(t\mid\x)=\frac{f_{Y,I,\X}(t,1,\x)}{R(t,\x)},
\end{equation}
where the \emph{at-risk (scalar) function} $R(t,\x)$ is a function of the distribution of the observables, typically a regression function times $f_\X(\x)$: $R(t,\x)=\PP(Y\geq t\! \mid\! \X\!=\!\x)f_\X(\x)$ under right censoring, and $R(t,\x)=\PP(L\leq t\leq Y\! \mid\! \X\!=\!\x)f_\X(\x)$ under left truncation and right censoring, with $L$ the observed truncation variable. For modified current status data, $R$ is \emph{affine} in a nuisance function $\eta(\x)$,
\begin{equation}\label{eq:R_affine_intro}
	R(t,\x)=R_1(t,\x)+\eta(\x)R_2(t,\x),
\end{equation}
with $R_1,R_2$ again regression-type functions of the observables. Cure models can be written in this form, but with $R_2$ not of regression type. They are handled through a smooth transformation of the hazard of the improper duration; see Section~\ref{subsec:cure}. The numerator is a density in time and in the covariates on the event $\{I=1\}$, whereas the denominator is a regression-type function which requires no smoothing in time; it therefore converges faster, and inference on the hazard is essentially inference on the numerator. The cumulative hazard and the survival function follow by $\Lambda_{T\mid\X}(t\mid\x)=\int_0^t\lambda_{T\mid\X}(s\mid\x)\dd s$ and $S_{T\mid\X}(t\mid\x)=\exp\{-\Lambda_{T\mid\X}(t\mid\x)\}$. The representations \eqref{eq:haz_rep} are derived in Section~\ref{sec:mechanisms} and Appendix~\ref{app:mechanisms}; see also \cite{AP2026}, where some of them first appeared, and \citet[Chapter~6]{E2018} for the case without covariates.   Note that the choice of the dominating measure for $f_{\X}$ is flexible, allowing for continuous and discrete covariates.  

The representation \eqref{eq:haz_rep} also delimits the scope of the approach: it requires exact observations, $I=1$ with positive probability and, for the inference, $f_{Y,I,\X}(t,1,\x)>0$ on the time interval of interest. This rules out interval censoring, where $T$ is never observed exactly, and makes the hazard, rather than the distribution function, the natural target: it is a density divided by a regression function, two standard objects for which one-pass learners are available.


\section{One-pass estimation}\label{sec:learners}

This section introduces the one-pass structure and the conditions on the learner under which the inference of Section~\ref{sec:inference} is derived, and describes three learners satisfying them: recursive kernel estimators, recursive series estimators and dyadic streaming histograms. The conditions are verified for these learners in Section~\ref{sec:examples}.   Below, $C$, $C'$ and $\mathfrak C$ denote generic constants, whose value may change from line to line. Constants with a subscript, such as $C_M$, $C_V$ or $C_q$, have a fixed meaning. The censoring variable is denoted $C$ or $C_i$; it only appears in the description of the models, where no confusion with a constant can arise.

\subsection{Setting and notation}\label{subsec:setting}

The observations $Z_1,Z_2,\ldots$ are independent copies of a random vector $Z$, arriving one at a time. The index $i$ is the order in which the observations are processed, and this order is assumed not to depend on their values; with staggered entry, this is the case when the records are processed in the order of enrollment; see Section~\ref{sec:discussion}. In the time-to-event context, $Z=(Y,I,\X)$, where $Y\geq0$ is the observed duration, $I$ is an indicator taking values in a finite set of non-negative integers containing $1$, with $Y=T$ if and only if $I=1$, and $\X=(\X_c,\X_d)$ is the covariate vector, with $\X_c$ taking values in a compact rectangle $\mathcal X_c\subset\RR^{d_c}$ and $\X_d$ in a finite set $\mathcal X_d$; $Z$ is augmented with a truncation variable $L\geq0$ or a censoring variable $C\geq0$ when these are observed. Hence $Z$ takes values in a subset $\mathcal Z$ of a Euclidean space. A generic point of $\mathcal Z$ is written $\mathfrak z=(y,\iota,\z)$, or $\mathfrak z=(y,\iota,\ell,\z)$ when a truncation variable is observed, with $\z=(\z_c,\z_d)$: the components $y$, $\iota$, $\ell$ and $\z$ are the values taken by $Y$, $I$, $L$ and $\X$. Functions on $\mathcal Z$, such as the contributions $\psi_i(t;\cdot)$ below or the functions $u_m(t;\cdot)$ of   Assumption~\ref{ass:data}, are written in terms of these components. The covariate value used in conditioning is always denoted $\x=(\x_c,\x_d)$, and it is \emph{fixed} throughout.  All functions of the time argument are considered on a compact interval
\[
\mathcal T=[t_0,t_1],\qquad 0\leq t_0<t_1<\tau_0,
\]
where $\tau_0$ is a known constant. The bands of Section~\ref{sec:inference} are constructed on intervals with $t_0>0$, for the hazard function as well as for the cumulative hazard, but the estimator of the latter integrates the hazard estimator from the origin, so that the conditions on the learner are then required on $[0,t_1]$. The targets of this section are generic real-valued functions $g$ on $\mathcal T$; in the following sections, $g$ will be $f_{Y,I,\X}(\cdot,1,\x)$, one of the components $R_1(\cdot,\x)$, $R_2(\cdot,\x)$ of the at-risk function, or an integrated version of these. We write $f(t,1,\x)$ for $f_{Y,I,\X}(t,1,\x)$ when no confusion is possible. Finally, for $a\in(0,1)$, let
\begin{equation}\label{eq:w_n}
	w_n(a)=n^{-(1-a)/2}\sqrt{\log(n+1)},
\end{equation}
the order of the stochastic error of an estimator whose $i$-th contribution has variance of order $i^{a}$, and, as shown in Section~\ref{sec:inference}, of the width of the anytime-valid bands.   Throughout, $\log$ denotes the natural logarithm and $\log_2$ the logarithm in base $2$.  

\subsection{One-pass estimators}\label{subsec:onepass}

\begin{definition}\label{def:onepass}
	For a given $\x$, consider a sequence of estimators $(\hat g_n)_{n\geq1}$ of $g$ and the representation
	\begin{equation}\label{eq:onepass}
		\hat g_n(t)=\frac1n\sum_{i=1}^n\psi_i(t;Z_i)+\mathfrak r_n(t),\qquad t\in\mathcal T,\ n\geq1,
	\end{equation}
	where, for each $i\geq1$, $\psi_i:\mathcal T\times\mathcal Z\to\RR$ is a deterministic measurable function.
	\begin{enumerate}
		\item[(a)] $(\hat g_n)$ is a \emph{one-pass average} if \eqref{eq:onepass} holds with $\mathfrak r_n\equiv0$ and known functions $\psi_i$, each fixed before $Z_i$ is observed.
		
		\item[(b)] $(\hat g_n)$ is a \emph{one-pass estimator} if $\hat g_n$ is a known function of finitely many scalar one-pass averages, all fixed in advance, whose number may grow with $n$.
		
		\item[(c)] $(\hat g_n)$ \emph{admits a one-pass expansion} if \eqref{eq:onepass} holds for some functions $\psi_i$, possibly depending on unknown features of the distribution of $Z$. The content of this notion lies in the conditions of Section~\ref{subsec:conditions}, which require a negligible remainder.
	\end{enumerate}
	The \emph{contribution} of the $i$-th observation is $\psi_i(\cdot;Z_i)$; the centered contributions, their partial sums and the smoothed target are
	\[
	\xi_i(t)=\psi_i(t;Z_i)-\EE[\psi_i(t;Z_i)],\qquad \mathfrak S_n(t)=\sum_{i=1}^n\xi_i(t),\qquad \bar g_n(t)=\frac1n\sum_{i=1}^n\EE[\psi_i(t;Z_i)].
	\]
\end{definition}

The functions $\psi_i$ depend on $\x$, but this dependence is omitted from Definition~\ref{def:onepass} for simplicity. The random functions of $t$ considered here are continuous except at finitely many points, which are measurable functions of the data or the grid points. Suprema over $\mathcal T$ are therefore measurable.

A one-pass average is updated recursively, $\hat g_n=(1-1/n)\,\hat g_{n-1}+n^{-1}\psi_n(\cdot;Z_n)$, so that each observation is processed once and its contribution is never revised. A one-pass estimator is computed from the recursions of its averages, without storing the past data; its memory is the number of scalar averages, such as the values of one-pass averages at fixed time points, the coefficients of a series expansion or the counts of a histogram. A one-pass expansion is a property of an estimator, used for the inference, not a way of computing it: its contributions may involve unknown quantities, and its remainder is defined by difference. The estimators proposed in this paper are smooth functions of such averages: ratio estimators, local linear smoothers or plug-in nuisance parameters; for the kernel learner,   the estimators of the at-risk function\ and the integrals of ratios, the averages are the values on a fixed grid of time points (Remark~\ref{rem:memory}). They admit a one-pass expansion with a negligible remainder under the conditions of Remark~\ref{rem:remainder} and Lemma~\ref{lem:int_rate}. An estimator may also admit such an expansion without being a one-pass estimator: the exact Nelson--Aalen estimator of Appendix~\ref{app:dc0} uses the at-risk average at the observed times, which are not fixed in advance, and becomes a one-pass estimator only on a fixed grid. Since the $\psi_i$ are deterministic and the $Z_i$ are independent, for each $t$ the centered contributions $\xi_1(t),\xi_2(t),\ldots$ form a \emph{single sequence} of independent variables, which is what the strong approximations of Section~\ref{sec:inference} require. For a standard estimator, whose tuning parameter chosen at time $n$ is applied to all the past observations, the $i$-th contribution changes with $n$.

\begin{remark}
	\label{rem:aux_rand}
	For simplicity, we allocate the learning rate $i^{-1}$ to the observation $Z_i$ and normalize the partial sum defining $\hat g_n$ by $n^{-1}$. Different learning rates can be considered with weights depending on $i$ in $\psi_i$ and suitably normalizing their sum to define $\hat g_n$; see \cite{AP2026}. Moreover, the $\psi_i$ may also depend on an auxiliary random element $\boldsymbol\varpi$ (random partitions, random features), independent of the data and drawn once at the beginning. The conditions below are then required, and the results of Section~\ref{sec:inference} hold, conditionally on $\boldsymbol\varpi$, hence also unconditionally.
\end{remark}

\subsection{Conditions}\label{subsec:conditions}

Let $a\in(0,1)$ and $\mathfrak b>0$ be two exponents, and write $w_n=w_n(a)$. The exponent $a$ measures the growth of the contributions; the exponent $\mathfrak b$ is the rate of the bias.   The target $g$ is bounded on $\mathcal T$.  

\medskip\noindent
\textbf{Condition (I)} (one-pass structure). $(\hat g_n)$   admits a one-pass expansion   \eqref{eq:onepass} for which there exist constants $C_M,C_V$ such that, 
\begin{equation}\tag{I}\label{cond:I}
	\sup_{t\in\mathcal T}|\xi_i(t)|\leq C_M i^{a}\ \text{ a.s.},\quad \sup_{t\in\mathcal T}\Var\{\xi_i(t)\}\leq C_V i^{a}\quad\forall i\geq1,\quad \text{and } \quad \sup_{t\in\mathcal T}|\mathfrak r_n(t)|=o(w_n)\ \text{ a.s.}
\end{equation}

\medskip\noindent
\textbf{Condition (V)} (variance profile). There exist a   given, deterministic sequence $(\mathfrak V_n)$,   with $c_0n^{1+a}\leq\mathfrak V_n\leq C_0n^{1+a}$ for some $0<c_0\leq C_0$, a continuous function $v:\mathcal T\to(0,\infty)$, and computable functions $\hat v_n$, such that, with $V_n(t)=\sum_{i=1}^n\Var\{\xi_i(t)\}$,
\begin{equation}\tag{V}\label{cond:V}
	\sup_{t\in\mathcal T}\left|\frac{V_n(t)}{\mathfrak V_n\,v(t)}-1\right|\to0,\qquad \sup_{t\in\mathcal T}\left|\frac{\hat v_n(t)}{v(t)}-1\right|\to0\ \text{ a.s.}
\end{equation}

\medskip\noindent
\textbf{Condition (B)} (bias).
\begin{equation}\tag{B}\label{cond:B}
	\sup_{t\in\mathcal T}|\bar g_n(t)-g(t)|=O(n^{-\mathfrak b}).
\end{equation}

\medskip\noindent
\textbf{Condition (E)} (uniform rate).
\begin{equation}\tag{E}\label{cond:E}
	\sup_{t\in\mathcal T}|\hat g_n(t)-g(t)|=O(n^{-\mathfrak b}+w_n)\ \text{ a.s.}
\end{equation}

\medskip
A last condition, \eqref{cond:M}, links the learner to the grid on which the bands are calibrated; it is stated in Section~\ref{subsubsec:grid}. Conditions \eqref{cond:I}, \eqref{cond:V}, \eqref{cond:B} and \eqref{cond:M} are required for the function whose band is constructed, that is, the numerator of the hazard representation and, in integrated form, the cumulative hazard. For the denominators, only \eqref{cond:E} and the Lipschitz property of the target are needed, with a rate negligible with respect to the width of the band. Condition \eqref{cond:E} follows from \eqref{cond:I} and \eqref{cond:B} under a mild regularity of the contributions in $t$: Lemma~\ref{lemma:chaining} in Appendix~\ref{app:learners} gives sufficient conditions for \eqref{cond:E} and \eqref{cond:M}, covering smooth contributions, contributions of bounded variation in $t$, such as the indicators of the at-risk function, and piecewise constant contributions, as in histograms. For the histogram learner of Section~\ref{subsubsec:histogram}, $\mathfrak V_n=\tau_0^{-1}\sum_{i\leq n}2^{k_ip}$ is not asymptotically equivalent to a power of $n$. This is why Condition~\eqref{cond:V} is formulated with a general deterministic sequence.   The order of magnitude of the tuning parameters, the memory and computation burden, and other learners are discussed in Appendix~\ref{app:learners_more}.

\begin{remark}[Estimators with a remainder]\label{rem:remainder}
	Let $\hat g_n=\Gamma(\hat g_{1,n},\ldots,\hat g_{q,n})$, where the $\hat g_{k,n}$ are   one-pass averages   of functions $g_k$, each satisfying \eqref{cond:I} with the same $a$ and \eqref{cond:E}, and $\Gamma$ is twice continuously differentiable near $\{(g_1(t),\ldots,g_q(t)):t\in\mathcal T\}$. A Taylor expansion shows that $\hat g_n$   admits a one-pass expansion for   $g=\Gamma(g_1,\ldots,g_q)$, with contributions $\psi_i=g+\sum_k\partial_k\Gamma(g_1,\ldots,g_q)(\psi_{k,i}-g_k)$ and a remainder $O(n^{-2\mathfrak b}+w_n^2)$ a.s., uniformly in $t$, which is $o(w_n)$ as soon as $2\mathfrak b \ge (1-a)/2$. This covers local linear smoothers, smooth functions of recursive weighted moments, ratio estimators such as $\hat f_n/\hat R_n$, and a plug-in nuisance parameter: if $\hat R_n=\hat R_{1,n}+\hat\eta_n\hat R_{2,n}$ with $|\hat\eta_n-\eta|=o(w_n)$ a.s., then $\hat R_n$ differs from the   one-pass average   $\hat R_{1,n}+\eta\hat R_{2,n}$ by $o(w_n)$, uniformly in $t$. A nuisance parameter whose error is of the order of a standard deviation is discussed in Remark~\ref{rem:nuisance_slow}.
\end{remark}

\subsection{Three learners}\label{subsec:learners}

  We now describe three   one-pass averages,   considering as target the sub-density $g(t)=f(t,1,\x)$, $t\in\mathcal T$.   In all three cases, the contribution of   $Z_i=(Y_i, I_i, (\X_{i,c},\X_{i,d}))$  is nonzero only if $I_i=1$.   The learners differ in the way the time argument is handled: kernel smoothing, series expansion or binning. The tuning parameters are indexed by $i$, and chosen before $Z_i$ arrives; they decrease (bandwidths) or increase (number of basis functions, resolution) at a polynomial rate governed by a single exponent $\gamma$. The estimators of   the at-risk function   are described in Section~\ref{subsec:denominator}; they only involve smoothing in   $\z_c$,   so that the discussion below applies with $d_c+1$ replaced by $d_c$. 

Throughout, $\boldsymbol K$ is a product kernel on $\RR^{d_c}$ built from a symmetric, Lipschitz continuous probability density $K$ supported on $[-1,1]$   with Lipschitz constant $L_K$,   $\boldsymbol K_h(\ub)=h^{-d_c}\boldsymbol K(\ub/h)$, $K_b(u)=b^{-1}K(u/b)$ and $c_K=\int K^2(u)\dd u$.   Let $p=d_c+1$ be   the number of continuous arguments of $f(\cdot,1,\cdot)$. In Section~\ref{sec:examples}, we verify that the three learners satisfy the conditions of Section~\ref{subsec:conditions}, and Condition~\eqref{cond:M} of Section~\ref{subsubsec:grid} for a grid finer than the current bandwidth, with $a=\gamma p$, $\mathfrak b=2\gamma$ for the kernel and series learners under second-order smoothness, $\mathfrak b=\gamma$ for the histogram under Lipschitz smoothness, and a variance profile proportional to $f(\cdot,1,\x)$, so that $\hat v_n=\mathfrak c\hat f_n$ for a known constant $\mathfrak c$.

\subsubsection{Recursive kernel estimator}\label{subsubsec:kernel}

  For given $\x=(\x_c,\x_d)$, let   $b(i)=c_bi^{-\gamma}$ and $h(i)=c_hi^{-\gamma}$, $i\geq1$, with $\gamma\in(0,1/p)$ and $c_b,c_h>0$, and
\begin{equation}\label{eq:psi_kernel}
\psi^{\mathrm{K}}_i(t;Z_i)=   \1\{I_i=1\}\,K_{b(i)} (Y_i -t)\,\boldsymbol K_{h(i)}(\X_{i,c}-\x_c)\,\1\{\X_{i,d}=\x_d\},  
\end{equation}
so that $\hat f^{\mathrm K}_n(t,1,\x)=n^{-1}\sum_{i\leq n}\psi^{\mathrm K}_i(t;Z_i)$,   with $\psi^{\mathrm K}_i $ depending on $\x$.   This is the estimator considered by \cite{AP2026}.  

\begin{remark}[Boundary correction]\label{rem:boundary}
When $t_0=0$, the kernel $K_{b(i)}(y-t)$ has a bias of order $1$, not $b(i)^2$, for $t<b(i)$: the window $[t-b(i),t+b(i)]$ extends beyond the origin, where the sub-density vanishes. A standard remedy is a boundary kernel, which replaces $K_b(y-t)$ for $t\in[0,b)$ by a linear combination of $K_b(y-t)$ and $\{(y-t)/b\}K_b(y-t)$ with coefficients depending on $t/b$, chosen so that the first two moments on $[-t/b,1]$ are those of a symmetric kernel \citep[see][]{J1993}. The boundary-corrected contributions are again of the form \eqref{eq:psi_kernel},
and Proposition~\ref{prop:kernel}(\ref{prop:kernel:ii}) shows that they satisfy Conditions \eqref{cond:I}, \eqref{cond:B} and \eqref{cond:E} on $[0,t_1]$ with the same exponents. Condition \eqref{cond:V} does not hold uniformly up to $0$, because the variance of the corrected kernel differs from $c_K$ in the boundary region, but \eqref{cond:V} is not needed for the cumulative hazard, whose variance profile is governed by the integrated contributions; see Section~\ref{subsec:cumhaz}. Without correction, the bias of the integrated estimator on $[0,b(i)]$ is of order $b(i)$, which forces the stronger undersmoothing $\gamma\geq1/(d_c+2)$ and a slower rate for the band of the cumulative hazard. The series and histogram learners below need no   such correction.   
\end{remark}

\subsubsection{Recursive series estimator}\label{subsubsec:series}

Let $\{\varphi_j\}_{j\geq0}$ be an orthonormal basis of $L^2[0,\tau_0]$ and $D_J(y,t)=\sum_{j=0}^{J}\varphi_j(y)\varphi_j(t)$ the kernel of the orthogonal projection $S_J$ on the span of $\varphi_0,\ldots,\varphi_J$. Let $J(i)=\lfloor c_Ji^{\gamma}\rfloor\vee 1$ and $h(i)=c_hi^{-\gamma}$, with $\gamma\in(0,1/p)$ and $c_J,c_h>0$, and
\begin{equation}\label{eq:psi_series}
  \psi^{\mathrm S}_i(t;Z_i)=\1\{I_i=1,\ Y_i\leq\tau_0\}\,D_{J(i)}(Y_i,t)\,\boldsymbol K_{h(i)}(\X_{i,c}-\x_c)\,\1\{\X_{i,d}=\x_d\}.   
\end{equation}
The resulting estimator can be written as a series,
\begin{align*}
\hat f^{\mathrm S}_n(t,1,\x)&=\sum_{j=0}^{J(n)}\hat\theta_{n,j}\,\varphi_j(t),\\
\hat\theta_{n,j}&=\frac1n\sum_{i\leq n:\,J(i)\geq j}\1\{I_i=1,\,Y_i\leq\tau_0\}\,\varphi_j(Y_i)\,\boldsymbol K_{h(i)}(\X_{i,c}-\x_c)\,\1\{\X_{i,d}=\x_d\},
\end{align*}
where the $j$-th coefficient only uses the observations which arrived after the number of basis functions reached $j$. The coefficients are updated in $O(J(n))$ operations per observation and the number of terms increases at the polynomial rate $i^{\gamma}$, the series analogue of the inverse bandwidth $1/b(i)$; a rule where $J(i)$ is constant on dyadic blocks of $i$ is equally admissible, since only the order of $J(i)$ matters. The observations with   $Y_i>\tau_0>t_1$ do not contribute,   which does not affect the target.

We use   an augmented half-cosine basis obtained by   the Gram--Schmidt orthonormalization in $L^2[0,\tau_0]$ of the dictionary
\[
(1,\ t,\ t^2,\ c_1,\ c_2,\ \ldots),\qquad c_j(t)=(2/\tau_0)^{1/2}\cos(\pi jt/\tau_0),\quad j\geq1,
\]
  see \cite{E1999}.   The first three elements span the quadratic polynomials, and $\varphi_j$, $j\geq3$, is $c_{j-2}$ orthogonalized against   the previous elements of the dictionary    and normalized. The properties of this basis used in Section~\ref{sec:examples} are the following, with $\mathcal C^{2,1}$ the set of functions on $[0,\tau_0]$ with two derivatives, the second one being Lipschitz continuous.
\begin{enumerate}
\item[(S1)]\condlabel{S1}{cond:S1} (Boundedness) $\sup_j\|\varphi_j\|_\infty<\infty$ and $\|\varphi_j'\|_\infty\leq \mathfrak C j$ for $j\geq1$.
\item[(S2)]\condlabel{S2}{cond:S2} (Approximation) For $q\in\mathcal C^{2,1}$, $\sup_{t\in[0,\tau_0]}|S_Jq(t)-q(t)|\leq C_qJ^{-2}$, where $C_q$ depends on $q$ only through $\|q\|_\infty$, $\|q'\|_\infty$, $\|q''\|_\infty$ and the Lipschitz constant of $q''$.
\item[(S3)]\condlabel{S3}{cond:S3} (Localization) For every compact $\mathcal T\subset(0,\tau_0)$ and every continuous $q$ on $[0,\tau_0]$,
\[
\sup_{t\in\mathcal T}\Big|\frac1J\int_0^{\tau_0}D_J(y,t)^2q(y)\dd y-\frac{q(t)}{\tau_0}\Big|\to0\qquad\text{as }J\to\infty,
\]
uniformly over $q$ in a set which is uniformly bounded and uniformly equicontinuous.
\end{enumerate}
  These properties are discussed, with references and proofs, in Appendix~\ref{app:basis}.

\subsubsection{Dyadic streaming histogram}\label{subsubsec:histogram}

Binning is the simplest learner, and the building block of tree-based methods and of the streaming summaries used in large-scale learning systems. Assume, after an affine transformation of $\X_c$, that $\mathcal X_c\subset[0,1)^{d_c}$. For $k\geq0$, let $\mathcal C_k$ be the partition of $[0,\tau_0)\times[0,1)^{d_c}$ into the cells
\[
A=[\tau_0m_02^{-k},\tau_0(m_0+1)2^{-k})\times\prod_{r=1}^{d_c}[m_r2^{-k},(m_r+1)2^{-k}),\qquad (m_0,\ldots,m_{d_c})\in\NN^{p},
\]
of volume $|A|=\tau_02^{-kp}$, and let $A_k(t)$ denote the cell of $\mathcal C_k$ containing $(t,\x_c)$. Let $k_i=\lfloor\gamma\log_2 i\rfloor$, with $\gamma\in(0,1/p)$, and
\begin{equation}\label{eq:psi_hist}
  	\psi^{\mathrm H}_i(t;Z_i)=\1\{I_i=1\}\,\frac{\1\{(Y_i,\X_{i,c})\in A_{k_i}(t)\}}{|A_{k_i}(t)|}\,\1\{\X_{i,d}=\x_d\}.  
\end{equation}
The estimator $\hat f^{\mathrm H}_n(t,1,\x)=n^{-1}\sum_{i\leq n}\psi^{\mathrm H}_i(t;Z_i)$ is a histogram whose resolution doubles each time $i$ is multiplied by $2^{1/\gamma}$. The counts of the current resolution are updated in $O(1)$ operations per observation, and the counts of the previous resolutions are not revised. The partitions are nested, and the estimator is piecewise constant on the cells of $\mathcal C_{k_n}$, the finest partition used up to time $n$. 

A histogram has a first order bias, so that the natural smoothness class is that of Lipschitz sub-densities, instead of twice differentiable ones.

\begin{remark}[Memory of the learners]\label{rem:memory}
	The series and histogram learners are one-pass estimators in the sense of Definition~\ref{def:onepass}(b). Their contributions are of the form $\psi_i(t;\mathfrak z)=\sum_{j\leq m_i}\theta_{i,j}(\mathfrak z)e_j(t)$, with fixed functions $e_j$ and $m_i$ nondecreasing in $i$: the basis functions $\varphi_j$, with $m_i=J(i)$, or the indicators $t\mapsto\1\{(t,\x_c)\in A\}$ of the cells $A$ of resolution at most $k_i$, enumerated by increasing resolution. Hence $\hat f_n(t,1,\x)=\sum_{j\leq m_n}\hat\theta_{n,j}e_j(t)$, where $\hat\theta_{n,j}=n^{-1}\sum_{i\leq n}\theta_{i,j}(Z_i)$ are scalar one-pass averages. Since $\theta_{i,j}=0$ for $j>m_i$, the average $\hat\theta_{n,j}$ only receives contributions from the observations arriving after $m_i$ has reached $j$, so that the number $m_n$ of averages grows with $n$ without requiring the past data. The kernel learner is a one-pass average, but storing it on the continuum $\mathcal T$ would require the data; its values at time points fixed in advance are scalar one-pass averages. The estimators of the at-risk function of Section~\ref{subsec:denominator} are stored in the same way, and the integrals of Section~\ref{subsec:cumhaz} are approximated by quadrature on such points; see Appendix~\ref{subsubsec:memory}.
\end{remark}


\section{Anytime-valid inference}\label{sec:inference}

This section contains the inferential results. Section~\ref{subsec:defs} recalls the asymptotic confidence sequences of \cite{WASKR2024} and extends them to bands uniform in $t$. Moreover, it states the three probabilistic tools: a strong approximation with exponential tails, valid for all $n$ simultaneously, for sums of independent variables of growing size; a gluing lemma, which puts countably many such approximations on one probability space; and Robbins' normal mixture boundary. The uniformity in $n$ of the tail bound allows a union bound, on the Gaussian side, over the points of a nested dyadic grid introduced in Section~\ref{subsubsec:grid}, so that no multivariate strong approximation is needed.
The general results for   estimators satisfying the conditions of Section~\ref{subsec:conditions}\ are given in Section~\ref{subsec:general}, and they are specialized to the hazard, the cumulative hazard and the survival function in Sections~\ref{subsec:hazard} and \ref{subsec:cumhaz}. All proofs are in Appendix~\ref{app:inference}.

\subsection{Confidence sequences and band sequences: definitions and tools}\label{subsec:defs}

Let $\theta$ be an unknown real parameter and $\hat\theta_n$ an estimator computed from $Z_1,\ldots,Z_n$. Following \citet[Definition 2.1]{WASKR2024}, a sequence of intervals $[\hat\theta_n-\hat L_n,\hat\theta_n+\hat U_n]$, $n\geq1$, with computable $\hat L_n,\hat U_n>0$, is a \emph{$(1-\alpha)$-asymptotic confidence sequence} for $\theta$ if, on a possibly enriched probability space, there exist $L^*_n,U^*_n$ such that
\begin{equation}\label{eq:def_CS}
\frac{L^*_n}{\hat L_n}\xrightarrow{a.s.}1,\qquad \frac{U^*_n}{\hat U_n}\xrightarrow{a.s.}1,\qquad\text{and}\qquad \PP\big(\theta\in[\hat\theta_n-L^*_n,\hat\theta_n+U^*_n]\ \text{ for all }n\geq1\big)\geq1-\alpha.
\end{equation}
The intervals $[\hat\theta_n-L^*_n,\hat\theta_n+U^*_n]$ form a nonasymptotic confidence sequence, which is typically not computable; the computable bounds $\hat L_n,\hat U_n$ are asymptotically equivalent to the unknown ones. The guarantee is uniform in $n$: the probability that the sequence with bounds $L^*_n,U^*_n$ ever excludes $\theta$ is at most $\alpha$, so that these intervals can be inspected and the data collection stopped at any data-dependent time without inflating the error probability. The computable intervals inherit this property asymptotically, in the sense of Proposition~\ref{prop:stopping} below.

For a function $\theta(t)$, $t\in\mathcal T$, we consider bands $[\hat\theta_n(t)-\hat L_n(t),\hat\theta_n(t)+\hat U_n(t)]$ with computable $\hat L_n(t),\hat U_n(t)>0$, and say that they form a \emph{$(1-\alpha)$-asymptotic confidence band sequence} for $\theta(\cdot)$ on $\mathcal T$ if, on a possibly enriched probability space, there exist $L^*_n(t),U^*_n(t)$ such that
\begin{equation}\label{eq:def_band_ratio}
\sup_{t\in\mathcal T}\Big|\frac{L^*_n(t)}{\hat L_n(t)}-1\Big|\xrightarrow{a.s.}0,\qquad \sup_{t\in\mathcal T}\Big|\frac{U^*_n(t)}{\hat U_n(t)}-1\Big|\xrightarrow{a.s.}0,
\end{equation}
and
\begin{equation}\label{eq:def_band}
\PP\Big(\theta(t)\in[\hat\theta_n(t)-L^*_n(t),\hat\theta_n(t)+U^*_n(t)]\ \text{ for all }n\geq1\text{ and all }t\in\mathcal T\Big)\geq1-\alpha.
\end{equation}
The guarantee is now uniform in $n$ \emph{and} in $t$. It concerns the unknown bounds $L^*_n,U^*_n$. 
For the computable bands, it has the following operational meaning.

\begin{proposition}[Monitoring and stopping]\label{prop:stopping}
Let $[\hat\theta_n(t)-\hat L_n(t),\hat\theta_n(t)+\hat U_n(t)]$ be a $(1-\alpha)$-asymptotic confidence band sequence for $\theta(\cdot)$ on $\mathcal T$. Then, for every $\varepsilon>0$,
\begin{equation}\label{eq:stopping}
\liminf_{m\to\infty}\PP\Big(\theta(t)\in[\hat\theta_n(t)-(1+\varepsilon)\hat L_n(t),\hat\theta_n(t)+(1+\varepsilon)\hat U_n(t)]\ \text{ for all }n\geq m\text{ and all }t\in\mathcal T\Big)\geq1-\alpha.
\end{equation}
In particular, the same lower bound holds at any almost surely finite random times $\tau_m\geq m$, stopping times or not.
\end{proposition}

Thus the bands can be monitored after a burn-in period, and the data collection stopped at a data-dependent time. No guarantee is given for small $n$. The finite-sample behavior of the uninflated
bands is assessed in Section~\ref{sec:numerics}.

Let us now present the probabilistic tools used for the construction of our anytime valid inference.

\paragraph{Strong approximation with exponential tails.} The following is Sakhanenko's refinement of the Koml\'os--Major--Tusn\'ady construction for independent, non-identically distributed variables with exponential moments; see \citet[Th.~1]{Z2013}, where the condition is written $\mathcal L(\xi_j)\in S_1(\tau)$ with $\tau=1/\lambda$.

\begin{lemma}[\citealp{S1984}]\label{lemma:sakhanenko}
Let $\xi_1,\ldots,\xi_n$ be independent, centered random variables, and let $\lambda>0$ be such that
\[
\lambda\,\EE\big[|\xi_j|^3\exp (\lambda|\xi_j|)\big]\leq\Var(\xi_j),\qquad j=1,\ldots,n.
\]
Then, on a possibly enriched probability space, there exist independent Gaussian variables $\eta_1,\ldots,\eta_n$, with $\EE\eta_j=0$ and $\Var(\eta_j)=\Var(\xi_j)$, such that
\[
\EE\exp\Big(c\lambda\max_{m\leq n}\Big|\sum_{j=1}^m(\xi_j-\eta_j)\Big|\Big)\leq 1+\lambda B_n,\qquad B_n^2=\sum_{j=1}^n\Var(\xi_j),
\]
where $c>0$ is an absolute constant. In particular, 
$$
\PP\bigg(\max_{m\leq n}|\sum_{j\leq m}(\xi_j-\eta_j)|\geq z\bigg)\leq(1+\lambda B_n)\exp (-c\lambda z), \quad  \text{for any $z>0$}.
$$
\end{lemma}

The construction of Lemma~\ref{lemma:sakhanenko} depends on $n$. A construction valid for all $n$ simultaneously is obtained by combining independent constructions on blocks of indices; see \citet[p.~727]{Z2013} and the proof of Theorem 2.3 in \cite{WLR2025}. Lemma~\ref{lemma:blocked} adapts this construction to summands growing polynomially with the index, with dyadic blocks instead of doubly exponential ones.

\begin{lemma}[Blocked coupling]\label{lemma:blocked}
Let $(\xi_j)_{j\geq1}$ be independent, centered random variables such that $|\xi_j|\leq C_Mj^{a}$ almost surely, for some $C_M>0$ and $a\in(0,1)$. Then, on a possibly enriched probability space, there exist independent Gaussian variables $(\eta_j)_{j\geq1}$, with $\EE\eta_j=0$ and $\Var(\eta_j)=\Var(\xi_j)$, and a constant $C_1$ depending only on $C_M$, $a$ and the absolute constant $c$ of Lemma~\ref{lemma:sakhanenko}, such that
\begin{equation}\label{eq:blocked}
\PP\Bigg(\exists n\geq 1:\ \Big|\sum_{j=1}^n(\xi_j-\eta_j)\Big|>C_1\,n^{a}\big\{\log(n+1)+u\big\}\Bigg)\leq \exp (-u),\qquad \forall u\ge 0.
\end{equation}
In particular, $\sum_{j\leq n}(\xi_j-\eta_j)= O(n^a\log (n+1))$ almost surely.
\end{lemma}

The approximation error $n^{a}\log n$ in \eqref{eq:blocked} is of the order of the largest summand, up to the logarithm. What matters for the sequel is the \emph{tail control, uniform in $n$}: it allows a union bound over a countable family of such couplings. The almost sure rates of the Strassen type used in \citet[Lemma~A.2]{WASKR2024}, which are well suited to the finite-dimensional problems considered there, do not provide such a tail control, so that they cannot be combined in this way. Lemma~\ref{lemma:blocked} is used for the pointwise results as well.

\paragraph{Gluing.} Lemma~\ref{lemma:blocked} is applied to the centered contributions $(\xi_i(s))_{i\geq1}$ at each point $s$ of the countable grid $G$ defined in Section~\ref{subsubsec:grid}. The couplings for different $s$ live, a priori, on different enriched spaces; the following `gluing lemma' puts them on a single one. Lemma~\ref{lemma:gluing} is a countable version of Lemma A1 of \cite{BP1979}.

\begin{lemma}\label{lemma:gluing}
Let $\mathbf Z=(Z_j)_{j\geq1}$ be a random sequence with values in a Polish space, with distribution $P_{\mathbf Z}$, and let $\mathcal S$ be a countable set. Assume that, for each $s\in\mathcal S$, there exists a probability measure $Q_s$ on the product of this space with $\RR^{\NN}$, with first marginal $P_{\mathbf Z}$. Then there is a probability space carrying $\mathbf Z\sim P_{\mathbf Z}$ and random sequences $\boldsymbol\eta^{(s)}$, $s\in\mathcal S$, such that $(\mathbf Z,\boldsymbol\eta^{(s)})\sim Q_s$ for each $s$.
\end{lemma}

The Gaussian sequences $\boldsymbol\eta^{(s)}$ attached to different points $s$ are conditionally independent given the data in this construction, and they need not be jointly Gaussian. This is the reason why no multivariate strong approximation is needed: the union bound over the grid is applied on the Gaussian side, point by point, where the boundary of the next lemma is nonasymptotic.

\paragraph{Normal mixture boundary.} For $V\geq0$, $\rho>0$ and $\alpha\in(0,1)$, let
\begin{equation}\label{eq:robbins_b}
\mathcal B_\rho(V,\alpha)=\sqrt{\frac{V\rho^2+1}{\rho^2}\log\left(\frac{V\rho^2+1}{\alpha^2}\right)}.
\end{equation}

\begin{lemma}[\citealp{R1970}; \citealp{HRMS2021}]\label{lemma:robbins}
  Let $(\zeta_i)_{i\geq1}$ be independent centered Gaussian variables with variances $\sigma_i^2\geq0$, $T_n=\sum_{i=1}^n\zeta_i$ and $V_n=\sum_{i=1}^n\sigma_i^2$.   Then, for any $\rho>0$ and $\alpha\in(0,1)$,
\[
\PP\big(\exists n\geq 1:\ |T_n|\geq \mathcal B_\rho(V_n,\alpha)\big)\leq \alpha.
\]
\end{lemma}

The boundary $\mathcal B_\rho(V,\alpha)$ is of order $\sqrt{V\log V}$ as $V\uparrow\infty$, wider than a fixed-$n$ Gaussian interval of width $\sqrt V$ by the factor $\sqrt{\log V}$. This factor is due to the mixture with a fixed $\rho$. By the law of the iterated logarithm, time-uniformity alone only costs a factor $\sqrt{\log\log V}$, which stitched boundaries attain \citep{HRMS2021}. These are, however, typically wider near the value of $V$ for which $\rho$ is tuned.   The parameter $\rho$ tunes   the sample size at which the boundary is tightest, and can be chosen as in \citet[Supplement, Appendix~B.2]{WASKR2024}, with a target value of $V$ in place of a target time.
For the bands of Theorem~\ref{theo:band} below, built with a fixed $\rho$, the probability in \eqref{eq:stopping} even tends to one, since the normal mixture boundary  eventually dominates the law of the iterated logarithm.  The bound $1-\alpha$ is thus conservative. For a scalar parameter, a coverage close to $1-\alpha$ for the monitoring from time $m$ on is obtained in \citet[Theorem~2.8]{WASKR2024} with $\rho$ tuned to $m$. For the bands, the union bound over the grid remains a source of conservativeness (item~\ref{rem:band:cov} of Remark~\ref{rem:band}).

\subsection{The calibration grid}\label{subsubsec:grid}

The bands are calibrated on a nested dyadic grid, the \emph{calibration grid}, which does not enter the estimators. The level $\alpha$ is split once among all its points, each grid point $s$ receiving a share $\alpha_s$, as specified in Section~\ref{subsec:general}. Then the band at $t$ uses the variance and the share $\alpha_s$  of the grid point $s=\pi_n(t)$ defined below. The error is  controlled at the grid points, by a union bound valid for all $n$, and in between through Condition~\eqref{cond:M}. For $\ell\geq0$, let 
\begin{equation}\label{eq:grid}
	G_\ell=\big(\{\tau_0 m2^{-\ell}:\ m\in\NN\}\cap\mathcal T\big)\cup\{t_0,t_1\},\qquad G=\bigcup_{\ell\geq0}G_\ell,
\end{equation}
so that $G_0\subset G_1\subset\cdots$ and $G$ is countable. The endpoints are added so that the projection $\pi_n$ below is defined on $\mathcal T$. Let $\ell_0\geq0$ be a fixed initial depth. The \emph{depth} of $s\in G$ is $\ell(s)=\min\{\ell\geq\ell_0: s\in G_\ell\}$, and $N_\ell$ denotes the number of points of depth $\ell$; since $t_1<\tau_0$, $N_\ell\leq2^{\ell+1}$. At time $n$, the grid is $G_{\ell_n}$, where $(\ell_n)$ is a nondecreasing sequence of integers, fixed in advance, with
\begin{equation}\label{eq:ell_n}
	\ell_n\geq\ell_0,\qquad \varsigma_n:=\tau_02^{-\ell_n}\to0,\qquad \ell_n\leq C\log(n+1),
\end{equation}
and $\pi_n(t)=\max\{s\in G_{\ell_n}: s\leq t\}$ is the projection of $t\in\mathcal T$ on $G_{\ell_n}$ by rounding down, so that $0\leq t-\pi_n(t)\leq\varsigma_n$. The grid has to be fine enough for the following condition.

\medskip\noindent
\textbf{Condition (M)}  (interpolation).   With $w_n=w_n(a)$ and $a$ the exponent of Condition~\eqref{cond:I},  
\begin{equation}\tag{M}\label{cond:M}
	\sup_{t\in\mathcal T}|\mathfrak S_n(t)-\mathfrak S_n(\pi_n(t))|=o(nw_n)\ \text{ a.s.}
\end{equation}

\medskip
Condition \eqref{cond:M} is an equicontinuity property of $t\mapsto\mathfrak S_n(t)$ at the scale $\varsigma_n$ of the grid. Together with \eqref{cond:B} and \eqref{cond:I}, it gives $\sup_{t}|\{\hat g_n(t)-g(t)\}-\{\hat g_n(\pi_n(t))-g(\pi_n(t))\}|=o(w_n)$ almost surely when $2\mathfrak b\geq1-a$, which is what the band construction uses. By Lemma~\ref{lemma:chaining}, it holds as soon as the grid is finer than the current bandwidth: $\varsigma_n=o(n^{-\gamma})$ for the kernel and series learners and, for the histogram, whose time-cells at resolution $k$, intersected with $\mathcal T$, have their endpoints in $G_k$, $\ell_n\geq k_n$.

The grid $G_{\ell_n}$ should be distinguished from the set of points where the practitioner evaluates the band, which may be as large as one wishes and has no effect on its validity or width. The grid serves to calibrate the band, through the levels $\alpha_s$ of Section~\ref{subsec:general}, and its size determines the constant of the width. We recommend the grid at the scale of the bandwidth,
\begin{equation}\label{eq:rule_grid}
	\ell_n=\max\big\{\ell_0,\ \lceil\log_2(\kappa_n\tau_0/b(n))\rceil\big\},
\end{equation}
with $\kappa_n\to\infty$ slowly, $\tau_0/J(n)$ replacing $b(n)$ for the series learner. In practice, $\kappa_n$ can be taken constant over the sample sizes of interest; $\kappa_n=3$ is used in   numerical experiments in   Section~\ref{sec:numerics}. A constant $\kappa_n$ is not covered by the theory: Condition \eqref{cond:M} then only holds with $O(nw_n)$. Our bounds then control the interpolation error only at a fraction $O(1/\kappa_n)$ of the half-width. However, one can theoretically take $\kappa_n=\max\{3,\log\log n\}$, and this is equal to $3$ for all $n\leq5\cdot10^{8}$, hence over all sample sizes of Section~\ref{sec:numerics}. The size of this grid   changes little: with $b(n)=n^{-1/5}$, $\tau_0=2$ and $\kappa_n=3$, it has $2^{6}$ cells at $n=10^4$ and $2^{7}$ at $n=10^6$. The fine grid $\ell_n=\max\{\ell_0,\lceil2\log_2(n+1)\rceil\}$ is also admissible, but gives much wider bands without any asymptotic benefit; see item~\ref{rem:band:width} of Remark~\ref{rem:band}.

\subsection{\texorpdfstring{  General results under a one-pass expansion}{General results under a one-pass expansion}}\label{subsec:general} 

Let $(\hat g_n)$ be   an estimator\ of a generic target $g$ on $\mathcal T$,   admitting a one-pass expansion,\ with exponents $a\in(0,1)$ and $\mathfrak b>0$. The results are applied to the hazard function in Section~\ref{subsec:hazard}, through Proposition~\ref{prop:ratio}, and to the cumulative hazard in Section~\ref{subsec:cumhaz}; the components of the at-risk function only need Condition~\eqref{cond:E}.   For $s\in G$ of depth $\ell(s)$, 
the level $\alpha\in(0,1)$ is allocated according to
\begin{equation}\label{eq:alpha_s}
\alpha_s=\alpha\,\varpi_s,\qquad \varpi_s=\frac{1}{N_{\ell(s)}\{\ell(s)-\ell_0+1\}\{\ell(s)-\ell_0+2\}},
\end{equation}
so that  $\sum_{s\in G}\varpi_s\leq\sum_{j\geq0}\{(j+1)(j+2)\}^{-1}=1$, with equality if there are grid points of every depth.   For $s\in G_{\ell_n}$, $\ell(s)\leq\ell_n\leq C\log(n+1)$ and $N_{\ell(s)}\leq2^{\ell(s)+1}$, hence
\begin{equation}\label{eq:log_alpha_s}
\log(1/\alpha_s)\leq\log(1/\alpha)+C\log(n+1),\qquad s\in G_{\ell_n}.
\end{equation}

\begin{remark}[Choice of the weights]\label{rem:weights}
The weights in \eqref{eq:alpha_s} split the level $\alpha$ in two ways. The factor $\{\ell(s)-\ell_0+1\}^{-1}\{\ell(s)-\ell_0+2\}^{-1}$ distributes it across the depths of the grid, which are reached one after the other as $n$ grows. This part is what makes a single union bound valid for all $n$. The factor $1/N_{\ell(s)}$ distributes it across the grid points $s$ of a given depth, that is, across locations $t$ in $\mathcal T$, and is a free choice. The proof of Theorem~\ref{theo:band} only uses $\sum_{s\in G}\varpi_s\leq1$ and the bound \eqref{eq:log_alpha_s}, so that any deterministic weights fixed in advance with these two properties can be used. A union bound over a grid refined with the sample size is used by \citet{HR2022} for quantiles. The allocation \eqref{eq:alpha_s} is a weighted union bound over a nested dyadic family, as in \citet{MH2023}. By \eqref{eq:robbins_b}, the level $\alpha_s$ enters the width through $\log(1/\alpha_s)\leq\log(1/\alpha)+\{\ell(s)+1\}\log2+2\log\{\ell(s)-\ell_0+2\}$. The second term bounds the logarithm of the number of grid points of depth $\ell(s)$, among which the share of this depth is divided. The third bounds the logarithm of the inverse of this share. The allocation \eqref{eq:alpha_s} is economical because the share of a depth is divided among its new points only and  the shares decrease polynomially across depths. Geometric shares, as in \citet{MH2023}, would cost $\ell(s)$ instead of $\log\ell(s)$. Since the grid points of small depth keep the share received when they appeared, the half-width varies irregularly along $t$; see item~\ref{rem:band:irreg} of Remark~\ref{rem:band} below. Dividing the share of each depth among all the points of $G_\ell$ gives a regular half-width, at a cost of about $\log2$ in $\log(1/\alpha_s)$; see Appendix~\ref{app:weights}, which also discusses target-oriented weights.
\end{remark}

\begin{theorem}[Confidence band  and pointwise confidence sequences]\label{theo:band}
Let $\rho>0$ and $\alpha\in(0,1)$.
\begin{enumerate}
\item[(i)]\condlabel{i}{theo:band:i} If $(\hat g_n)$ admits a one-pass expansion \eqref{eq:onepass} satisfying Conditions \eqref{cond:I}, \eqref{cond:V}, \eqref{cond:B} 
and \eqref{cond:M} on $\mathcal T$, with $2\mathfrak b\geq1-a$, then the bands $[\hat g_n(t)\pm\hat L_n(t)]$, $t\in\mathcal T$, $n\geq1$, with
\begin{equation}\label{eq:band_width}
\hat L_n(t)=\frac1n\,\mathcal B_\rho\big(\mathfrak V_n\hat v_n(\pi_n(t)),\alpha_{\pi_n(t)}\big),
\end{equation}
form a $(1\!-\!\alpha)$-asymptotic confidence band sequence for $g$ on $\mathcal T$, and $\sup_{t\in\mathcal T}\hat L_n(t)=O(w_n)$ a.s.
\item[(ii)]\condlabel{ii}{theo:band:ii} Let $t\in\mathcal T$ be fixed. If $(\hat g_n)$   admits a one-pass expansion   \eqref{eq:onepass} and, at the point $t$ only, $|\xi_i(t)|\leq C_Mi^{a}$ a.s. for all $i$, $V_n(t)=\mathfrak V_nv(t)\{1+o(1)\}$ with $v(t)>0$ and $c_0n^{1+a}\leq\mathfrak V_n\leq C_0n^{1+a}$, $\hat v_n(t)\to v(t)$ a.s., $\bar g_n(t)-g(t)=O(n^{-\mathfrak b})$ with $2\mathfrak b\geq1-a$, and $\mathfrak r_n(t)=o(w_n)$ a.s., then
\begin{equation}\label{eq:CS_pointwise}
\Big[\hat g_n(t)\pm\frac1n\,\mathcal B_\rho\big(\mathfrak V_n\hat v_n(t),\alpha\big)\Big],\qquad n\geq1,
\end{equation}
is a $(1\!-\!\alpha)$-asymptotic confidence sequence for $g(t)$, with half-width of order $w_n$ a.s.
\end{enumerate}
\end{theorem}

\begin{remark}\label{rem:band}
\begin{enumerate}
\item\label{rem:band:conv} (Conventions.) The half-widths in Theorem~\ref{theo:band} require $\hat v_n\geq0$, while, for small $n$, $\hat v_n$ can be negative, since the series and boundary-corrected kernels are not nonnegative. We therefore read $\hat v_n$ as $\max(\hat v_n,\mathfrak e_n)$, for a deterministic $\mathfrak e_n\downarrow0$. The same convention applies to the estimator $\hat R_n$ of the at-risk function when necessary, for example in the ratios, variance and integral formulas of Sections~\ref{subsec:hazard} and \ref{subsec:cumhaz}. Since the limits are bounded away from zero, this convention changes nothing for $n$ large, almost surely.

\item\label{rem:band:width} (Width.) The allocation \eqref{eq:alpha_s} replaces $2\log(1/\alpha)$ in the boundary \eqref{eq:robbins_b} by $2\log(1/\alpha_s)$, which is of order $O(\log n)$ by \eqref{eq:log_alpha_s}. This changes the squared boundary only by a constant factor. The half-width thus has the order $w_n$ of the pointwise sequence \eqref{eq:CS_pointwise} and, when $2\mathfrak b\geq1-a$, of the uniform rate in Condition~\eqref{cond:E}: uniformity in $t$ entails no loss in the rate. The constant depends on the grid: if $\ell_n\sim c\log_2n$, then $n^2\hat L_n(t)^2\sim\mathfrak V_n\hat v_n(t)\{1+a+2c\}\log n$ at the locations $t$ calibrated at grid points of depth $\ell_n$. The fine grid has $c=2$ and the grid \eqref{eq:rule_grid} has $c=\gamma$, which reduces the width by the factor $\{(1+a+2\gamma)/(5+a)\}^{1/2}$, that is $(5/16)^{1/2}\approx0.56$ for the kernel and series learners with $d_c=1$, at the rate-optimal $\gamma=1/6$, for which $a=\gamma p=1/3$.

\item\label{rem:band:cov} (Coverage.) The union bound over the grid makes \eqref{eq:def_band} conservative. With the grid \eqref{eq:rule_grid}, its cost $2\gamma\log n$ in the squared width matches, for the hazard, the leading constant suggested by extreme value theory for the supremum of a kernel-type process \citep{BR1973}. The remaining part $(1+a)\log n$ is due to the normal mixture boundary; see Section~\ref{subsec:defs}. Exact coverage would require a joint Gaussian approximation of $t\mapsto\mathfrak S_n(t)$, uniform in $n$, and a time-uniform boundary for its supremum; see Appendix~\ref{app:refinement}.

\item\label{rem:band:bias} (Bias.) Since the width carries the factor $\sqrt{\log n}$, the condition $2\mathfrak b\geq1-a$ allows a bias of order $n^{-(1-a)/2}$: for the three learners of Section~\ref{subsec:learners}, the rate-optimal $\gamma$ is admissible, and the band shrinks at the optimal uniform rate, up to a logarithmic factor; see Corollary~\ref{cor:exponents} and Remark~\ref{rem:exponents}.   In finite samples, the bias is then smaller than the width only by this logarithmic factor, and the proof of Theorem~\ref{theo:band} only gives the rate $1/\sqrt{\log n}$ for the ratios in \eqref{eq:def_band_ratio}. A slight undersmoothing is advisable, as the stress design of Appendix~\ref{app:stress} confirms.

\item\label{rem:band:var} (Variance at $t$.) Using $\hat v_n(t)$ instead of $\hat v_n(\pi_n(t))$ in \eqref{eq:band_width}, with the same level $\alpha_{\pi_n(t)}$, gives an asymptotically equivalent band, since $v$ is uniformly continuous and $\varsigma_n\to0$.

\item\label{rem:band:irreg} (Irregular half-width.) Since $\alpha_{\pi_n(t)}$ depends on the depth of the grid point $\pi_n(t)$, the half-width is larger at the locations $t$ calibrated at grid points of depth $\ell_n$ than at the others, and it varies irregularly along $t$, at the scale $\varsigma_n$. A single level at all locations $t$, not larger than $\min_{s\in G_{\ell_n}}\alpha_s$, such as $\alpha/\{2^{\ell_n+1}(\ell_n+1)(\ell_n+2)\}$, gives a band with a regular half-width, which contains the previous one, has the same order of width and satisfies Theorem~\ref{theo:band}, with the same proof; it is used for the fine grid in Section~\ref{sec:numerics}. See also Appendix~\ref{app:weights}.
	\end{enumerate}
\end{remark}
 
The hazard function is a ratio, and the cumulative hazard an integral of a ratio. The next result transfers the conditions of Sections~\ref{subsec:conditions} and \ref{subsubsec:grid} from   the numerator\ to the ratio, when the denominator is estimated at a faster rate. 

\begin{proposition}[Ratio]\label{prop:ratio}
Assume that $(\hat g_n)$ satisfies Conditions \eqref{cond:I}, \eqref{cond:V}, \eqref{cond:B}, \eqref{cond:E} and \eqref{cond:M} on $\mathcal T$ for $g$, with exponents $a,\mathfrak b$ such that $2\mathfrak b\geq1-a$. Let $R$ be a Lipschitz function on $\mathcal T$ with $\inf_{\mathcal T}R>0$, and let $(\hat R_n)$ be estimators such that $\sup_{t\in\mathcal T}|\hat R_n(t)-R(t)|=o(w_n)$   a.s.   Then $(\hat g_n/\hat R_n)$ satisfies Conditions \eqref{cond:I}, \eqref{cond:V}, \eqref{cond:B}, \eqref{cond:E} and \eqref{cond:M} on $\mathcal T$ for $g/R$, with the same exponents $a,\mathfrak b$, the same sequence $\mathfrak V_n$, and
\[
v^{(R)}(t)=\frac{v(t)}{R(t)^2},\qquad \hat v^{(R)}_n(t)=\frac{\hat v_n(t)}{\hat R_n(t)^2}.
\]
\end{proposition}

\subsection{The conditional hazard function}\label{subsec:hazard}

Let $\lambda(t)=\lambda_{T\mid\X}(t\mid\x)=f(t,1,\x)/R(t,\x)$ on $\mathcal T=[t_0,t_1]$, $t_0>0$, with the at-risk function
\begin{equation}\label{eq:R_affine}
R(t,\x)=R_1(t,\x)+\eta(\x)R_2(t,\x),
\end{equation}
where $\eta(\x)$ is a scalar nuisance parameter, known or unknown, and $R_1,R_2$ are functions of the distribution of the observables; the examples of Section~\ref{sec:mechanisms} have $\eta\equiv0$ (right censoring; left truncation and right censoring, Appendix~\ref{app:ltrc}) or an unknown $\eta$ (modified current status; cure models can also be written in this form, see Section~\ref{subsec:cure}). The estimator is
\begin{equation}\label{eq:hat_lambda}
\hat\lambda_n(t)=\frac{\hat f_n(t)}{\hat R_n(t)},\qquad \hat R_n(t)=\hat R_{1,n}(t)+\hat\eta_n\hat R_{2,n}(t),
\end{equation}
with $\hat f_n$   an estimator\ of $f(\cdot,1,\x)$, $\hat R_{1,n},\hat R_{2,n}$ estimators of $R_1(\cdot,\x),R_2(\cdot,\x)$ and $\hat\eta_n$ an estimator of $\eta(\x)$, with $\hat\eta_n=\eta$ when the nuisance is known.

\begin{assum}\label{ass:hazard}
\begin{enumerate}
\item[(i)]\condlabel{i}{ass:hazard:i} $(\hat f_n)$ satisfies Conditions \eqref{cond:I}, \eqref{cond:V}, \eqref{cond:B}, \eqref{cond:E} and \eqref{cond:M} on $\mathcal T$ for $f(\cdot,1,\x)$, with exponents $a,\mathfrak b$ such that $2\mathfrak b\geq1-a$.
\item[(ii)]\condlabel{ii}{ass:hazard:ii} $R_1(\cdot,\x)$ and $R_2(\cdot,\x)$ are Lipschitz on $\mathcal T$, and $\inf_{t\in\mathcal T}R(t,\x)>0$.
\item[(iii)]\condlabel{iii}{ass:hazard:iii} $\sup_{t\in\mathcal T}|\hat R_{m,n}(t)-R_m(t,\x)|=o(w_n)$    a.s.,   $m=1,2$, and $|\hat\eta_n-\eta(\x)|=o(w_n)$   a.s.  
\end{enumerate}
\end{assum}

Condition (\ref{ass:hazard:iii}) holds, for instance, when $\hat R_{m,n}$ are   one-pass averages   satisfying \eqref{cond:E} with exponents $a_R<a$ and $\mathfrak b_R\geq(1-a)/2$: the stochastic part of the denominator is then of smaller order than that of the numerator, because the denominator only involves smoothing in $\x_c$, in dimension $d_c$, while the numerator involves smoothing in $(t,\x_c)$, in dimension $d_c+1$; see Section~\ref{sec:examples}. The same applies to $\hat\eta_n$ in the examples of Section~\ref{sec:mechanisms}, where $\eta(\x)$ is a functional of the conditional distribution of the observables given $\X=\x$, estimated by smoothing in $\x_c$ only.

\begin{theorem}[Hazard]\label{theo:hazard}
Under Assumption~\ref{ass:hazard}, let
\begin{equation}\label{eq:band_hazard}
\hat L^{\lambda}_n(t)=\frac1n\,\mathcal B_\rho\Big(\mathfrak V_n\frac{\hat v_n(\pi_n(t))}{\hat R_n(\pi_n(t))^2},\alpha_{\pi_n(t)}\Big),\qquad t\in\mathcal T.
\end{equation}
Then, for any $\rho>0$ and $\alpha\in(0,1)$, the bands $[\hat\lambda_n(t)\pm\hat L^{\lambda}_n(t)]$, $t\in\mathcal T$, $n\geq1$, form a $(1-\alpha)$-asymptotic confidence band sequence for $\lambda_{T\mid\X}(\cdot\mid\x)$ on $\mathcal T$, with $\sup_t\hat L^\lambda_n(t)=O(w_n)$   a.s.   For a fixed $t\in\mathcal T$, $[\hat\lambda_n(t)\pm n^{-1}\mathcal B_\rho(\mathfrak V_n\hat v_n(t)/\hat R_n(t)^2,\alpha)]$ is a $(1-\alpha)$-asymptotic confidence sequence for $\lambda_{T\mid\X}(t\mid\x)$.
\end{theorem}

For the three learners of Section~\ref{subsec:learners}, $\hat v_n=\mathfrak c\hat f_n$ for a known constant $\mathfrak c$, so that $\mathfrak V_n\hat v_n/\hat R_n^2=\mathfrak c\,\mathfrak V_n\hat\lambda_n/\hat R_n$: the variance of the hazard estimator is, as in the fixed-$n$ theory, proportional to the hazard divided by the at-risk function.

\begin{remark}[Local-maximum variance]\label{rem:localmax}
With the plug-in $\hat v_n=\mathfrak c\hat f_n$, the band is narrow where $\hat f_n$ is small, that is exactly where the estimate is below the target, and this happens in the regions of low density when $n$ is small. Condition \eqref{cond:V} leaves room for a better choice. Let $b_n\to0$ be the current resolution of the learner, $b(n)$ for the kernel learner, $\tau_0/J(n)$ for the series learner and the length of the time-cells of $\mathcal C_{k_n}$ for the histogram, and let
\[
\hat v^{\max}_n(t)=\mathfrak c\max_{s\in\mathcal T,\,|s-t|\leq b_n}\hat f_n(s).
\]
Since $f(\cdot,1,\x)$ is uniformly continuous and bounded away from zero on $\mathcal T$, Condition \eqref{cond:E} gives $\sup_{\mathcal T}|\hat v^{\max}_n/v-1|\to0$ almost surely. A window extending beyond $\mathcal T$ can be used if Condition \eqref{cond:E} holds on a neighborhood of $\mathcal T$. The local-maximum version therefore satisfies Condition \eqref{cond:V}, all the results hold with it, and it requires no additional storage. The same applies to $\hat v^\Lambda_n$ below, with the local-maximum estimator in the integral. We recommend it: its effect in finite samples is substantial; see Section~\ref{sec:numerics}.
\end{remark}

\subsection{The conditional cumulative hazard and the survival function}\label{subsec:cumhaz}

Let $\Lambda(t)=\Lambda_{T\mid\X}(t\mid\x)=\int_0^t\lambda(s)\dd s$ and $\hat\Lambda_n(t)=\int_0^t\hat\lambda_n(s)\dd s$, where now $\hat f_n$ and $\hat R_n$ are defined on $[0,t_1]$. The integration in time changes the structure of the problem in two ways. First, it removes the smoothing in $t$: the contributions
\begin{equation}\label{eq:psi_Lambda}
\psi^{\Lambda}_i(t;\mathfrak z)=\int_0^t\frac{\psi_i(s;\mathfrak z)}{R(s,\x)}\dd s
\end{equation}
of the   one-pass expansion of $\hat\Lambda_n$, whose leading term   $\tilde\Lambda_n(t)=n^{-1}\sum_{i\leq n}\psi_i^\Lambda(t;Z_i)$ uses the true at-risk function, are of a smaller order than the $\psi_i$, with a growth exponent $a_\Lambda<a$, so that the bands for $\Lambda$ are narrower than those for $\lambda$, of order $w_n(a_\Lambda)$. Second, and as a consequence, the error of the denominator has to be compared with $w_n(a_\Lambda)$ instead of $w_n(a)$: the denominator has to be estimated with a smaller stochastic error, which, for kernel-type learners, means a larger bandwidth than the one used for the numerator. The band is constructed on $\mathcal T^\Lambda=[t_0,t_1]$ with $t_0>0$, where the variance of $\tilde\Lambda_n$ is bounded away from zero, while the learners are used on $[0,t_1]$.

\begin{assum}\label{ass:cumhaz} 
\begin{enumerate}
\item[(i)]\condlabel{i}{ass:cumhaz:i} $(\hat f_n)$   admits a one-pass expansion   \eqref{eq:onepass} and satisfies Conditions \eqref{cond:I}, \eqref{cond:B} and \eqref{cond:E} on $[0,t_1]$ for $f(\cdot,1,\x)$, with exponents $a,\mathfrak b$.
\item[(ii)]\condlabel{ii}{ass:cumhaz:ii} $R_1(\cdot,\x)$ and $R_2(\cdot,\x)$ are Lipschitz on $[0,t_1]$, and $\inf_{t\in[0,t_1]}R(t,\x)>0$.
\item[(iii)]\condlabel{iii}{ass:cumhaz:iii} (Integrated contributions) There exist   $a_\Lambda\in(0,a)$   and constants $C_M',C_V'$ with
\[
\sup_{t\in[0,t_1]}|\xi^\Lambda_i(t)|\leq C_M'i^{a_\Lambda}\ \text{    a.s.,  }\qquad \sup_{t\in[0,t_1]}\Var\{\xi_i^\Lambda(t)\}\leq C_V'i^{a_\Lambda},
\]
where $\xi_i^\Lambda(t)=\psi_i^\Lambda(t;Z_i)-\EE\psi^\Lambda_i(t;Z_i)$.
Moreover, Condition \eqref{cond:V} holds on $\mathcal T^\Lambda$ for the contributions $\psi_i^\Lambda$, with a known sequence $\mathfrak V^\Lambda_n$, $c_0n^{1+a_\Lambda}\leq\mathfrak V_n^\Lambda\leq C_0n^{1+a_\Lambda}$, a continuous $v^\Lambda>0$ and computable $\hat v^\Lambda_n$ with $\sup_{\mathcal T^\Lambda}|\hat v^\Lambda_n/v^\Lambda-1|\to0$ almost surely. Finally,
\[
\sup_{t\in[0,t_1]}\Big|\int_0^t\frac{\mathfrak r_n(s)}{R(s,\x)}\dd s\Big|=o(w_n(a_\Lambda))\quad\text{   a.s.  }
\]
\item[(iv)]\condlabel{iv}{ass:cumhaz:iv} $\sup_{t\in[0,t_1]}|\hat R_{m,n}(t)-R_m(t,\x)|=o(w_n(a_\Lambda))$    a.s.,   $m=1,2$, and $|\hat\eta_n-\eta(\x)|=o(w_n(a_\Lambda))$    a.s.  
\item[(v)]\condlabel{v}{ass:cumhaz:v} $2\mathfrak b\geq1-a_\Lambda$ and $\varsigma_nn^{(a-a_\Lambda)/2}\to0$.
\end{enumerate}
\end{assum}

Condition (\ref{ass:cumhaz:iii}) is a property of the learner, verified in Section~\ref{sec:examples} for the three learners of Section~\ref{subsec:learners}, with $a_\Lambda=\gamma d_c$ and $v^\Lambda(t)=\mathfrak c_\Lambda\int_0^tf(s,1,\x)R(s,\x)^{-2}\dd s$, for a constant $\mathfrak c_\Lambda$ depending on the learner and given in Section~\ref{subsec:numerator}: the integrated contributions behave like the contributions of a one-pass estimator of a regression-type quantity, smoothed in $\x_c$ only. The last requirement of (\ref{ass:cumhaz:iii}) concerns the integrated remainder, the only part of $\mathfrak r_n$ compared with the smaller width $w_n(a_\Lambda)$. It is trivial for the three learners, where $\mathfrak r_n\equiv0$, and a nonzero remainder is allowed: for a remainder of the quadratic type of Remark~\ref{rem:remainder}, of order $n^{-2\mathfrak b}+w_n(a)^2$, it holds as soon as $2(1-a)>1-a_\Lambda$. Note that Condition \eqref{cond:V} is not required for $\hat f_n$ itself, so that the boundary behavior of the variance of the learners at $t=0$ plays no role; see Remark~\ref{rem:boundary}. Condition (\ref{ass:cumhaz:v}) gathers the restrictions on the exponents and on the grid; the second one holds for the grids of Section~\ref{subsubsec:grid}, where $a-a_\Lambda=\gamma$. It can be relaxed to $\varsigma_n\to0$, which allows narrower bands for $\Lambda$ and $S$; see Appendix~\ref{app:grid_cumhaz}. Condition (\ref{ass:cumhaz:iv}) is the counterpart of Assumption~\ref{ass:hazard}(\ref{ass:hazard:iii}) with the smaller rate $w_n(a_\Lambda)$; for one-pass denominators with exponents $a_R,\mathfrak b_R$, it holds as soon as $a_R<a_\Lambda$ and $\mathfrak b_R\geq(1-a_\Lambda)/2$.

\begin{theorem}[Cumulative hazard]\label{theo:cumhaz}
Under Assumption~\ref{ass:cumhaz}, let
\[
\hat L^{\Lambda}_n(t)=\frac1n\,\mathcal B_\rho\big(\mathfrak V^\Lambda_n\hat v^\Lambda_n(\pi_n(t)),\alpha_{\pi_n(t)}\big),\qquad t\in\mathcal T^\Lambda.
\]
Then, for any $\rho>0$ and $\alpha\in(0,1)$, the bands $[\hat\Lambda_n(t)\pm\hat L^{\Lambda}_n(t)]$, $t\in\mathcal T^\Lambda$, $n\geq1$, form a $(1-\alpha)$-asymptotic confidence band sequence for $\Lambda_{T\mid\X}(\cdot\mid\x)$ on $\mathcal T^\Lambda$, with $\sup_t\hat L^\Lambda_n(t)=O(w_n(a_\Lambda))$   a.s.   For a fixed $t\in\mathcal T^\Lambda$, $[\hat\Lambda_n(t)\pm n^{-1}\mathcal B_\rho(\mathfrak V^\Lambda_n\hat v^\Lambda_n(t),\alpha)]$ is a $(1-\alpha)$-asymptotic confidence sequence for $\Lambda_{T\mid\X}(t\mid\x)$.
\end{theorem}

\begin{corollary}[Survival function]\label{cor:survival}
Under Assumption~\ref{ass:cumhaz}, let $\hat S_n(t)=\exp\{-\hat\Lambda_n(t)\}$. Then the bands
\[
\big[\hat S_n(t)\,e^{-\hat L^\Lambda_n(t)},\ \hat S_n(t)\,e^{\hat L^\Lambda_n(t)}\big],\qquad t\in\mathcal T^\Lambda,\ n\geq1,
\]
form a $(1-\alpha)$-asymptotic confidence band sequence for $S_{T\mid\X}(\cdot\mid\x)$ on $\mathcal T^\Lambda$.
   The pointwise statement holds likewise.  
\end{corollary}

The bands we propose can be adapted to the intrinsic nature of the targets, such as range and shape constraints, at no cost in coverage. See Appendix~\ref{subsec:app_range_shape}.   

\begin{remark}[Nuisance parameters estimated at a single point]\label{rem:nuisance_slow}
In the cure models of Section~\ref{subsec:cure}, $-\eta(\x)$ is the cure proportion, estimated by the survival estimator at the end $\tau$ of the follow-up. Under the conditions on $(\gamma,\beta,\nu)$ of Corollary~\ref{cor:exponents}(\ref{cor:exponents:ii}), its error is $o(w_n(a_\Lambda))$ almost surely. Indeed, Lemma~\ref{lemma:blocked} at the single point $\tau$ and the law of the iterated logarithm give a stochastic part of order $n^{-(1-a_\Lambda)/2}\sqrt{\log\log n}$, the bias is $O(n^{-(1-a_\Lambda)/2})$, and the error of the denominator is $o(w_n(a_\Lambda))$. The requirement on $\hat\eta_n$ in Assumption~\ref{ass:cumhaz}(\ref{ass:cumhaz:iv}) thus holds, although this error is of the order of a standard deviation: the factor $\sqrt{\log n}$ in $w_n(a_\Lambda)$ absorbs it. A band ignoring it would be asymptotically valid, but with an incorrect variance profile. Theorem~\ref{theo:cure}(\ref{theo:cure:ii}) includes this error in the contributions, and gives the correct variance profile.
\end{remark}


\section{The three learners: verification of the conditions}\label{sec:examples} 

We now verify Conditions \eqref{cond:I}, \eqref{cond:V}, \eqref{cond:B}, \eqref{cond:E} and \eqref{cond:M} of Sections~\ref{subsec:conditions} and \ref{subsubsec:grid}, and the integrated conditions of Assumption~\ref{ass:cumhaz}(\ref{ass:cumhaz:iii}), for the three learners of Section~\ref{subsec:learners}, and we give one-pass kernel estimators of the components of the at-risk function satisfying the rate conditions of Assumptions~\ref{ass:hazard}(\ref{ass:hazard:iii}) and \ref{ass:cumhaz}(\ref{ass:cumhaz:iv}). This gives explicit conditions on the exponents $\gamma$ and $\beta$ under which Theorems~\ref{theo:hazard} and \ref{theo:cumhaz} apply. The proofs are in Appendix~\ref{app:examples}. Throughout, $d_c\geq1$; see Appendix~\ref{app:dc0} for the case $d_c=0$.

\subsection{Assumptions on the distribution of the observables}\label{subsec:data_ass}

  Recall that $\x=(\x_c,\x_d)$ denotes the fixed covariate value used in conditioning, $\z=(\z_c,\z_d)$ denotes a generic value in the space of values of the covariate vector $\X = (\X_c,\X_d)$,  
$\mathcal T=[t_0,t_1]$ with $0\leq t_0<t_1<\tau_0$, and $f(y,1,\z)=f_{Y,I,\X}(y,1,\z)$ is the joint density of $(Y,\X)$ on the event $\{I=1\}$, with respect to the product of the Lebesgue measure on $[0,\infty)\times\mathcal X_c$ and the counting measure on $\mathcal X_d$. Let $\mathcal N$ be a closed ball centered at $\x_c$ and contained in the interior of $\mathcal X_c$.

\begin{assum}\label{ass:data}
\begin{enumerate}
\item[(D1)]\condlabel{D1}{ass:D1} $\mathcal X_d$ is finite, $\PP(\X_d=\x_d)>0$, and $\x_c$ is an interior point of the compact rectangle $\mathcal X_c\subset\RR^{d_c}$.
\item[(D2)]\condlabel{D2}{ass:D2} The function $(y,\z_c)\mapsto f(y,1,(\z_c,\x_d))$ is bounded on $[0,\tau_0]\times\mathcal N$ and satisfies one of the following smoothness conditions, according to the learner:
\begin{enumerate}
\item[(D2a)]\condlabel{D2a}{ass:D2a} (kernel) it has partial derivatives of order two which are bounded and uniformly continuous on $[0,\tau_0]\times\mathcal N$;
\item[(D2b)]\condlabel{D2b}{ass:D2b} (series) it has partial derivatives of order two which are bounded and Lipschitz on $[0,\tau_0]\times\mathcal N$;
\item[(D2c)]\condlabel{D2c}{ass:D2c} (histogram) it is Lipschitz on $[0,\tau_0]\times\mathcal N$.
\end{enumerate}
\item[(D3)]\condlabel{D3}{ass:D3} $\inf_{t\in\mathcal T}f(t,1,\x)>0$ for the hazard function, and $\int_0^{t_0}f(s,1,\x)\dd s>0$ with $t_0>0$ for the cumulative hazard.
\item[(D4)]\condlabel{D4}{ass:D4} (At-risk components) For $m=1,2$, $R_m(t,\x)=\EE[u_m(t;Z)\mid\X=\x]f_\X(\x)$, where $f_\X$ is the density of $\X$ with respect to the product of the Lebesgue and counting measures, $f_\X(\cdot,\x_d)$ is bounded on $\mathcal N$, and $u_m=u_m^+-u_m^-$ for known functions $u_m^\pm:[0,\tau_0]\times\mathcal Z\to\RR$, bounded, jointly measurable and non-decreasing in $t$, such that $t\mapsto\EE[u_m^{\pm}(t;Z)\mid\X=(\z_c,\x_d)]f_\X(\z_c,\x_d)$ is Lipschitz on $[0,\tau_0]$, uniformly in $\z_c\in\mathcal N$; and $\z_c\mapsto R_m(t,(\z_c,\x_d))$ has partial derivatives of order two bounded on $\mathcal N$, uniformly in $t\in[0,\tau_0]$. Moreover, $R(\cdot,\x)=R_1(\cdot,\x)+\eta(\x)R_2(\cdot,\x)$ is Lipschitz on $[0,\tau_0]$ and $\inf_{t\in[0,t_1]}R(t,\x)>0$.
\end{enumerate}
\end{assum}

Condition (\ref{ass:D2}) concerns the sub-density of the exactly observed durations only. Nothing is assumed directly on the distribution of the censoring or truncation variables, but since $f(\cdot,1,\x)$ involves them, for instance $f(t,1,\x)=S_{C\mid\X}(t-\mid\x)f_{T\mid\X}(t\mid\x)f_\X(\x)$ under right censoring, the smoothness in (\ref{ass:D2}) is a joint condition on the duration of interest and on the incompleteness mechanism. For instance, administrative censoring at a time $c^*$ creates a jump of $f(\cdot,1,\x)$ at $c^*$, so that $\tau_0<c^*$ is needed. This matters most for the series learner: a jump inside $(0,\tau_0)$   degrades   its approximation on the whole interval, while for the kernel and histogram learners it only affects a neighborhood of the jump. The requirements in (\ref{ass:D4}) are checked for each mechanism in Section~\ref{sec:mechanisms} and Appendix~\ref{app:mechanisms}.
In these examples, $u_m(t;\mathfrak z)$ is an indicator such as $\1\{y\geq t\}$ or $\1\{\ell\leq t\leq y\}$, or a difference of such indicators, so that the decomposition in (\ref{ass:D4}) is explicit and its Lipschitz requirement reduces to the existence of bounded conditional densities for $Y$, $L$ or $C$ given $\X$. Condition (\ref{ass:D3}) ensures that the variance profile is bounded away from zero on $\mathcal T$, which is where the band is constructed; the series learner also needs $t_0>0$ for the hazard, see Section~\ref{subsubsec:series}. The smoothness in $y$ in (\ref{ass:D2}) is required on $[0,\tau_0]$, including at the origin, because the cumulative hazard integrates from $0$; the constant $\tau_0$ is chosen by the user, and smoothness is required only up to $\tau_0$. When the at-risk function is not bounded away from zero near the origin, under truncation, the interval starts at a known $s_0$ instead of $0$ in all these conditions; see Remark~\ref{rem:s0}.

\subsection{The numerator}\label{subsec:numerator}

We write $\mathfrak c_K=c_K^{p}c_b^{-1}c_h^{-d_c}$ and $\mathfrak c_\Lambda=c_K^{d_c}c_h^{-d_c}$, where $c_K=\int K^2$.

\begin{proposition}[Recursive kernel estimator]\label{prop:kernel}
Let Assumption~\ref{ass:data} (\ref{ass:D1}), (\ref{ass:D2a}), (\ref{ass:D3}) and (\ref{ass:D4}) hold, let $K$ be a symmetric, Lipschitz continuous probability density supported on $[-1,1]$, let $\gamma\in(0,1/p)$ and $\varsigma_n=o(n^{-\gamma})$.
\begin{enumerate}
\item[(i)]\condlabel{i}{prop:kernel:i} If $t_0>0$, the estimator $\hat f^{\mathrm K}_n$ of Section~\ref{subsubsec:kernel} satisfies Conditions \eqref{cond:I}, \eqref{cond:V}, \eqref{cond:B}, \eqref{cond:E} and \eqref{cond:M} on $\mathcal T$ for $f(\cdot,1,\x)$, with
\[
a=\gamma p,\quad \mathfrak b=2\gamma,\quad \mathfrak V_n=\frac{n^{1+\gamma p}}{1+\gamma p},\quad v(t)=\mathfrak c_Kf(t,1,\x),\quad \hat v_n(t)=\mathfrak c_K\hat f^{\mathrm K}_n(t,1,\x).
\]
\item[(ii)]\condlabel{ii}{prop:kernel:ii} Let the boundary-corrected contributions be given by \eqref{eq:psi_kernel} with $K_{b(i)}(y-t)$ replaced by $K^{t}_{B,b(i)}(y-t)$, where, for $t\in[0,b)$, $K^t_{B,b}(s)=\{\mathsf k_0(t/b)-\mathsf k_1(t/b)s/b\}K_b(s)$ and $(\mathsf k_0(\theta),\mathsf k_1(\theta))$ solves $\mathsf k_0\int_{-1}^{\theta}K+\mathsf k_1\int_{-1}^\theta uK=1$, $\mathsf k_0\int_{-1}^\theta uK+\mathsf k_1\int_{-1}^\theta u^2K=0$, while $K^t_{B,b}=K_b$ for $t\geq b$. Then $\hat f^{\mathrm K}_n$ satisfies Conditions \eqref{cond:I}, \eqref{cond:B} and \eqref{cond:E} on $[0,t_1]$ with the same $a$ and $\mathfrak b$, and the integrated contributions \eqref{eq:psi_Lambda} satisfy Assumption~\ref{ass:cumhaz}(\ref{ass:cumhaz:iii}) with
\begin{align*}
a_\Lambda&=\gamma d_c,\qquad \mathfrak V^\Lambda_n=\frac{n^{1+\gamma d_c}}{1+\gamma d_c},\\
v^\Lambda(t)&=\mathfrak c_\Lambda\int_0^t\frac{f(s,1,\x)}{R(s,\x)^2}\dd s,\qquad \hat v^\Lambda_n(t)=\mathfrak c_\Lambda\int_0^t\frac{\hat f^{\mathrm K}_n(s,1,\x)}{\hat R_n(s)^2}\dd s,
\end{align*}
for any estimator $\hat R_n$ with $\sup_{[0,t_1]}|\hat R_n-R(\cdot,\x)|\to0$ almost surely.
\end{enumerate}
\end{proposition}

The boundary kernel in (\ref{prop:kernel:ii}) is the linear boundary kernel discussed by \cite{J1993}; see also \cite{GMM1985}. The linear system defining $(\mathsf k_0,\mathsf k_1)$ is nonsingular for every $\theta\in[0,1]$ by the Cauchy--Schwarz inequality, and $\mathsf k_0,\mathsf k_1$ are bounded on $[0,1]$. The only role of the correction is to keep the bias of order $b(i)^2$ on $[0,b(i))$; it does not change the orders of the contributions. The variance profile $v^\Lambda$ no longer involves the time bandwidth: integration over $t$ removes the smoothing in time, and the integrated contributions behave like those of a kernel estimator of a regression function in $\x_c$, with the growth exponent $\gamma d_c$ instead of $\gamma(d_c+1)$.

\begin{proposition}[Recursive series estimator]\label{prop:series}
Let Assumption~\ref{ass:data} (\ref{ass:D1}), (\ref{ass:D2b}), (\ref{ass:D3}) and (\ref{ass:D4}) hold, let $\{\varphi_j\}$ be the augmented half-cosine basis of Section~\ref{subsubsec:series}, let $K$ be as in Proposition~\ref{prop:kernel}, $\gamma\in(0,1/p)$ and $\varsigma_n=o(n^{-\gamma})$.
\begin{enumerate}
\item[(i)]\condlabel{i}{prop:series:i} If $t_0>0$, the estimator $\hat f^{\mathrm S}_n$ satisfies Conditions \eqref{cond:I}, \eqref{cond:V}, \eqref{cond:B}, \eqref{cond:E} and \eqref{cond:M} on $\mathcal T$ for $f(\cdot,1,\x)$, with
\[
a=\gamma p,\quad \mathfrak b=2\gamma,\quad \mathfrak V_n=\sum_{i=1}^nJ(i)\,i^{\gamma d_c},\quad v(t)=\frac{c_K^{d_c}f(t,1,\x)}{\tau_0c_h^{d_c}},\quad \hat v_n(t)=\frac{c_K^{d_c}\hat f^{\mathrm S}_n(t,1,\x)}{\tau_0c_h^{d_c}}.
\]
\item[(ii)]\condlabel{ii}{prop:series:ii} Without any modification, $\hat f^{\mathrm S}_n$ satisfies Conditions \eqref{cond:I}, \eqref{cond:B} and \eqref{cond:E} on $[0,t_1]$, and the integrated contributions \eqref{eq:psi_Lambda} satisfy Assumption~\ref{ass:cumhaz}(\ref{ass:cumhaz:iii}) with $a_\Lambda=\gamma d_c$ and the same $\mathfrak V^\Lambda_n$, $v^\Lambda$ and $\hat v^\Lambda_n$ as in Proposition~\ref{prop:kernel}(\ref{prop:kernel:ii}), with $\hat f^{\mathrm S}_n$ in place of $\hat f^{\mathrm K}_n$.
\end{enumerate}
\end{proposition}

Two features of the series learner deserve a comment. First, the variance profile of the hazard estimator involves the basis only through the constant $1/\tau_0$, the counterpart of $c_K/c_b$ for the kernel learner; this is the localization property \eqref{cond:S3}. Second, the variance profile of the integrated contributions does not involve the basis at all: once integrated in time, the series and kernel learners have the same asymptotic variance, as the smoothing in $t$ disappears in both cases. The proof uses the uniform boundedness of the partial sums of the series expansion of functions of bounded variation (Lemma~\ref{lem:S3}(\ref{lem:S3:ii})), which is where the structure of the augmented basis enters.

\begin{proposition}[Dyadic streaming histogram]\label{prop:histogram}
Let Assumption~\ref{ass:data} (\ref{ass:D1}), (\ref{ass:D2c}), (\ref{ass:D3}) and (\ref{ass:D4}) hold, with $\mathcal X_c\subset[0,1)^{d_c}$, let $\gamma\in(0,1/p)$ and $\ell_n\geq k_n$.
\begin{enumerate}
\item[(i)]\condlabel{i}{prop:histogram:i} The estimator $\hat f^{\mathrm H}_n$ of Section~\ref{subsubsec:histogram} satisfies Conditions \eqref{cond:I}, \eqref{cond:V}, \eqref{cond:B}, \eqref{cond:E} and \eqref{cond:M} on $\mathcal T$ for $f(\cdot,1,\x)$, including when $t_0=0$, with
\[
a=\gamma p,\qquad \mathfrak b=\gamma,\qquad \mathfrak V_n=\tau_0^{-1}\sum_{i=1}^n2^{k_ip},\qquad v(t)=f(t,1,\x),\qquad \hat v_n(t)=\hat f^{\mathrm H}_n(t,1,\x).
\]
\item[(ii)]\condlabel{ii}{prop:histogram:ii} The integrated contributions \eqref{eq:psi_Lambda} satisfy Assumption~\ref{ass:cumhaz}(\ref{ass:cumhaz:iii}) with
\[
a_\Lambda=\gamma d_c,\qquad \mathfrak V^\Lambda_n=\sum_{i=1}^n2^{k_id_c},\qquad v^\Lambda(t)=\int_0^t\frac{f(s,1,\x)}{R(s,\x)^2}\dd s,\qquad \hat v^\Lambda_n(t)=\int_0^t\frac{\hat f^{\mathrm H}_n(s,1,\x)}{\hat R_n(s)^2}\dd s.
\]
for any estimator $\hat R_n$ with $\sup_{[0,t_1]}|\hat R_n-R(\cdot,\x)|\to0$ almost surely.
\end{enumerate}
\end{proposition}

\subsection{The denominator}\label{subsec:denominator}

The components $R_m(\cdot,\x)$ of the at-risk function are estimated by one-pass kernel estimators smoothing in $\x_c$ only: with $g(i)=c_gi^{-\beta}\{\log(i+1)\}^{\nu}$, $\beta\in(0,1/d_c)$, $\nu\in[0,1/4)$, $c_g>0$,
\begin{equation}\label{eq:psi_R}
\psi^{R}_{m,i}(t;\mathfrak z)=u_m(t;\mathfrak z)\,\boldsymbol K_{g(i)}(\z_c-\x_c)\,\1\{\z_d=\x_d\},\qquad \hat R_{m,n}(t)=\frac1n\sum_{i=1}^n\psi^{R}_{m,i}(t;Z_i).
\end{equation}
The logarithmic factor only matters at the end of the range of exponents for the cumulative hazard, see Corollary~\ref{cor:exponents}(\ref{cor:exponents:ii}); otherwise, one can take $\nu=0$.

\begin{proposition}[At-risk components]\label{prop:denominator}
Under Assumption~\ref{ass:data} (\ref{ass:D1}) and (\ref{ass:D4}), with $K$ as in Proposition~\ref{prop:kernel}, the estimators \eqref{eq:psi_R} satisfy, for $m=1,2$,
\begin{equation}\label{eq:rate_R}
\sup_{t\in[0,\tau_0]}|\hat R_{m,n}(t)-R_m(t,\x)|=O(r^R_n)\quad\text{almost surely},
\end{equation}
with
\[
r^R_n=\frac{\{\log(n+1)\}^{2\nu}}{n^{2\beta}}+\frac{\{\log(n+1)\}^{(1-\nu d_c)/2}}{n^{(1-\beta d_c)/2}}.
\]
For $\nu=0$ and $\beta<1/2$, they satisfy Conditions \eqref{cond:I}, \eqref{cond:B} and \eqref{cond:E} on $[0,\tau_0]$ for $R_m(\cdot,\x)$, with $a_R=\beta d_c$ and $\mathfrak b_R=2\beta$.
\end{proposition}

Combining Propositions~\ref{prop:kernel}--\ref{prop:denominator} with Theorems~\ref{theo:hazard} and \ref{theo:cumhaz} gives explicit conditions on the exponents. In the next statement, the nuisance parameter is assumed known or estimated at the required rate; the estimators of $\eta(\x)$ are discussed mechanism by mechanism in Section~\ref{sec:mechanisms}.

\begin{corollary}[Exponents]\label{cor:exponents}
Let Assumption~\ref{ass:data} hold with the smoothness condition (\ref{ass:D2}) corresponding to the learner, let the grid be as in Propositions~\ref{prop:kernel}--\ref{prop:histogram}, and let $\hat R_{m,n}$ be given by \eqref{eq:psi_R}.
\begin{enumerate}
\item[(i)]\condlabel{i}{cor:exponents:i} (Hazard function, kernel and series learners) Assumption~\ref{ass:hazard} holds, and the bands of Theorem~\ref{theo:hazard} are valid on $\mathcal T=[t_0,t_1]$, $t_0>0$, provided $|\hat\eta_n-\eta(\x)|=o(w_n(\gamma p))$ almost surely and
\[
\frac{1}{d_c+5}\leq\gamma<\frac{1}{d_c+1},\qquad \frac{1-\gamma p}{4}\leq\beta<\frac{\gamma p}{d_c}.
\]
The width of the bands is of order $n^{-(1-\gamma p)/2}\sqrt{\log n}$, which equals $n^{-2/(d_c+5)}\sqrt{\log n}$ for $\gamma=1/(d_c+5)$.
\item[(ii)]\condlabel{ii}{cor:exponents:ii} (Cumulative hazard and survival function, kernel with boundary correction and series learners) Assumption~\ref{ass:cumhaz} holds, and the bands of Theorem~\ref{theo:cumhaz} and Corollary~\ref{cor:survival} are valid on $\mathcal T^\Lambda=[t_0,t_1]$, $t_0>0$, provided $|\hat\eta_n-\eta(\x)|=o(w_n(\gamma d_c))$ almost surely and
\[
\frac{1}{d_c+4}\leq\gamma<\frac{1}{d_c+1},\qquad \frac{1-\gamma d_c}{4}\leq\beta\leq\gamma,\qquad \nu>0\ \text{ if }\beta=\gamma.
\]
The width of the bands is of order $n^{-(1-\gamma d_c)/2}\sqrt{\log n}$, which equals $n^{-2/(d_c+4)}\sqrt{\log n}$ for $\gamma=1/(d_c+4)$, where these sufficient conditions force $\beta=\gamma$ and $\nu\in(0,1/4)$.
\item[(iii)]\condlabel{iii}{cor:exponents:iii} (Histogram learner) The statements (\ref{cor:exponents:i}) and (\ref{cor:exponents:ii}) hold for the histogram learner with the lower bounds on $\gamma$ replaced by $\gamma\geq1/(d_c+3)$ for the hazard and $\gamma\geq1/(d_c+2)$ for the cumulative hazard, and widths $n^{-1/(d_c+3)}\sqrt{\log n}$ and $n^{-1/(d_c+2)}\sqrt{\log n}$ at these values.
\end{enumerate}
\end{corollary}

\begin{remark}\label{rem:exponents}
\begin{enumerate}
\item The same exponent for the numerator and the denominator, $\beta=\gamma$, as in \cite{AP2026}, is always admissible for the hazard function. For the cumulative hazard, the error of $\hat R_n$ is compared with the smaller width $w_n(\gamma d_c)$, so that our bound requires a smoother denominator: $\beta<\gamma$, or $\beta=\gamma$ with $\nu>0$, the latter being the only option at $\gamma=1/(d_c+4)$. These conditions are only sufficient. The rate of Proposition~\ref{prop:denominator} carries a factor $\sqrt{\log n}$, while a factor $\sqrt{\log\log n}$ is expected, as in Remark~\ref{rem:nuisance_slow}; this would allow $\beta=\gamma$ with $\nu=0$.

\item The lower bounds on $\gamma$ are the rate-optimal exponents of the corresponding fixed-$n$ problems, balancing a bias of order $n^{-2\gamma}$, or $n^{-\gamma}$ for the histogram, with the stochastic error. At these exponents, the widths are the optimal uniform rates for a function of $d_c+1$ variables (hazard) or $d_c$ variables (cumulative hazard, after integration in time), up to the logarithmic factor, which is the price of uniformity in $t$ and in $n$ with our construction.

\item In practice, we recommend $\gamma$ slightly above $1/(d_c+4)$, with $\nu=0$ and $(1-\gamma d_c)/4\leq\beta<\gamma$. A single numerator then serves the hazard and the survival function, with a bias of smaller order than the width for both. At $\gamma=1/(d_c+4)$, the bias of the cumulative hazard is of the order of its standard deviation, and the margin is only logarithmic; see Section~\ref{sec:numerics}. If only the hazard is of interest, $\gamma$ slightly above $1/(d_c+5)$ suffices, and gives a narrower band.
\end{enumerate}
\end{remark}


\section{Incompleteness mechanisms}\label{sec:mechanisms}

We now specialize the results to right censoring, to the modified current status data of \cite{PR2006b}, where a duration is observed exactly, or only known to be larger, or only known to be smaller than the observed duration, and to cure models, where a fraction of the population never experiences the event \citep{PVK2020,AVK2018}. For each of them we give the latent-to-observed map and its inverse (in Section~\ref{subsec:observables} for right censoring), the at-risk function in the affine form \eqref{eq:R_affine}, the one-pass estimators of its components and of the nuisance parameter, and the resulting bands. Left truncation and right censoring, competing risks and further mechanisms are treated in Appendix~\ref{app:mechanisms}.

The verification of Assumption~\ref{ass:data}(\ref{ass:D4}) reduces, in every case, to the existence of bounded conditional densities for the observed durations and to the smoothness of a few regression-type functions in $\x_c$; we state the corresponding condition once.

\begin{assum}\label{ass:mech}
\begin{enumerate}
\item[(A1)]\condlabel{A1}{ass:mech:A1} For every observed duration $W\in\{Y,L,C\}$ entering the at-risk function, the conditional distribution of $W$ given $\X=(\z_c,\x_d)$ has a density on $[0,\tau_0]$ bounded uniformly in $\z_c\in\mathcal N$.

\item[(A2)]\condlabel{A2}{ass:mech:A2} For every indicator $u(t; \mathfrak z)$ entering the at-risk function and every $t\in[0,\tau_0]$, the map $\z_c\mapsto\EE[u(t;Z)\mid\X=(\z_c,\x_d)]f_\X(\z_c,\x_d)$ has partial derivatives of order two, bounded on $\mathcal N$ uniformly in $t$; $f_\X(\cdot,\x_d)$ is bounded on $\mathcal N$.
\end{enumerate}
\end{assum}

In each mechanism, $u_m$ is an indicator which is monotone in $t$, or, under left truncation, the difference of two such indicators displayed in \eqref{eq:R_ltrc}, so that the decomposition in Assumption~\ref{ass:data}(\ref{ass:D4}) is explicit, with $u_m^+=0$ and $u_m^-=-u_m$ when $u_m$ is non-increasing. Under (\ref{ass:mech:A1}), $t\mapsto\EE[\1\{Y\geq t\}\mid\X=(\z_c,\x_d)]$, $t\mapsto\EE[\1\{L\leq t\}\mid\X=(\z_c,\x_d)]$ and their analogues are Lipschitz uniformly in $\z_c$, and the same holds for events such as $\{Y\geq t,I=k\}$, since $\PP(t\leq Y<t',I=k\mid\X)\leq\PP(t\leq Y<t'\mid\X)$; together with (\ref{ass:mech:A2}) this gives Assumption~\ref{ass:data}(\ref{ass:D4}) for all the mechanisms considered, the positivity of the at-risk function being stated separately in each case. Throughout, $\hat f_n$ denotes any of the three learners of Section~\ref{subsec:learners}, applied with the indicator $\1\{I=1\}$ of the exactly observed durations, and $\hat R_{m,n}$ denotes the one-pass kernel estimator \eqref{eq:psi_R}.

\subsection{Right censoring}\label{subsec:rc}

Under right censoring, $Y=T\wedge C$, $I=\1\{T\leq C\}$ and $T\perp C\mid\X$, with the map and its inverse given in Section~\ref{subsec:observables}. In the form \eqref{eq:R_affine},
\begin{equation}\label{eq:R_rc}
R(t,\x)=R_1(t,\x)=\PP(Y\geq t\mid\X=\x)f_\X(\x),\qquad u_1(t;\mathfrak z)=\1\{y\geq t\},\qquad \eta\equiv0.
\end{equation}
The function $u_1$ is non-increasing in $t$, and $R_1(t,\x)=S_{T\mid\X}(t-\mid\x)S_{C\mid\X}(t-\mid\x)f_\X(\x)$ is bounded away from zero on $[0,t_1]$ as soon as $S_{T\mid\X}(t_1\mid\x)S_{C\mid\X}(t_1\mid\x)>0$. Hence, under Assumptions~\ref{ass:data} and \ref{ass:mech}, Corollary~\ref{cor:exponents} applies without any nuisance parameter, and gives the bands of Theorem~\ref{theo:hazard} for $\lambda_{T\mid\X}(\cdot\mid\x)$, and those of Theorem~\ref{theo:cumhaz} and Corollary~\ref{cor:survival} for $\Lambda_{T\mid\X}(\cdot\mid\x)$ and $S_{T\mid\X}(\cdot\mid\x)$.

\subsection{Modified current status}\label{subsec:mcs}

This is Model 1 of \cite{PR2006b}, extended to covariates. Besides $T$ and the censoring variable $C$, there is a latent Bernoulli variable $\Delta$ with $\PP(\Delta=1\mid\X=\x)=\mathfrak p(\x)\in(0,1]$, with $T\perp C\mid\X$ and $\Delta\perp(T,C)\mid\X$, and the observed variables are
\begin{equation}\label{eq:mcs_obs}
(Y,I)=(C,0)\ \text{ if }C<T,\qquad (Y,I)=(T,1)\ \text{ if }T\leq C,\ \Delta=1,\qquad (Y,I)=(C,2)\ \text{ if }T\leq C,\ \Delta=0.
\end{equation}
The case $\mathfrak p\equiv1$ is right censoring, and $\mathfrak p\equiv0$ gives current status data, for which no explicit inverse exists. The latent-to-observed map is
\begin{gather*}
H_0(\dd t\mid\x)=S_{T\mid\X}(t\mid\x)F_{C\mid\X}(\dd t\mid\x),\qquad H_1(\dd t\mid\x)=\mathfrak p(\x)S_{C\mid\X}(t-\mid\x)F_{T\mid\X}(\dd t\mid\x),\\
H_2(\dd t\mid\x)=\{1-\mathfrak p(\x)\}F_{T\mid\X}(t\mid\x)F_{C\mid\X}(\dd t\mid\x).
\end{gather*}
Integrating $H_1(\dd s\mid\x)+\mathfrak p(\x)H_0(\dd s\mid\x)$ over $[t,\infty)$ gives $\mathfrak p(\x)S_{T\mid\X}(t-\mid\x)S_{C\mid\X}(t-\mid\x)=H_1([t,\infty)\mid\x)+\mathfrak p(\x)H_0([t,\infty)\mid\x)$, which is the denominator of Equation (2.5) of \cite{PR2006b}, where the indices $0$ and $1$ are exchanged, and $\mathfrak p(\x)=\PP(I=1\mid\X=\x)/\PP(I\neq0\mid\X=\x)$. Hence \eqref{eq:R_affine} holds with
\begin{equation}\label{eq:R_mcs}
R_1(t,\x)=\PP(Y\geq t,I=1\mid\X=\x)f_\X(\x),\quad R_2(t,\x)=\PP(Y\geq t,I=0\mid\X=\x)f_\X(\x),\quad \eta=\mathfrak p(\x),
\end{equation}
that is $u_1(t; \mathfrak z)=\1\{y\geq t,\iota=1\}$ and $u_2(t;\mathfrak z)=\1\{y\geq t,\iota=0\}$, both non-increasing in $t$. The nuisance parameter is estimated by the ratio of two one-pass kernel estimators,
\begin{equation}\label{eq:hat_p}
\hat\eta_n=\hat{\mathfrak p}_n(\x)=\frac{\sum_{i\leq n}\1\{I_i=1\}\boldsymbol K_{g(i)}(\X_{i,c}-\x_c)\1\{\X_{i,d}=\x_d\}}{\sum_{i\leq n}\1\{I_i\neq0\}\boldsymbol K_{g(i)}(\X_{i,c}-\x_c)\1\{\X_{i,d}=\x_d\}},
\end{equation}
whose numerator and denominator satisfy \eqref{eq:rate_R}, by the proof of Proposition~\ref{prop:denominator}, under (\ref{ass:mech:A2}) for the indicators $\1\{\iota=1\}$ and $\1\{\iota\neq0\}$ and $\PP(I\neq0\mid\X=\x)>0$. Hence $|\hat\eta_n-\eta(\x)|=O(r^R_n)$ almost surely, the same rate as the at-risk components, and Corollary~\ref{cor:exponents} applies directly: the bands for $\lambda_{T\mid\X}$, $\Lambda_{T\mid\X}$ and $S_{T\mid\X}$ are valid under the same conditions on $(\gamma,\beta)$ as under right censoring, with $\inf_{[0,t_1]}R(\cdot,\x)>0$ guaranteed by $\mathfrak p(\x)S_{T\mid\X}(t_1\mid\x)S_{C\mid\X}(t_1\mid\x)>0$. The estimator $\hat R_n=\hat R_{1,n}+\hat{\mathfrak p}_n\hat R_{2,n}$ is updated from four one-pass recursions, so that the procedure remains online. The simulations of Section~\ref{sec:numerics} use this mechanism.

\subsection{Cure models}\label{subsec:cure}

In a cure model, a fraction of the population never experiences the event: $T$ takes values in $[0,\infty]$, $\phi(\x)=\PP(T<\infty\mid\X=\x)\in(0,1)$ is the conditional probability of being susceptible, and
\begin{equation}\label{eq:cure_mixture}
S_{T\mid\X}(t\mid\x)=1-\phi(\x)+\phi(\x)S_{0}(t\mid\x),\qquad t\geq0,
\end{equation}
where $S_0(\cdot\mid\x)$ is the proper conditional survival function of the susceptibles, with hazard $\lambda_0(\cdot\mid\x)$ and cumulative hazard $\Lambda_0(\cdot\mid\x)$; see \cite{AVK2018} for a review and \cite{PVK2020} for the nonparametric identification. The targets are $\lambda_0$, $\Lambda_0$ and $S_0$. The cure fraction is identified under the sufficient follow-up assumption, which we state with a known constant $\tau$:
\begin{equation}\label{eq:suff_followup}
S_0(\tau\mid\x)=0
\qquad\text{and}\qquad \inf_{t\in[0,\tau]}R_1(t,\x)>0,
\end{equation}
where $R_1$ is the at-risk function of the mechanism without cure, applied to the improper variable $T$. Then $S_{T\mid\X}(\tau\mid\x)=1-\phi(\x)$. Tests of sufficient follow-up are available without covariates \citep{YM2024,XEBVK2024,MRS2024} and with categorical covariates \citep{YMVK2025}.

The incompleteness mechanism acts on $T$ whether or not $T$ is proper, so that the hazard $\lambda_{T\mid\X}=f_{T\mid\X}/S_{T\mid\X}$ of the improper variable has the representation $f_{Y,I,\X}(t,1,\x)/R_1(t,\x)$ of the mechanism without cure, with the nuisance-free denominator \eqref{eq:R_rc} under right censoring and \eqref{eq:R_ltrc} under left truncation and right censoring. The cure enters only through the relation between the hazards of $T$ and of the susceptibles: by \eqref{eq:cure_mixture} and \eqref{eq:suff_followup}, $f_{T\mid\X}=\phi f_0$ and $S_{T\mid\X}(t\mid\x)-S_{T\mid\X}(\tau\mid\x)=\phi(\x)S_0(t\mid\x)$ on $[0,\tau]$, whence, for $t\in[0,\tau)$,
\begin{equation}\label{eq:cure_hazard}
\lambda_0(t\mid\x)=\frac{\lambda_{T\mid\X}(t\mid\x)}{1-\exp\{-\varrho(t\mid\x)\}},\qquad \varrho(t\mid\x)=\int_t^{\tau}\lambda_{T\mid\X}(s\mid\x)\dd s=\Lambda_{T\mid\X}(\tau\mid\x)-\Lambda_{T\mid\X}(t\mid\x),
\end{equation}
and
\begin{equation}\label{eq:cure_cumhaz}
\Lambda_0(t\mid\x)=-\log\frac{S_{T\mid\X}(t\mid\x)-S_{T\mid\X}(\tau\mid\x)}{1-S_{T\mid\X}(\tau\mid\x)}=\mathfrak G\big(\Lambda_{T\mid\X}(t\mid\x),\Lambda_{T\mid\X}(\tau\mid\x)\big),
\end{equation}
with $\mathfrak G(u,w)=\log\{(1-e^{-w})/(e^{-u}-e^{-w})\}$. Equation \eqref{eq:cure_hazard} is \eqref{eq:R_affine} with $R=R_1(1-e^{-\varrho})$, that is $\eta=-\{1-\phi(\x)\}$ and $R_2=R_1/S_{T\mid\X}$.

\paragraph{Estimators.} Let $\hat\lambda_n=\hat f_n/\hat R_{1,n}$ be the estimator of $\lambda_{T\mid\X}$ of the mechanism without cure, $\hat\Lambda_n(t)=\int_0^t\hat\lambda_n$, $\hat\varrho_n(t)=\hat\Lambda_n(\tau)-\hat\Lambda_n(t)$, and
\begin{equation}\label{eq:cure_estimators}
\hat\lambda_{0,n}(t)=\frac{\hat\lambda_n(t)}{1-e^{-\hat\varrho_n(t)}},\qquad \hat\Lambda_{0,n}(t)=\mathfrak G\big(\hat\Lambda_n(t),\hat\Lambda_n(\tau)\big),\qquad \hat S_{0,n}(t)=e^{-\hat\Lambda_{0,n}(t)}.
\end{equation}
They are computed from the one-pass recursions of $\hat f_n$ and $\hat R_{1,n}$, with $\hat\Lambda_n$ evaluated on $[0,\tau]$. For small $n$, they may be undefined, since $1-e^{-\hat\varrho_n(t)}$ or $e^{-\hat\Lambda_n(t)}-e^{-\hat\Lambda_n(\tau)}$ may be nonpositive; they then take arbitrary bounded values, as in item~\ref{rem:band:conv} of Remark~\ref{rem:band}, which changes nothing for $n$ large, almost surely. For the cumulative hazard we also need the derivatives of $\mathfrak G$, evaluated at $(\Lambda_{T\mid\X}(t\mid\x),\Lambda_{T\mid\X}(\tau\mid\x))$ and estimated by plug-in:
\begin{equation}\label{eq:cure_c1c2}
c_1(t)=\frac{S_{T\mid\X}(t\mid\x)}{S_{T\mid\X}(t\mid\x)-S_{T\mid\X}(\tau\mid\x)},\quad c_2(t)=\frac{1}{\phi(\x)}-c_1(t)\leq0,\quad w_t(s)=c_1(t)\1\{s\leq t\}+c_2(t),
\end{equation}
so that $w_t(s)=1/\phi(\x)$ for $s\leq t$.

\begin{theorem}[Cure models]\label{theo:cure}
Let the mechanism without cure be right censoring or left truncation and right censoring (Appendix~\ref{app:ltrc}), let Assumptions~\ref{ass:data} and \ref{ass:mech} hold on $[0,\tau]$ in place of $[0,t_1]$, with $R=R_1$ in (\ref{ass:D4}) and with \eqref{eq:suff_followup}, and let $\mathcal T=[t_0,t_1]$ with $0<t_0<t_1<\tau<\tau_0$ and $S_0(t_1\mid\x)>0$.
\begin{enumerate}
\item[(i)]\condlabel{i}{theo:cure:i} (Hazard of the susceptibles) Under the conditions of Corollary~\ref{cor:exponents}(\ref{cor:exponents:i}), or (\ref{cor:exponents:iii}) for the histogram, the bands $[\hat\lambda_{0,n}(t)\pm\hat L^{\lambda}_n(t)]$, $t\in\mathcal T$, $n\geq1$, with $\hat L^\lambda_n$ from \eqref{eq:band_hazard} computed with the denominator
\[
\hat R_n(t)=\hat R_{1,n}(t)\big\{1-e^{-\hat\varrho_n(t)}\big\},
\]
form a $(1-\alpha)$-asymptotic confidence band sequence for $\lambda_0(\cdot\mid\x)$ on $\mathcal T$.

\item[(ii)]\condlabel{ii}{theo:cure:ii} (Cumulative hazard and survival of the susceptibles) Under the conditions of Corollary~\ref{cor:exponents}(\ref{cor:exponents:ii}), or (\ref{cor:exponents:iii}) for the histogram, let $\hat c_1,\hat c_2$ be the plug-in versions of \eqref{eq:cure_c1c2} with $\hat S_n=e^{-\hat\Lambda_n}$, let $\hat w_t(s)=\hat c_1(t)\1\{s\leq t\}+\hat c_2(t)$, and
\[
\hat v^{0}_n(t)=\mathfrak c_\Lambda\int_0^{\tau}\hat w_t(s)^2\frac{\hat f_n(s,1,\x)}{\hat R_{1,n}(s)^2}\dd s,\qquad \hat L^{0}_n(t)=\frac1n\mathcal B_\rho\big(\mathfrak V^\Lambda_n\hat v^0_n(\pi_n(t)),\alpha_{\pi_n(t)}\big),
\]
with $\mathfrak c_\Lambda$ and $\mathfrak V^\Lambda_n$ as in Propositions~\ref{prop:kernel}--\ref{prop:histogram} ($\mathfrak c_\Lambda=1$ for the histogram). Then the bands $[\hat\Lambda_{0,n}(t)\pm\hat L^0_n(t)]$ form a $(1-\alpha)$-asymptotic confidence band sequence for $\Lambda_0(\cdot\mid\x)$ on $\mathcal T$, and the bands $[\hat S_{0,n}(t)e^{-\hat L^0_n(t)},\ \hat S_{0,n}(t)e^{\hat L^0_n(t)}]$ one for $S_0(\cdot\mid\x)$. The pointwise statements hold likewise.
\end{enumerate}
\end{theorem}

The widths are of the same orders as without cure; the cure fraction only changes the variance profile. In (\ref{theo:cure:i}), the profile is that of $\lambda_{T\mid\X}$ divided by $\{1-e^{-\varrho}\}^2$, which   diverges   as $t\uparrow\tau$, where $S_0$ vanishes; this is why $t_1<\tau$. In (\ref{theo:cure:ii}), the negative weight $c_2(t)$ on $(t,\tau]$ is the price of estimating the cure fraction, and it vanishes as $t\downarrow0$. The proof, in Appendix~\ref{app:cure_more}, rests on a delta method for one-pass estimators: the errors of $\hat\Lambda_n(t)$ and $\hat\Lambda_n(\tau)$ are sums of the same contributions integrated against different weights, so that $\hat\Lambda_{0,n}$ is, up to a negligible remainder and a deterministic term, the one-pass estimator with contributions $\int_0^\tau w_t(s)\psi_i(s;\mathfrak z)R_1(s,\x)^{-1}\dd s$, and no joint approximation is needed. The smoothness of $f(\cdot,1,\x)$ at $\tau$ and the choice of $\tau$ (including the case of a threshold depending on the covariates) are discussed in Appendix~\ref{app:cure_more}. 


\section{Numerical experiments}\label{sec:numerics}

We study the bands for the hazard and the survival function under modified current status data, which contain right censoring as a special case, with the recursive kernel estimator. The purpose is to compare them with a fixed-sample band inspected repeatedly, and to measure the effect of the nuisance parameter, of the grid, of the variance estimator and of the tuning.

\vspace{-.3cm}

\paragraph{Design.} Modified current status data arise with recall: in a survey on the age at first drug use, or on the age at onset of a disease, a subject who experienced the event before the interview may not recall when \citep{HKJ1975,CRFM1991}. The covariates are $\X=(X_1,X_2,X_3)$, where $(X_1,X_2)$ has a Gaussian copula with correlation $0.6$ and uniform margins, and $X_3\mid(X_1,X_2)$ is Bernoulli with parameter $0.3+0.4X_1$. The hazard is $\lambda(t\mid\x)=\exp\{g(x_1,x_2)\}h_{x_3}(t\mid x_2)$, with $g(x_1,x_2)=-0.9+0.8\sin(\pi x_1)-0.8x_2^2$, $h_0(t\mid x_2)=\exp\{(0.2+0.8x_2)t\}$, which is increasing, and $h_1(t)=0.4+1.6\exp\{-(t-1)^2/(2\cdot0.35^2)\}$, a peak. The censoring variable is $C=\min(E,2.5)$, with $E$ exponential with rate $0.25+0.35X_1$ given $\X$, and $T\perp C\mid\X$, which gives about $43\%$ of observations with $I=0$. An event time $T\leq C$ is recalled with probability $\mathfrak p(\x)=1-q(0.75+0.5x_1)$, independently of $(T,C)$ given $\X$, as in \eqref{eq:mcs_obs}. The case $q=0$ is right censoring; with $q=0.15$, about $15\%$ of the event times are not recalled, that is $8\%$ of the observations. The recall indicators are generated after the other variables, so that, in each replication, the values of $q$ share the durations and the covariates. We consider four covariate values: $(0.5,0.5,0)$, $(0.5,0.5,1)$, $(0.25,0.25,1)$, along the correlation of $(X_1,X_2)$, and $(0.25,0.75,0)$, against it. A stress design, with a stronger curvature, a narrow   peak   and $64\%$ of observations with $I=0$, is studied in Appendix~\ref{app:stress}.

\vspace{-.3cm}

\paragraph{Procedures.} We take the Epanechnikov kernel with the boundary kernel of Proposition~\ref{prop:kernel}, $b(i)=i^{-\gamma}$, $h(i)=0.6i^{-\gamma}$, $g(i)=0.8i^{-\beta}\{\log(i+1)\}^{\nu}$, one numerator for the two targets, $\tau_0=2$, $\mathcal T=[0.25,1.75]$, and the time grid of mesh $1/128$ for the estimators and the coverage. The nuisance parameter is estimated by \eqref{eq:hat_p}, with the same bandwidth $g(i)$, and $\hat R_n=\hat R_{1,n}+\hat{\mathfrak p}_n\hat R_{2,n}$, as in Section~\ref{subsec:mcs}; for $q=0$, $\hat{\mathfrak p}_n=1$ and the procedure is the one for right censoring. Two tunings are used: $(\gamma,\beta,\nu)=(1/6,1/6,1/8)$, for which the cumulative hazard is at its rate-optimal exponent, and $(\gamma,\beta,\nu)=(1/5,1/6,0)$, which undersmooths both targets. We compare the bands with the fine grid, $\ell_n=\lceil2\log_2(n+1)\rceil$, using at all $t$ the single level $\alpha/\{2^{\ell_n+1}(\ell_n+1)(\ell_n+2)\}$ of item~\ref{rem:band:irreg} of Remark~\ref{rem:band}, and the bands with the grid \eqref{eq:rule_grid}, $\kappa_n=3$, using the levels \eqref{eq:alpha_s} with $\ell_0=\lceil\log_2(3\tau_0/b(n_0))\rceil$, so that the grid at $n_0$ is $G_{\ell_0}$, each with the plug-in variance and the local-maximum variance of Remark~\ref{rem:localmax}; the variance is evaluated at $t$, as in item~\ref{rem:band:var} of Remark~\ref{rem:band}. The inflation factor $1+\varepsilon$ of Proposition~\ref{prop:stopping} is not used. The parameter $\rho$ is tuned at $n_0=10^4$, with the true variance at $n_0$ as target value; in practice, the estimated variance is used instead, which only changes the target value. A target value computed from a pilot sample, not used afterwards, is covered by the theory, conditionally on the pilot, as in Remark~\ref{rem:aux_rand}. The benchmark is an oracle fixed-sample band: at each $n$, the critical value is the empirical $95\%$ quantile, over the replications, of the supremum over $t$ of the absolute error divided by the standard deviation estimated with the local-maximum variance, so that its coverage at a single look is $0.95$ by construction; this favors the benchmark. The streams have length $10^6$, the bands are inspected at 40 geometrically spaced times from $n_0$ on, and the level is $\alpha=0.05$. With 1000 replications, the standard error of a miscoverage near $0.05$ is about $0.007$.

\vspace{-.3cm}

\begin{figure}[ht!]
	\centering
	\includegraphics[width=0.9\textwidth]{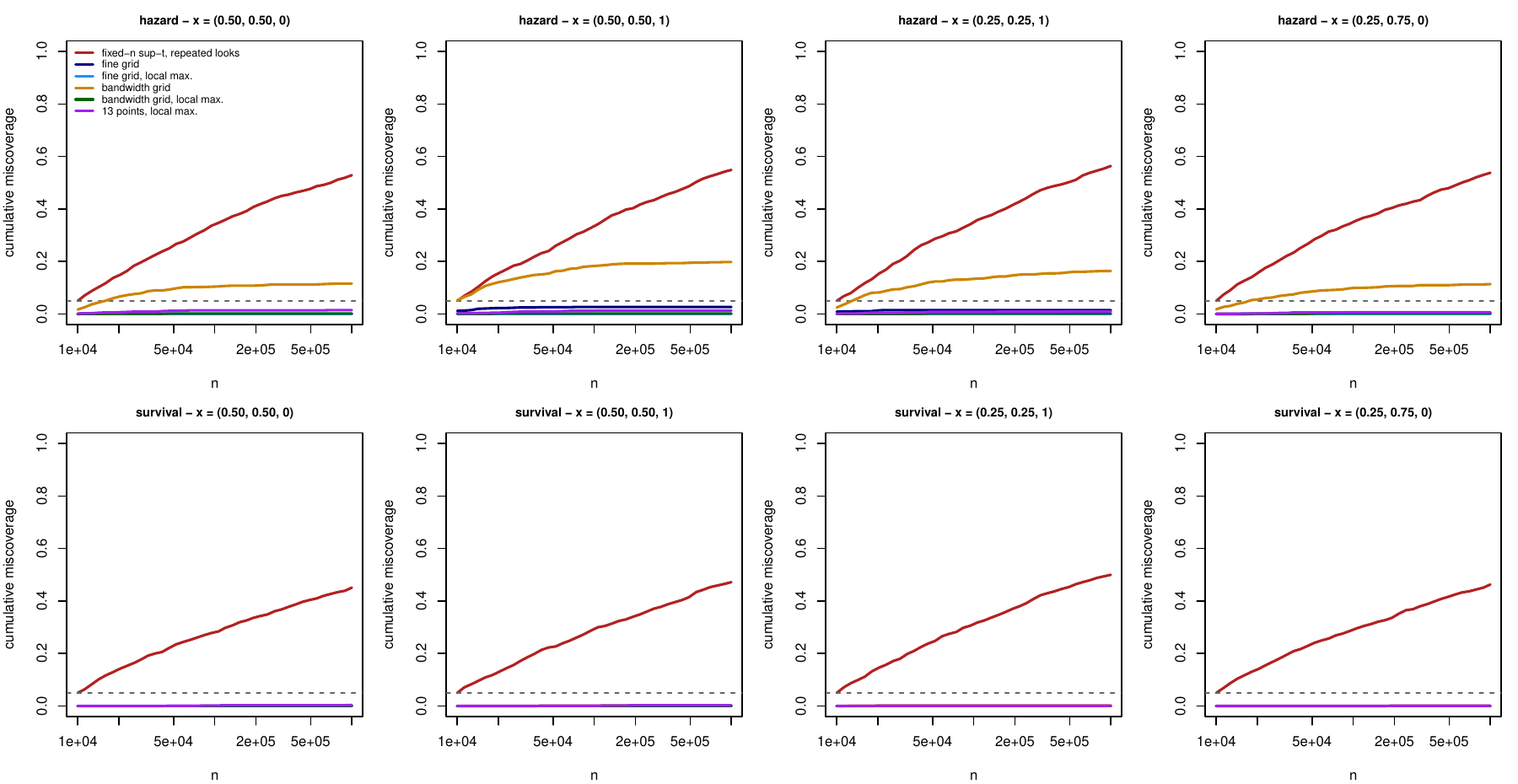}
	\caption{\small Cumulative miscoverage along the stream, $\gamma=1/5$, $q=0.15$. Top: hazard; bottom: survival function. Columns: the four covariate values. Curves as in the legend of the first panel; ``13 points'' corresponds to a fixed calibration grid as discussed below in the paragraph Recommendation. The dashed line is the nominal level.}\label{fig:sim_miss}
\end{figure}

\paragraph{Results.} Table~\ref{tab:sim} gives the cumulative miscoverage at $n=10^6$, that is the proportion of streams for which the band excluded the target at some inspected time and some $t$, and Figure~\ref{fig:sim_miss} displays it along the stream for $\gamma=1/5$ and $q=0.15$. A fixed-sample band inspected repeatedly fails in about half of the streams, and its miscoverage keeps increasing at $n=10^6$. This is the main   finding of this   section. Our bands behave differently according to the choices made.

\begin{table}[ht]
	\centering
	\footnotesize
	\setlength{\tabcolsep}{4pt}
	\begin{tabular}{lll|ccccc|cc}
		\hline
		& & & \multicolumn{5}{c|}{Largest miscoverage over the four $\x$} & \multicolumn{2}{c}{Ratio to the oracle width}\\
		$\gamma$ & $q$ & Target & Fixed-$n$ & Fine & Fine, l.m. & Band. & Band., l.m. & Fine, l.m. & Band., l.m.\\
		\hline
		$1/6$ & $0$    & $\lambda$ & 0.513 & 0.000 & 0.000 & 0.053 & 0.003 & 3.07--3.18 & 1.61--1.68\\
		$1/6$ & $0$    & $S$       & 0.453 & 0.000 & 0.000 & 0.002 & 0.000 & 3.30--3.86 & 1.71--2.01\\
		$1/6$ & $0.15$ & $\lambda$ & 0.512 & 0.000 & 0.000 & 0.062 & 0.004 & 3.07--3.29 & 1.61--1.73\\
		$1/6$ & $0.15$ & $S$       & 0.455 & 0.000 & 0.000 & 0.002 & 0.000 & 3.33--3.95 & 1.73--2.05\\
		$1/5$ & $0$    & $\lambda$ & 0.551 & 0.016 & 0.000 & 0.169 & 0.003 & 3.04--3.12 & 1.66--1.70\\
		$1/5$ & $0$    & $S$       & 0.486 & 0.000 & 0.000 & 0.002 & 0.000 & 3.73--3.86 & 2.01--2.07\\
		$1/5$ & $0.15$ & $\lambda$ & 0.564 & 0.027 & 0.000 & 0.198 & 0.003 & 3.05--3.21 & 1.66--1.75\\
		$1/5$ & $0.15$ & $S$       & 0.500 & 0.000 & 0.000 & 0.003 & 0.000 & 3.68--3.85 & 1.98--2.07\\
		\hline
	\end{tabular}
	\caption{Cumulative miscoverage at $n=10^6$ for the nominal level $0.05$, and range over the four covariate values of the ratio, at $n=10^6$, of the median half-widths to the oracle fixed-sample half-width; for $S$, these are the half-widths of the bands for $\Lambda$. $q=0$: right censoring; $q=0.15$: about $15\%$ of the event times are not recalled. Fixed-$n$: oracle band inspected repeatedly; Fine: fine grid with a single level; Band.: grid at the scale of the bandwidth with the levels \eqref{eq:alpha_s}; l.m.: local-maximum variance.}\label{tab:sim}
\end{table}

\begin{figure}[ht!]
	\centering
	\includegraphics[width=0.8\textwidth]{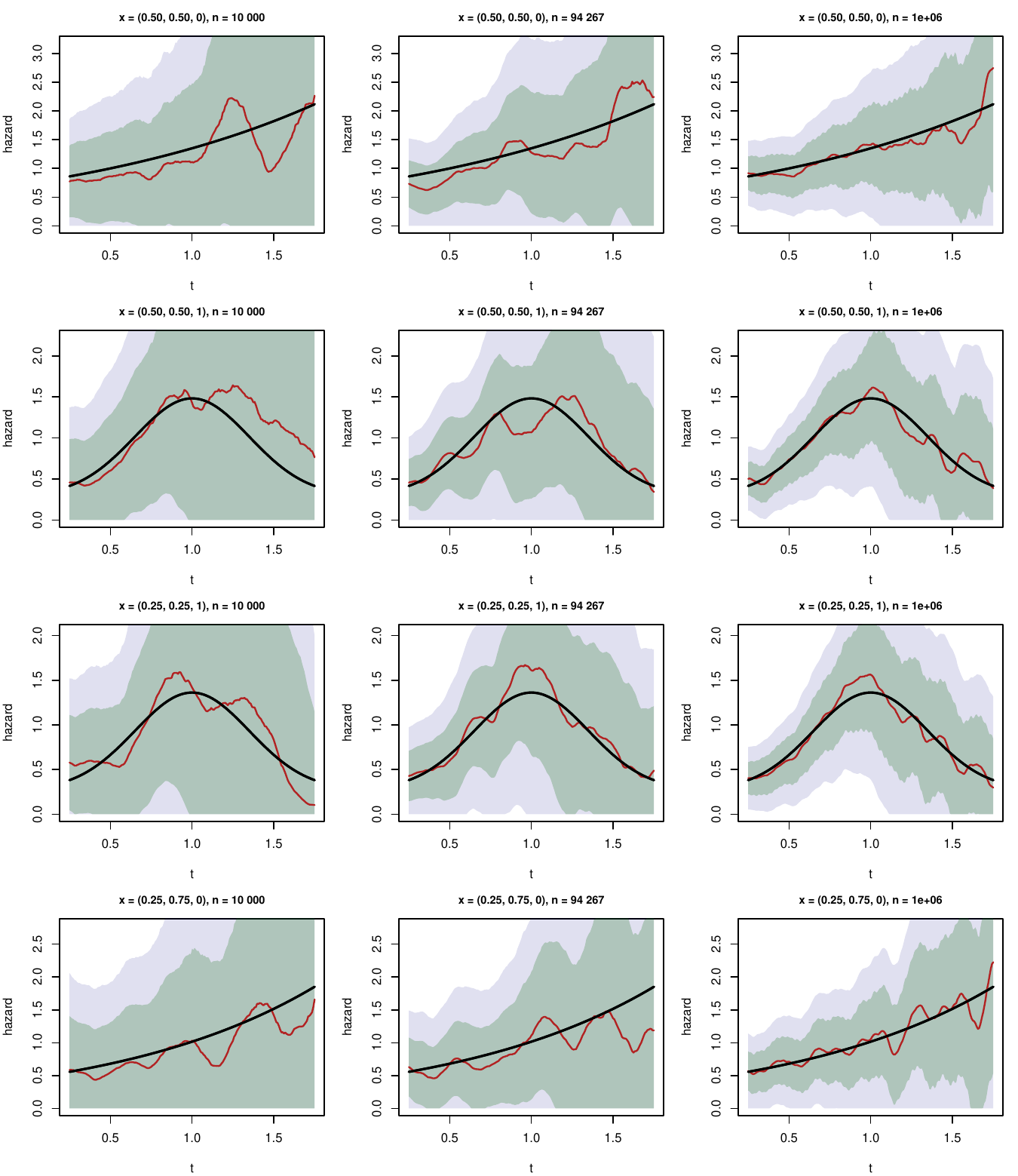}
	\caption{\small Bands for the hazard on one stream, $\gamma=1/5$, $q=0.15$, local-maximum variance: true hazard (black), estimator (red), band with the grid at the scale of the bandwidth (dark) and with the fine grid (light).}\label{fig:sim_bands}
\end{figure}

\begin{figure}[ht!]
	\centering
	\includegraphics[width=0.8\textwidth]{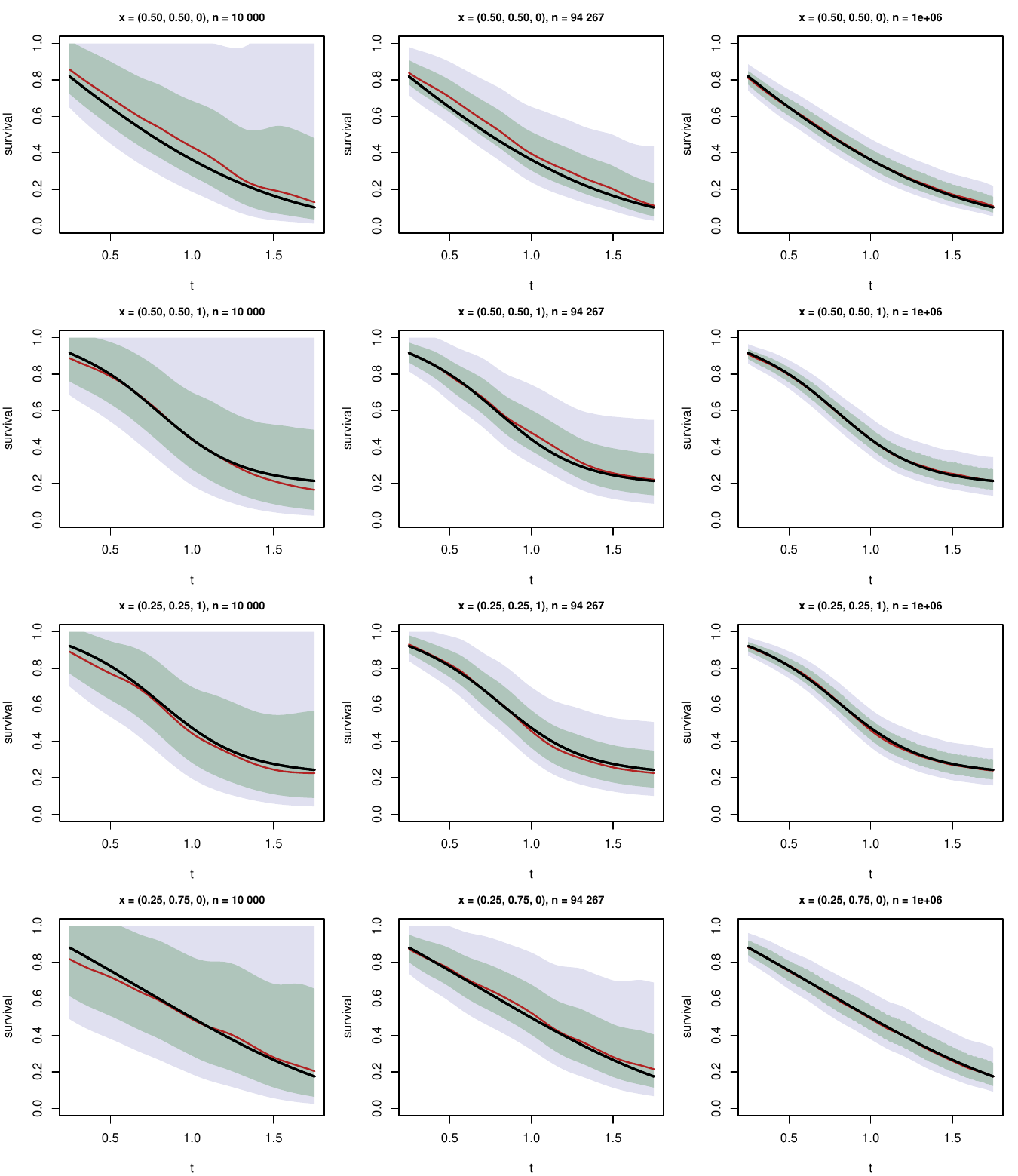}
	\caption[Bands for the survival function on one stream]{\small Bands for the survival function on the same stream as in Figure~\ref{fig:sim_bands}, $\gamma=1/5$, $q=0.15$, local-maximum variance: true survival function (black), estimator (red), band with the grid at the scale of the bandwidth (dark) and with the fine grid (light), truncated at $1$.}\label{fig:sim_bands_S}
\end{figure}

\emph{The nuisance parameter.} Its estimation has no visible cost. The variance profile is proportional to $f/R^2$, and both $f(\cdot,1,\x)$ and $R(\cdot,\x)$ are proportional to $\mathfrak p(\x)$, so that the loss of information due to recall multiplies the width by $\mathfrak p(\x)^{-1/2}$. The ratio of the widths with $q=0.15$ and $q=0$ matches this factor to the third decimal, for instance $1.084$ against $0.85^{-1/2}=1.085$, for both targets and both tunings, and the coverage of the bands with the local-maximum variance is unchanged. A level $q=0.30$, not shown, confirms these facts.

\emph{The grid.} The bands with the fine grid are very conservative, about three times as wide as the oracle band for the hazard and up to four times for the survival function. With the grid at the scale of the bandwidth, the width is reduced by a factor of about $0.53$ for $\gamma=1/6$ and $0.55$ for $\gamma=1/5$, close to the asymptotic values $0.58$ and $0.60$ of item~\ref{rem:band:width} of Remark~\ref{rem:band} for $d_c=2$. The ratio to the oracle width, $1.6$ to $2.1$, combines the cost of time-uniformity, the conservativeness of the union bound, and the bias, which the critical value of the oracle absorbs. For the hazard, the asymptotic ratio with this grid is $\{(1+a+2\gamma)/(2\gamma)\}^{1/2}$, about $2.2$ to $2.3$ here, and it is approached slowly.

\emph{The variance.} The plug-in variance is the main weakness. On the grid at the scale of the bandwidth, with $\gamma=1/5$, the miscoverage of the band for the hazard goes from $0.169$ to $0.003$ with the local-maximum variance under right censoring, and from $0.198$ to $0.003$ with $q=0.15$: the plug-in version deteriorates when the information decreases. Its failures occur early in the stream, with the estimate below the target. The price of the local-maximum variance is an increase of the width by $3$ to $14\%$ at $n=10^6$.

\emph{The tuning.} In this design, both tunings hold the nominal level with the local-maximum variance. In the stress design of Appendix~\ref{app:stress}, the exponent $\gamma=1/6$, rate-optimal for the survival function, fails for both targets because of the bias, mainly that of the denominator at one covariate value, while $\gamma=1/5$ stays close to the nominal level.

\emph{Recommendation.} We recommend the grid \eqref{eq:rule_grid}, with the constant $\kappa_n=3$ as a practical approximation of a slowly diverging $\kappa_n$, the local-maximum variance and an exponent $\gamma$ slightly above $1/(d_c+4)$, as in Remark~\ref{rem:exponents}. In this design, this gives a miscoverage of at most $0.003$ for the hazard and $0$ for the survival function, with or without recall. Figure~\ref{fig:sim_bands} shows that the band for the hazard narrows slowly, as expected with two continuous covariates, while the bands for the survival function, of a smaller order, are informative earlier  (Figure~\ref{fig:sim_bands_S}). We also computed the band with a fixed calibration grid of 13 equally spaced points in $\mathcal T$ (mesh $0.125$, not refined with $n$) and the level $\alpha/13$ at every point. This grid is outside the theory, since its mesh does not even tend to zero, so that the error between grid points is not controlled once the bandwidth falls below the mesh. It reduces the width by about $10\%$, but its miscoverage for the hazard, close to that of the recommended band at $n=10^4$, reaches $0.015$ in this design and $0.088$ in the stress design  at $n=10^6$. The grid at the scale of the bandwidth is thus close to the coarsest grid that preserves the coverage.

\emph{Limits.} Coverage is checked from $n_0=10^4$ on, on a finite time grid and over a finite horizon, and $\rho$ uses the true variance at $n_0$. The few failures of the recommended bands occur early in the stream, mostly near the right end of $\mathcal T$, where few subjects are at risk and the error of the at-risk estimator is not yet negligible; this effect is stronger in the stress design. The estimators are computed on $[0,t_1]$, so the local-maximum window is one-sided at $t_1$, which may contribute to these failures. These findings concern the designs considered, and the guarantee of Section~\ref{sec:inference} remains asymptotic.


\section{Discussion}\label{sec:discussion}

We proposed a general framework for the continuously updated estimation of conditional hazard functions, and for anytime-valid inference on the hazard, the cumulative hazard and the survival function. The framework rests on two ideas. The first one is the map viewpoint: a model is a map between the latent and the observed world, and for a large class of incompleteness mechanisms the inverse of this map is explicit, with the conditional hazard written as a density of the exactly observed durations divided by an at-risk function which is possibly affine in a nuisance parameter. The second one is the one-pass structure of the learners: each observation contributes once, and its contribution is never revised. This is what makes the estimators updatable at a small cost, without storing the past data, on a fixed grid of time points when the learner requires it, and, more importantly for us, what makes the contributions a single sequence of independent terms, to which a strong Gaussian approximation valid for all sample sizes simultaneously applies. The confidence bands are uniform in $t$ and in $n$, can shrink at the optimal rate up to a logarithmic factor, and are obtained with one-dimensional couplings at the points of a dyadic grid and a union bound, without any multivariate strong approximation. An extension to the reverse hazard, the natural object under left censoring, is outlined in Appendix~\ref{app:further_mech}. Below we discuss some limitations and further extensions.

\vspace{-.35cm}
  
\paragraph{Data arrival.} The observations are completed records, entering the stream once their follow-up has ended, in an order not depending on their values. This fits registries which report closed records, or subjects processed in the order of enrollment with a fixed maximal follow-up, at the price of a delay equal to this follow-up. Under calendar-time monitoring with staggered entry, the record of a subject changes until its follow-up ends, so that its contribution would be revised, which breaks the one-pass structure. Processing the records in the order of their completion would make the order depend on the durations. The safe logrank tests of \cite{GLS2021} and \cite{SPLG2024} allow calendar-time monitoring, but they concern testing, not the estimation of a hazard curve. Bands for the conditional hazard under calendar-time monitoring remain an open problem.

\vspace{-.35cm}

\paragraph{Uniformity in the covariates.} The bands are uniform in $t$ and in $n$, for a fixed covariate value $\x$. Uniformity in $\x_c$, over a compact subset of the interior of $\mathcal X_c$, can be obtained with the same tools: the grid becomes a nested dyadic grid in $(t,\x_c)$, with $2^{p\ell}$ points at depth $\ell$ instead of $2^{\ell}$, the allocation \eqref{eq:alpha_s} is adjusted accordingly, and the interpolation condition \eqref{cond:M} has to hold in both arguments, which is the case for the three learners. The order of the width is unchanged, the price being only in its constant. A uniformity in the discrete covariates is immediate, since $\mathcal X_d$ is finite.

\vspace{-.35cm}

\paragraph{Without continuous covariates.} When $d_c=0$, the bands for the hazard function are unchanged, with $p=1$ and at-risk components estimated without smoothing, while those for the cumulative hazard and survival function shrink at the rate $n^{-1/2}\sqrt{\log n}$. This is the bridge with time-uniform bands for distribution functions.   Two routes lead to these bands: the integrated learners, which require tuning and have a bias, and the Nelson--Aalen estimator, free of both, which we recommend. Computed exactly, it requires storing the past data, but a version binned on a fixed grid is computed in one pass. Under right censoring, it covers an unknown censoring law, with   the factor $\sqrt{\log n}$ of the normal mixture boundary instead of $\sqrt{\log\log n}$ in \cite{HR2022}. See Appendix~\ref{app:dc0}.

\vspace{-.35cm}

\paragraph{Exact coverage.} Because the Gaussian partial sums at the grid points are controlled separately and combined by a union bound, \eqref{eq:def_band} holds with an inequality which is not asymptotically sharp. With a grid at the scale of the bandwidth this loses little for the hazard, since the cost of the union bound then has the order suggested by extreme value theory; see item~\ref{rem:band:cov} of Remark~\ref{rem:band}. Asymptotically exact coverage would require a joint Gaussian approximation of the process $t\mapsto\sum_{i\leq n}\xi_i(t)$, uniform in $n$, and a time-uniform boundary for the supremum of the approximating Gaussian process; see Appendix~\ref{app:refinement}. Both ingredients are substantial, and the second one is, as far as we know, not available even for the simplest kernel density estimator with streaming data. On the other hand, we expect that the pointwise confidence sequences of   Theorem~\ref{theo:band}(\ref{theo:band:ii})   can be made sharp, in the sense of \citet[Def.~2.7 and Th.~2.8]{WASKR2024}, that is with a start time $m\to\infty$ and a parameter $\rho$ tuned to $m$; see also \cite{GK2025} for sharp asymptotic confidence sequences for a mean.

\vspace{-.35cm}

\paragraph{Data-dependent tuning and martingale contributions.} Definition~\ref{def:onepass} requires the tuning parameters of the $i$-th contribution to be deterministic, which makes the contributions independent. A bandwidth, a number of basis functions or a resolution selected from $Z_1,\ldots,Z_{i-1}$ is computationally admissible, and the contributions, centered conditionally on $Z_1,\ldots,Z_{i-1}$, then form a martingale difference sequence instead of independent variables. The same holds for stochastic gradient algorithms in function spaces \citep{DB2016}. The strong approximations we used have martingale counterparts, in which a Brownian motion and the Robbins boundary are evaluated at the predictable quadratic variation; see \cite{S1967} for the almost sure rates and \cite{WASKR2024} for asymptotic confidence sequences. Extending Lemma~\ref{lemma:blocked} to martingales, with exponential tails uniformly in $n$, is the step needed for bands uniform in $t$ with data-driven tuning. We leave it for future work.

\vspace{-.35cm}

\paragraph{Other learners.} The three learners detailed in   Sections~\ref{subsec:learners} and \ref{sec:examples}   were chosen for their simplicity. Local linear smoothers fit the framework, random partitions refined over time do so if their cut points are added to the calibration grid, whereas random features are not suited to sup-norm inference at nonparametric rates; see Appendix~\ref{subsec:further}; one-pass alternatives to the kernel for the denominator are studied in Appendix~\ref{app:denominator_alt}. We believe that the one-pass structure, rather than a specific learner, is the right level of generality for anytime-valid inference in this problem, and that it can serve as a guide for the design of streaming survival analysis procedures in the machine learning literature, where the continuous updating of a model is the rule and the control of the error at any data-dependent time is rarely addressed.

\section*{AI Use}
\vspace{-.3cm}
{\small Claude Opus 5.5 and Fable 5.1 (Anthropic) have been used for language improvement and proofreading; for proof and general error/flaw checking and correction; for guiding some of the theoretical derivations; and for prompt-based drafting of certain lemma and theorem statements and of parts of the accompanying proofs. The same tools were also used to assist with the code used in numerical experiments. The authors retain full responsibility for all aspects of this article. We certify that we have verified the details of all AI generated output and have converted any AI generated text to a personal style. To the best of our knowledge, all AI-influenced theoretical derivations, code and references are correct. Where AI assistance guided or contributed to proofs, we have identified and cited all relevant sources, including any not disclosed by the AI output.}\vspace{-.3cm}

\bibliographystyle{abbrvnat}
\bibliography{ref_clean_EM}

\newpage


\appendix
\addtocontents{toc}{\protect\setcounter{tocdepth}{2}}

\section*{Appendix}

The appendix contains complements and proofs, in the order of the sections of the paper. Appendix~\ref{app:sec_learners} gives sufficient conditions for Conditions \eqref{cond:E} and \eqref{cond:M} and complements on one-pass learners. Appendix~\ref{app:inference} contains the proofs of Section~\ref{sec:inference}, and discusses the choice of the weights, a coarser grid for the cumulative hazard, and the range and shape constraints. Appendix~\ref{app:refinement} indicates what exact coverage would require. Appendix~\ref{app:examples} contains the proofs of Section~\ref{sec:examples}   and other one-pass averages for the at-risk components, Appendix~\ref{app:mechanisms} further mechanisms and the proofs of Section~\ref{sec:mechanisms},   Appendix~\ref{app:dc0} the case without continuous covariates, and Appendix~\ref{app:stress} the stress design of the numerical experiments.

\renewcommand{\contentsname}{Contents of the appendix}
\tableofcontents

\section{Complements on Section~\ref{sec:learners}}\label{app:sec_learners}

\subsection{Sufficient conditions for (\ref{cond:E}) and (\ref{cond:M})}\label{app:learners}

\begin{lemma}\label{lemma:chaining}
	Assume that Condition \eqref{cond:I} holds.
	\begin{enumerate}
		\item[(a)]\condlabel{a}{lemma:chaining:a} If there exist $\kappa<(3+a)/2$ and $C$ such that $|\psi_i(t;\mathfrak z)-\psi_i(s;\mathfrak z)|\leq Ci^{\kappa}|t-s|$ for all $i\geq1$, $\mathfrak z\in\mathcal Z$ and $s,t\in\mathcal T$, then $\sup_{t\in\mathcal T}|\hat g_n(t)-\bar g_n(t)|=O(w_n)$ almost surely.

		\item[(b)]\condlabel{b}{lemma:chaining:b} If $\psi_i=\psi_i^{(1)}-\psi_i^{(2)}$, where, for each $\mathfrak z$, $t\mapsto\psi^{(j)}_i(t;\mathfrak z)$ is non-decreasing, $\sup_{t,\mathfrak z}|\psi^{(j)}_i(t;\mathfrak z)|\leq Ci^{a}$, $\sup_t\Var\{\psi^{(j)}_i(t;Z)\}\leq Ci^{a}$, and $t\mapsto\EE[\psi_i^{(j)}(t;Z)]$ is Lipschitz on $\mathcal T$ with a constant $Ci^{\kappa}$, $\kappa<(3+a)/2$, $j=1,2$, then $\sup_{t\in\mathcal T}|\hat g_n(t)-\bar g_n(t)|=O(w_n)$ almost surely.

		\item[(c)]\condlabel{c}{lemma:chaining:c} If, for each $n$, the functions $\hat g_n-\mathfrak r_n$ and $\bar g_n$ are constant on each cell of a partition of $\mathcal T$ into sets of the form $[s,s')\cap\mathcal T$ or $\{t_1\}$, with $s,s'\in G_{\ell_n}$, then $\mathfrak S_n(t)=\mathfrak S_n(\pi_n(t))$ for all $t\in\mathcal T$, so that Condition~\eqref{cond:M} holds, and $\sup_{t\in\mathcal T}|\hat g_n(t)-\bar g_n(t)|=O(w_n)$ almost surely.

		\item[(d)]\condlabel{d}{lemma:chaining:d} Assume that there exist a positive nonincreasing sequence $(\vartheta_i)$ and constants $C$ and $\kappa<(3+a)/2$ such that $i^{a}/\vartheta_i\leq Ci^{\kappa}$ and, for all $i\geq1$, $\mathfrak z\in\mathcal Z$ and $s,t\in\mathcal T$,
		\begin{align}
			|\psi_i(t;\mathfrak z)-\psi_i(s;\mathfrak z)|&\leq\frac{Ci^{a}}{\vartheta_i}|t-s|,\label{eq:incr_size}\\
			\Var\{\xi_i(t)-\xi_i(s)\}&\leq\frac{Ci^{a}}{\vartheta_i^2}|t-s|^2.\label{eq:incr_var}
		\end{align}
		Then the assumption of (\ref{lemma:chaining:a}) holds, and Condition \eqref{cond:M} holds as soon as $\varsigma_n=o(\vartheta_n)$.
	\end{enumerate}
	In each case, Condition \eqref{cond:B} implies Condition \eqref{cond:E}.
\end{lemma}

Case (\ref{lemma:chaining:a}) covers smooth contributions, case (\ref{lemma:chaining:b}) contributions of indicator type, such as those appearing in the at-risk function, which are of bounded variation in $t$, and case (\ref{lemma:chaining:c}) piecewise constant contributions, as in histograms. In case (\ref{lemma:chaining:b}), the function $t\mapsto\psi_i(t;\mathfrak z)$ is of bounded variation, written as a difference of non-decreasing pieces; the condition is formulated on the pieces, and the proof reduces to the monotone case. Case (\ref{lemma:chaining:b}) does not provide Condition \eqref{cond:M}, which is not needed for the denominators. Case (\ref{lemma:chaining:d}) is the sharp form of the condition on the grid: $\vartheta_i$ is the scale in $t$ of the $i$-th contribution. By \eqref{eq:incr_size} and \eqref{eq:incr_var}, the increment of the $i$-th contribution over a distance $\varsigma$ has the size and the variance allowed by Condition \eqref{cond:I}, multiplied by $\varsigma/\vartheta_i$ and $(\varsigma/\vartheta_i)^2$, respectively. The grid therefore only has to be finer than the scale of the latest contribution, which is the smallest one.

\begin{proof}[Proof of Lemma~\ref{lemma:chaining}]
Recall that $\mathfrak S_n(t)=\sum_{i\leq n}\xi_i(t)=n\{\hat g_n(t)-\bar g_n(t)-\mathfrak r_n(t)\}$. By Condition \eqref{cond:I}, for each fixed $t$, the variables $\xi_i(t)$ are independent, centered, bounded by $C_Mn^{a}$ for $i\leq n$, with $\sum_{i\leq n}\Var\{\xi_i(t)\}\leq C_Vn^{1+a}$. Bernstein's inequality gives, for $v_n=C_* \sqrt{n^{1+a}\log(n+1)}$,
\begin{equation}\label{eq:bernstein}
\PP\{|\mathfrak S_n(t)|\geq v_n\}\leq2\exp\left(-\frac{v_n^2/2}{C_Vn^{1+a}+C_Mn^{a}v_n/3}\right)\leq 2(n+1)^{-4},
\end{equation}
for $C_*$ large enough, because $n^{a}v_n=o(n^{1+a})$ when $a<1$, and $C_*$ does not depend on $t$. Since $v_n/n=C_*w_n$, it remains to control the supremum over $t$.

(\ref{lemma:chaining:a}) Let $\mathcal T_n$ be a grid of $n^2+1$ equidistant points in $\mathcal T$, including $t_0$ and $t_1$. By the union bound and \eqref{eq:bernstein}, $\PP\{\max_{s\in\mathcal T_n}|\mathfrak S_n(s)|\geq v_n\}\leq 4n^{-2}$, which is summable, so that $\max_{s\in\mathcal T_n}|\mathfrak S_n(s)|\leq v_n$ for $n$ large enough, almost surely, by the Borel--Cantelli lemma. For $t\in\mathcal T$ and $s$ the closest point of $\mathcal T_n$, $|\mathfrak S_n(t)-\mathfrak S_n(s)|\leq 2\sum_{i\leq n}Ci^{\kappa}|t-s|\leq Cn^{\kappa-1}$, deterministically, and $n^{\kappa-1}=o(nw_n)$ since $\kappa<(3+a)/2$. Hence $\sup_t|\mathfrak S_n(t)|/n=O(w_n)$ almost surely, and $\sup_t|\hat g_n-\bar g_n|=O(w_n)$ by Condition \eqref{cond:I}.

(\ref{lemma:chaining:b}) It suffices to treat $\psi^{(1)}_i$ and $\psi^{(2)}_i$ separately; let $\psi_i=\psi^{(1)}_i$, non-decreasing in $t$, with $\xi_i$, $\mathfrak S_n$ defined accordingly. By the bounds on the size and the variance of $\psi_i^{(1)}$, the bound \eqref{eq:bernstein} and the union bound over $\mathcal T_n$ apply as in (\ref{lemma:chaining:a}). For $t\in\mathcal T$, let $s_-\leq t\leq s_+$ be consecutive points of $\mathcal T_n$. By monotonicity, $\psi_i(s_-;Z_i)\leq\psi_i(t;Z_i)\leq\psi_i(s_+;Z_i)$, so that
\[
\mathfrak S_n(t)\leq \mathfrak S_n(s_+)+\sum_{i\leq n}\{\EE\psi_i(s_+;Z)-\EE\psi_i(t;Z)\}\leq \mathfrak S_n(s_+)+Cn^{\kappa-1},
\]
and similarly $\mathfrak S_n(t)\geq \mathfrak S_n(s_-)-Cn^{\kappa-1}$, by the Lipschitz property of $t\mapsto\EE\psi_i(t;Z)$. The conclusion follows as in (\ref{lemma:chaining:a}).

(\ref{lemma:chaining:c})   Since the endpoints of the cells belong to $G_{\ell_n}$ and $\pi_n(t)$ is the largest point of $G_{\ell_n}$ not exceeding $t$, the points $t$ and $\pi_n(t)$ belong to the same cell, so that $\mathfrak S_n(t)=\mathfrak S_n(\pi_n(t))$ and $\sup_t|\mathfrak S_n(t)|=\max_{s\in G_{\ell_n}}|\mathfrak S_n(s)|$. Since $|G_{\ell_n}|\leq2^{\ell_n+1}\leq Cn^C$ by \eqref{eq:ell_n}, the union bound and \eqref{eq:bernstein}, with the exponent $4$ increased to $C+2$ at the cost of a larger $C_\ast$, give a summable bound; the Borel--Cantelli lemma and Condition \eqref{cond:I} conclude.  

(\ref{lemma:chaining:d}) The bound \eqref{eq:incr_size} and $i^{a}/\vartheta_i\leq Ci^{\kappa}$ give the assumption of (\ref{lemma:chaining:a}). For Condition \eqref{cond:M}, put $\delta_n=\varsigma_n/\vartheta_n$, so that $\delta_n\to0$. We first reduce the supremum to a finite set of points, and then apply Bernstein's inequality to the increments of $\mathfrak S_n$.

\textit{Reduction to a finite set.} Let $\ell'_n=\max\{\ell_n,\lceil2\log_2(n+1)\rceil\}$ and let $\pi'_n$ be the projection on $G_{\ell'_n}$ by rounding down. For $t\in\mathcal T$, let $s=\pi'_n(t)$. Since $G_{\ell_n}\subset G_{\ell'_n}$, we have $\pi_n(s)=\pi_n(t)$, hence
\[
\mathfrak S_n(t)-\mathfrak S_n(\pi_n(t))=\{\mathfrak S_n(t)-\mathfrak S_n(s)\}+D_n(s),\qquad D_n(s)=\mathfrak S_n(s)-\mathfrak S_n(\pi_n(s)).
\]
Since $|t-s|\leq\tau_02^{-\ell'_n}\leq\tau_0(n+1)^{-2}$, the bound \eqref{eq:incr_size} and $i^{a}/\vartheta_i\leq Ci^{\kappa}$ give, deterministically, $|\mathfrak S_n(t)-\mathfrak S_n(s)|\leq Cn^{\kappa-1}=o(nw_n)$, as in (\ref{lemma:chaining:a}). It remains to bound $D_n(s)$ for $s\in G_{\ell'_n}$. If $\ell'_n=\ell_n$, then $\pi_n(s)=s$ and $D_n(s)=0$. Otherwise, $G_{\ell'_n}$ has at most $2(n+1)^2+1$ points.

\textit{Bernstein's inequality for the increments.} Fix $s\in G_{\ell'_n}$. The variable $D_n(s)=\sum_{i\leq n}\{\xi_i(s)-\xi_i(\pi_n(s))\}$ is a sum of independent centered variables. Since $|s-\pi_n(s)|\leq\varsigma_n$, and since $i^{a}/\vartheta_i$ and $i^{a}/\vartheta_i^2$ are nondecreasing in $i$, the bounds \eqref{eq:incr_size} and \eqref{eq:incr_var} give, for $i\leq n$,
\[
|\xi_i(s)-\xi_i(\pi_n(s))|\leq\frac{2Cn^{a}}{\vartheta_n}\varsigma_n=2Cn^{a}\delta_n,\qquad \sum_{i\leq n}\Var\{\xi_i(s)-\xi_i(\pi_n(s))\}\leq n\,\frac{Cn^{a}}{\vartheta_n^2}\varsigma_n^2=Cn^{1+a}\delta_n^2.
\]
These are the bounds of Condition \eqref{cond:I} used in \eqref{eq:bernstein}, with the size multiplied by $\delta_n$ and the variance by $\delta_n^2$. Taking the threshold $\delta_nv_n$, the factor $\delta_n^2$ cancels in the exponent of Bernstein's inequality:
\[
\begin{aligned}
\PP\{|D_n(s)|\geq\delta_nv_n\}&\leq2\exp\left(-\frac{\delta_n^2v_n^2/2}{Cn^{1+a}\delta_n^2+2Cn^{a}\delta_n^2v_n/3}\right)\\
&=2\exp\left(-\frac{v_n^2/2}{Cn^{1+a}+2Cn^{a}v_n/3}\right)\leq2(n+1)^{-4},
\end{aligned}
\]
for $C_*$ large enough, as in \eqref{eq:bernstein}. By the union bound over $G_{\ell'_n}$ and the Borel--Cantelli lemma, $\max_{s\in G_{\ell'_n}}|D_n(s)|\leq\delta_nv_n=C_*\delta_nnw_n$ for $n$ large enough, almost surely. This is $o(nw_n)$ because $\delta_n\to0$, which is the only place where the condition $\varsigma_n=o(\vartheta_n)$ is used.

In each case, $\sup_t|\hat g_n-g|\leq\sup_t|\hat g_n-\bar g_n|+\sup_t|\bar g_n-g|=O(w_n)+O(n^{-\mathfrak b})$ under \eqref{cond:B}.
\end{proof}


\subsection{Complements on one-pass learners}\label{app:learners_more}

\subsubsection{Orders of magnitude}\label{subsubsec:numbers}

Table~\ref{tab:numbers} gives the tuning parameters of the three learners of Section~\ref{subsec:learners} at the rate-optimal $\gamma$ of   Corollary~\ref{cor:exponents},   for one and two continuous covariates, with unit constants $c_b=c_h=c_J=1$ and the continuous arguments scaled to $[0,1]$; in practice, $c_b$, $c_h$ and $c_J$ are chosen by the usual rules of thumb, as discussed in \cite{AP2026}. The numbers of basis functions and of bins are small: with $d_c=1$, $J=10$ (eleven basis functions, eight of them cosines) and sixteen bins per axis suffice at $i=10^6$. This is the familiar observation that, under second order smoothness, the effective dimension of the problem grows slowly with the sample size; it also indicates that the choice of the constants $c_b$, $c_h$, $c_J$ matters more than the exponent in practice.

\begin{table}[ht]
\centering
\begin{tabular}{llcccc}
\hline
 & & $i=10^3$ & $i=10^4$ & $i=10^5$ & $i=10^6$\\
\hline
$d_c=1$, $\gamma=1/6$ & bandwidth $b(i)=h(i)=i^{-1/6}$ & 0.32 & 0.22 & 0.15 & 0.10\\
                      & series order  $J(i)=\lfloor i^{1/6}\rfloor$ & 3 & 4 & 6 & 10\\
$d_c=1$, $\gamma=1/4$ & bins per axis $2^{k_i}$, $k_i=\lfloor\log_2(i)/4\rfloor$ & 4 & 8 & 16 & 16\\
\hline
$d_c=2$, $\gamma=1/7$ & bandwidth $b(i)=h(i)=i^{-1/7}$ & 0.37 & 0.27 & 0.19 & 0.14\\
                      & series order $J(i)=\lfloor i^{1/7}\rfloor$ & 2 & 3 & 5 & 7\\
$d_c=2$, $\gamma=1/5$ & bins per axis $2^{k_i}$, $k_i=\lfloor\log_2(i)/5\rfloor$ & 2 & 4 & 8 & 8\\
\hline
\end{tabular}
\caption{Tuning parameters of the three learners at the rate-optimal exponent, with unit constants. The width of the band is of order $n^{-1/3}\sqrt{\log n}$ for $d_c=1$ and $n^{-2/7}\sqrt{\log n}$ for $d_c=2$ under second order smoothness, and $n^{-1/4}\sqrt{\log n}$, $n^{-1/5}\sqrt{\log n}$ under Lipschitz smoothness.}\label{tab:numbers}
\end{table}

\subsubsection{Memory and computation}\label{subsubsec:memory}

Concerning the memory, what has to be stored depends on the learner. The kernel learner stores its values at $M$ time points fixed in advance, with $M\gtrsim\tau_0/b(N)$ for a planned horizon $N$, and requires $O(Mb(i))$ operations per observation; at these points the estimator is exact, and the bands hold there as a restriction of the statement for all $t$. The series and histogram learners store $J(n)\asymp n^{\gamma}$ coefficients or $\asymp2^{k_n}$ counts, updated in $O(J(i))$ and $O(1)$ operations, and can be evaluated at any $t$ at any time. The components of the at-risk function are step functions of $t$; they are stored on a fixed grid of time points, with $O(1)$ operations per observation if the increments between grid points are stored, and the integrals defining $\hat\Lambda_n$ and $\hat\varrho_n$ are computed on the same grid. The results concern the exact integrals, so that the quadrature error on this grid must be $o(w_n(\gamma d_c))$. This constrains the number of grid points for the planned horizon. For the hazard, the series and histogram learners need no planned horizon.

\subsubsection{Other learners and limitations}\label{subsec:further}

\paragraph{Random partitions.} 
The deterministic dyadic partitions can be replaced by random partitions independent of the data, such as the Mondrian partitions with lifetime of order $i^{\gamma}$ \citep{LRT2014,MGS2020}, generated once and refined as $i$ increases, in the spirit of Remark~\ref{rem:aux_rand}. Averaging over a number of independent partitions growing with $n$ (a Mondrian forest) reduces the bias to the second order for twice differentiable functions \citep{MGS2020}. Conditionally on the partitions, the contributions are independent and computed in one pass. Condition \eqref{cond:M}, however, fails if a cut in $t$ is not a grid point. Such a cut, created at time $i_0$, separates $t$ from $\pi_n(t)$ in every later contribution, and $\mathfrak S_n(t)-\mathfrak S_n(\pi_n(t))$ is then of the order of $nw_n$. The cut points can instead be added to the calibration grid, conditionally on the partitions, and case (\ref{lemma:chaining:c}) of Lemma~\ref{lemma:chaining} then applies. Condition \eqref{cond:I} also needs care, since the smallest cells of a Mondrian partition are much smaller than the typical one. We do not pursue this extension.

\paragraph{Local linear smoothing.} For the components of the at-risk function, which are regression-type quantities, a recursive local linear estimator is obtained by solving, at time $n$, the weighted least squares problem with the recursive moments
\[
\frac1n\sum_{i\leq n}u_m(t;Z_i)\Big\{\frac{\X_{i,c}-\x_c}{g(i)}\Big\}^{\otimes r}\boldsymbol K_{g(i)}(\X_{i,c}-\x_c)\1\{\X_{i,d}=\x_d\},\qquad r=0,1,2,
\]
with the bandwidth $g(i)$ of Section~\ref{subsec:denominator}. Each moment is a one-pass average, and the local linear estimator is a smooth function of these moments, so that Remark~\ref{rem:remainder} applies. The gain is a smaller bias at the boundary of $\mathcal X_c$; in the interior, the leading terms coincide with those of the local constant estimator.

\paragraph{Other learners for the denominator.} The components of the at-risk function can also be estimated by one-pass averages of the other two types, smoothing in $\x_c$ only; see Appendix~\ref{app:denominator_alt}. A histogram in $\x_c$ has a first-order bias, due to the position of $\x_c$ in its cell. Unless the numerator is itself the histogram learner, this excludes the rate-optimal exponent for the cumulative hazard and, when $d_c\geq2$, for the hazard. A cell centered at $\x_c$, or a frequency polygon, which interpolates linearly between the cell counts, restores a second-order bias and remains a one-pass average with nonnegative weights; the conditions on the exponents are then those of Corollary~\ref{cor:exponents}. A cosine series in $\x_c$ also has a second-order bias at the interior point $\x_c$, without any adaptation to the boundary of $\mathcal X_c$, but its weights are signed and not local, so that the assumptions on the at-risk components are required on the whole of $\mathcal X_c$. We keep the kernel learner for simplicity.

\paragraph{Random features.} Random Fourier features  \citep{RR2007} replace a kernel by a Monte Carlo average of $\mathfrak M$ cosines with random frequencies drawn once, which gives a representation of the estimator in terms of $2 \, \mathfrak M$ fixed trigonometric functions and $O(\mathfrak M)$ updates per observation. In the present setting, the contributions are still computed in one pass, conditionally on the frequencies, but the approximation of a scaled kernel $K_b$ by $\mathfrak M$ features carries a Monte Carlo error. A heuristic computation gives, after the convolution with the sub-density, an error of order $b^{-p/2}\mathfrak M^{-1/2}$ in the supremum norm, the low frequencies being dominant. Keeping this random bias below the width of the band requires $\mathfrak M(i)\gg i$ features at time $i$, which   removes the computational advantage of   the approximation. Random features are thus not a suitable learner for sup-norm inference at the nonparametric rates considered here, and we omit them.


\section{Complements and proofs for Section~\ref{sec:inference}}\label{app:inference}

\subsection{Proofs of the tools}

\begin{proof}[Proof of Lemma~\ref{lemma:blocked}]
Let $I_k=\{2^k,\ldots,2^{k+1}-1\}$, $k\geq0$, and $\lambda_k=\{2C_M2^{(k+1)a}\}^{-1}$. For $j\in I_k$, $\lambda_k|\xi_j|\leq 1/2$, so that $\lambda_k\EE[|\xi_j|^3e^{\lambda_k|\xi_j|}]\leq\frac12e^{1/2}\EE\xi_j^2\leq\Var(\xi_j)$. Moreover, $B_k^2=\sum_{j\in I_k}\Var(\xi_j)\leq 2^kC_M^22^{2(k+1)a}$, so that $\lambda_kB_k\leq2^{k/2-1}$ and $\log(1+\lambda_kB_k)\leq k+1$. For each $k$, apply Lemma~\ref{lemma:sakhanenko} to $(\xi_j)_{j\in I_k}$, with $\lambda=\lambda_k$, and let
\[
\Delta^{(k)}=\max_{m\in I_k}\Big|\sum_{j=2^k}^m(\xi_j-\eta_j)\Big|,\qquad z_k(u)=\frac{1}{c\lambda_k}\big\{2(k+1)+u\big\}.
\]
Then $\PP(\Delta^{(k)}\geq z_k(u))\leq(1+\lambda_kB_k)e^{-2(k+1)-u}\leq e^{-(k+1)-u}$. Since the blocks of variables are independent, the constructions for the different $k$ can be taken independent, on the product of the probability spaces, as in Steps 1--2 of the proof of Theorem~2.3 in \cite{WLR2025}. On this space, $(\xi_j)_{j\geq1}$ has its original distribution and the $\eta_j$ are independent. By the union bound, the event $E_u=\{\Delta^{(k)}<z_k(u),\ \forall k\geq0\}$ has probability at least $1-e^{-u}\sum_{k\geq0}e^{-(k+1)}\geq1-e^{-u}$. On $E_u$, for $n\in I_K$, splitting the sum over the blocks,
\[
\Big|\sum_{j=1}^n(\xi_j-\eta_j)\Big|\leq\sum_{k=0}^K\Delta^{(k)}\leq\frac{2C_M2^a}{c}\{2(K+1)+u\}\sum_{k=0}^K2^{ka}\leq\frac{2C_M2^{2a}}{c(2^a-1)}\{2(K+1)+u\}\,2^{Ka}.
\]
Since $2^K\leq n$ and $K+1\leq\log_2(n)+1\leq 3\log(n+1)$, the first statement follows. For the second one, $\PP(\cup_{u\in\NN}E_u)=1$.
\end{proof}

\begin{proof}[Proof of Lemma~\ref{lemma:gluing}]
Disintegrate $Q_s(\dd\mathbf z,\dd\boldsymbol\eta)=P_{\mathbf Z}(\dd\mathbf z)\,\kappa_s(\mathbf z,\dd\boldsymbol\eta)$, which is possible since the spaces are Polish, and define the joint distribution of $(\mathbf Z,(\boldsymbol\eta^{(s)})_{s\in\mathcal S})$ as $P_{\mathbf Z}(\dd\mathbf z)\bigotimes_{s\in\mathcal S}\kappa_s(\mathbf z,\cdot)$, which exists by the Ionescu-Tulcea theorem.
\end{proof}
\noindent In our applications, $\xi_j(s)=F_{j,s}(Z_j)$ for measurable maps $F_{j,s}$. Lemma~\ref{lemma:blocked} provides a coupling of $(\xi_j(s))_{j\geq1}$ with Gaussian variables; gluing it with the joint distribution of $(\mathbf Z,(\xi_j(s))_{j\geq1})$ along $(\xi_j(s))_{j\geq1}$, in the same way, gives a measure $Q_s$ as in Lemma~\ref{lemma:gluing}.

\begin{proof}[Proof of Lemma~\ref{lemma:robbins}]
For each $\vartheta\in\RR$, $\exp(\vartheta T_n-\vartheta^2V_n/2)$ is a nonnegative martingale with respect to the natural filtration of $(\zeta_i)$, with initial value 1. Mixing over $\vartheta\sim N(0,\rho^2)$ and using Fubini's theorem,
\[
M_n=\int \exp\left(\vartheta T_n-\frac{\vartheta^2V_n}{2}\right)\frac{e^{-\vartheta^2/(2\rho^2)}}{\sqrt{2\pi}\rho}\dd\vartheta=\frac{1}{\sqrt{V_n\rho^2+1}}\exp\left(\frac{\rho^2T_n^2}{2(V_n\rho^2+1)}\right)
\]
is also a nonnegative martingale with $M_0=1$. By Ville's inequality, $\PP(\exists n: M_n\geq 1/\alpha)\leq\alpha$, and $M_n\geq 1/\alpha$ is equivalent to $|T_n|\geq\mathcal B_\rho(V_n,\alpha)$.
\end{proof}

\subsection{Proofs of the general results}

\begin{proof}[Proof of Proposition~\ref{prop:stopping}]
Let $\mathcal A$ be the event in \eqref{eq:def_band} and, for $m\geq1$, let $\mathcal A_m$ be the event that $L^*_n(t)\leq(1+\varepsilon)\hat L_n(t)$ and $U^*_n(t)\leq(1+\varepsilon)\hat U_n(t)$ for all $n\geq m$ and $t\in\mathcal T$. By \eqref{eq:def_band_ratio}, $\PP(\mathcal A_m)\to1$. On $\mathcal A\cap\mathcal A_m$, the inflated intervals contain the intervals with bounds $L^*_n(t),U^*_n(t)$ for all $n\geq m$ and $t\in\mathcal T$, hence they contain $\theta(t)$, and $\PP(\mathcal A\cap\mathcal A_m)\geq1-\alpha-\PP(\mathcal A_m^c)$. The statement on $\tau_m$ follows since $\tau_m\geq m$.
\end{proof}

Throughout, $\x$ is fixed and omitted from the notation when possible; $s$ denotes a point of $G$ and $t$ a point of $\mathcal T$. We write $\mathfrak B_n(s)=\mathcal B_\rho(V_n(s),\alpha_s)$ for the deterministic boundary built on the exact variances, and $\hat{\mathfrak B}_n(s)=\mathcal B_\rho(\mathfrak V_n\hat v_n(s),\alpha_s)$ for its computable counterpart, so that $\hat L_n(t)=\hat{\mathfrak B}_n(\pi_n(t))/n$.

\begin{proof}[Proof of Theorem~\ref{theo:band}(\ref{theo:band:i})]
\textit{Step 0 (Orders).} By Condition \eqref{cond:V}, $c_0'n^{1+a}\leq V_n(s)\leq C_0'n^{1+a}$ for all $s\in\mathcal T$ and $n$ large. Since $\mathcal B_\rho(V,\alpha')^2=\rho^{-2}(V\rho^2+1)\{\log(V\rho^2+1)+2\log(1/\alpha')\}$, using \eqref{eq:log_alpha_s} we   obtain   constants $0<c\leq C$ such that
\begin{equation}\label{eq:boundary_order}
c\,w_n\leq\frac{\mathfrak B_n(s)}{n}\leq C\,w_n\qquad\text{for all }s\in G_{\ell_n}\text{ and }n\text{ large}.
\end{equation}

\textit{Step 1 (Couplings at the grid points).} For each $s\in G$, Condition \eqref{cond:I} allows to apply Lemma~\ref{lemma:blocked} to $(\xi_i(s))_{i\geq1}$ with the constants $a$ and $C_M$, which do not depend on $s$; the couplings for all $s\in G$ are glued by Lemma~\ref{lemma:gluing}. This gives, on a single probability space, the data and, for each $s\in G$, independent Gaussian variables $\eta_i(s)\sim N(0,\Var\{\xi_i(s)\})$, $i\geq1$. Let $T_n(s)=\sum_{i\leq n}\eta_i(s)$. For $\varepsilon\in(0,1)$, let $u_s(\varepsilon)=\log(1/\varepsilon)+\log(1/\varpi_s)$ and
\[
A_\varepsilon=\Big\{\forall s\in G,\ \forall n\geq1:\ |\mathfrak S_n(s)-T_n(s)|\leq C_1n^{a}\big\{\log(n+1)+u_s(\varepsilon)\big\}\Big\}.
\]
By \eqref{eq:blocked} and the union bound, $\PP(A_\varepsilon^c)\leq\sum_{s\in G}e^{-u_s(\varepsilon)}=\varepsilon\sum_s\varpi_s\leq\varepsilon$. Hence $\PP(\cup_{m\geq2}A_{1/m})=1$ and, since $\log(1/\varpi_s)\leq C\log(n+1)$ for $s\in G_{\ell_n}$, almost surely
\begin{equation}\label{eq:unif_coupling}
\sup_{s\in G_{\ell_n}}|\mathfrak S_n(s)-T_n(s)|=O\big(n^{a}\log n\big).
\end{equation}

\textit{Step 2 (Gaussian side).} By Lemma~\ref{lemma:robbins}, applied for each $s\in G$ with the deterministic standard deviations $\sigma_i(s)=\Var\{\xi_i(s)\}^{1/2}$, $\PP(\exists n\geq1: |T_n(s)|\geq \mathfrak B_n(s))\leq\alpha_s$. By the union bound, the event
\[
\mathcal E=\{\forall s\in G,\ \forall n\geq1:\ |T_n(s)|<\mathfrak B_n(s)\}
\]
has probability at least $1-\sum_s\alpha_s\geq1-\alpha$.

\textit{Step 3 (Decomposition and the unknown band).} For $s\in G$, by \eqref{eq:onepass},
\[
g(s)-\hat g_n(s)=-\frac{T_n(s)}{n}+r_n(s),\qquad r_n(s)=-\frac{\mathfrak S_n(s)-T_n(s)}{n}-\{\bar g_n(s)-g(s)\}-\mathfrak r_n(s).
\]
For $t\in\mathcal T$, let $s=\pi_n(t)$ and $\omega_n(t)=\{g(t)-g(s)\}-\{\hat g_n(t)-\hat g_n(s)\}$, so that $g(t)-\hat g_n(t)=-T_n(s)/n+r_n(s)+\omega_n(t)$. Define
\[
L^*_n(t)=U^*_n(t)=\frac{\mathfrak B_n(\pi_n(t))}{n}+|r_n(\pi_n(t))+\omega_n(t)|,
\]
which are positive for all $n$. Then $|g(t)-\hat g_n(t)|\leq|T_n(\pi_n(t))|/n+|r_n(\pi_n(t))+\omega_n(t)|$, and on $\mathcal E$ this is smaller than $L^*_n(t)$ for all $n\geq1$ and all $t\in\mathcal T$, which proves \eqref{eq:def_band}.

\textit{Step 4 (The remainders are negligible).} With $s=\pi_n(t)$,
\[
\frac{L^*_n(t)}{\hat L_n(t)}=\frac{\mathfrak B_n(s)}{\hat{\mathfrak B}_n(s)}+\frac{n|r_n(s)+\omega_n(t)|}{\hat{\mathfrak B}_n(s)}.
\]
By \eqref{eq:boundary_order}, condition \eqref{eq:def_band_ratio} follows from the three almost sure statements\condlabel{a}{proof:band:a}\condlabel{b}{proof:band:b}\condlabel{c}{proof:band:c}
\[
\text{(\ref{proof:band:a}) }\sup_{s\in G_{\ell_n}}\Big|\frac{\hat{\mathfrak B}_n(s)}{\mathfrak B_n(s)}-1\Big|\to0,\qquad \text{(\ref{proof:band:b}) }\sup_{s\in G_{\ell_n}}|r_n(s)|=o(w_n),\qquad\text{(\ref{proof:band:c}) }\sup_{t\in\mathcal T}|\omega_n(t)|=o(w_n).
\]
(\ref{proof:band:a}) Condition \eqref{cond:V} gives $\sup_{s\in\mathcal T}|\mathfrak V_n\hat v_n(s)/V_n(s)-1|\to0$ almost surely, since $\hat v_n/v\to1$ and $V_n/(\mathfrak V_nv)\to1$ uniformly. Hence $\{\mathfrak V_n\hat v_n(s)\rho^2+1\}/\{V_n(s)\rho^2+1\}\to1$ uniformly, the difference of the logarithms in $\hat{\mathfrak B}_n(s)^2$ and $\mathfrak B_n(s)^2$ tends to $0$ uniformly, while these logarithms are larger than $c\log n$; this gives (\ref{proof:band:a}).

(\ref{proof:band:b}) By \eqref{eq:unif_coupling}, the first term of $r_n$ is $O(n^{a-1}\log n)=o(w_n)$ uniformly, since $a<1$. By Condition \eqref{cond:B}, the second term is $O(n^{-\mathfrak b})$ uniformly, and $n^{-\mathfrak b}\leq n^{-(1-a)/2}=o(w_n)$ since $2\mathfrak b\geq1-a$. The third term is $o(w_n)$ uniformly by Condition \eqref{cond:I}.

(\ref{proof:band:c}) By \eqref{eq:onepass}, with $s=\pi_n(t)$,
\[
\omega_n(t)=\{(g-\bar g_n)(t)-(g-\bar g_n)(s)\}-\frac{\mathfrak S_n(t)-\mathfrak S_n(s)}{n}-\{\mathfrak r_n(t)-\mathfrak r_n(s)\}.
\]
The first term is $O(n^{-\mathfrak b})=o(w_n)$ uniformly, by Condition \eqref{cond:B} and $2\mathfrak b\geq1-a$, and the other two are $o(w_n)$ uniformly in $t$, almost surely, by Conditions \eqref{cond:M} and \eqref{cond:I}.

Finally, the order of $\sup_t\hat L_n(t)$ follows from (\ref{proof:band:a}) and \eqref{eq:boundary_order}.
\end{proof}

\begin{proof}[Proof of Theorem~\ref{theo:band}(\ref{theo:band:ii})]
The proof is that of   Theorem~\ref{theo:band}(\ref{theo:band:i})   with the single point $s=t$ in place of the grid, $\alpha$ in place of $\alpha_s$ and $\omega_n\equiv0$. Lemma~\ref{lemma:blocked} applied to $(\xi_i(t))_{i\geq1}$ gives $\mathfrak S_n(t)-T_n(t)=O(n^a\log n)$ almost surely; Lemma~\ref{lemma:robbins} gives $\PP(\exists n: |T_n(t)|\geq\mathcal B_\rho(V_n(t),\alpha))\leq\alpha$; and the arguments of Step 4 apply at the point $t$, where $\mathcal B_\rho(V_n(t),\alpha)/n\asymp w_n$ because $V_n(t)\asymp n^{1+a}$ and $\log V_n(t)\asymp\log n$.
\end{proof}

\begin{proof}[Proof of Proposition~\ref{prop:ratio}]
Write $r=\inf_{\mathcal T}R$. By assumption, almost surely, $\inf_{\mathcal T}\hat R_n\geq r/2$ for $n$ large, and, by Condition \eqref{cond:E} for $\hat g_n$, $\sup_{\mathcal T}|\hat g_n|\leq C$ for $n$ large. Then
\[
\frac{\hat g_n(t)}{\hat R_n(t)}=\frac1n\sum_{i=1}^n\frac{\psi_i(t;Z_i)}{R(t)}+\mathfrak r_n^{(R)}(t),\qquad \mathfrak r^{(R)}_n(t)=\frac{\mathfrak r_n(t)}{R(t)}+\hat g_n(t)\,\frac{R(t)-\hat R_n(t)}{R(t)\hat R_n(t)},
\]
and $\sup_t|\mathfrak r^{(R)}_n(t)|\leq r^{-1}\sup_t|\mathfrak r_n(t)|+2Cr^{-2}\sup_t|\hat R_n(t)-R(t)|=o(w_n)$ almost surely. This is   a one-pass expansion   \eqref{eq:onepass} for $g/R$, with contributions $\psi_i/R$ and centered contributions $\xi_i/R$. Condition \eqref{cond:I} follows with the constants $C_M/r$ and $C_V/r^2$. Condition \eqref{cond:V} holds with $V_n^{(R)}(t)=V_n(t)/R(t)^2$, $v^{(R)}=v/R^2$ and $\hat v^{(R)}_n=\hat v_n/\hat R_n^2$, since $\hat R_n/R\to1$ uniformly. Condition \eqref{cond:B} holds since $|\bar g_n-g|/R\leq r^{-1}|\bar g_n-g|$, and Condition \eqref{cond:E} since
\[
\Big|\frac{\hat g_n}{\hat R_n}-\frac gR\Big|\leq\frac{|\hat g_n-g|}{\hat R_n}+|g|\,\frac{|\hat R_n-R|}{R\hat R_n}=O(n^{-\mathfrak b}+w_n)+o(w_n).
\]
For Condition \eqref{cond:M}, with $s=\pi_n(t)$,
\[
\frac{\mathfrak S_n(t)}{R(t)}-\frac{\mathfrak S_n(s)}{R(s)}=\frac{\mathfrak S_n(t)-\mathfrak S_n(s)}{R(t)}+\mathfrak S_n(s)\Big\{\frac1{R(t)}-\frac1{R(s)}\Big\},
\]
where the first term is $o(nw_n)$ by Condition \eqref{cond:M}. Since $\mathfrak S_n=n(\hat g_n-\bar g_n-\mathfrak r_n)$, Conditions \eqref{cond:E}, \eqref{cond:B}, \eqref{cond:I} and $2\mathfrak b\geq1-a$ give $\sup_s|\mathfrak S_n(s)|=O(nw_n)$ almost surely. As $1/R$ is Lipschitz, the second term is $O(\varsigma_nnw_n)=o(nw_n)$.
\end{proof}

\begin{proof}[Proof of Theorem~\ref{theo:hazard}]
By Assumption~\ref{ass:hazard}(\ref{ass:hazard:ii})--(\ref{ass:hazard:iii}), $R(\cdot,\x)$ is Lipschitz on $\mathcal T$, bounded away from zero, and
\[
\sup_t|\hat R_n(t)-R(t,\x)|\leq\sup_t|\hat R_{1,n}-R_1|+|\hat\eta_n|\sup_t|\hat R_{2,n}-R_2|+|\hat\eta_n-\eta(\x)|\sup_t|R_2|=o(w_n)
\]
almost surely, since $\hat\eta_n$ is bounded for $n$ large and $R_2$ is bounded. Proposition~\ref{prop:ratio} shows that $\hat\lambda_n=\hat f_n/\hat R_n$ satisfies the conditions of   Theorem~\ref{theo:band}(\ref{theo:band:i})   for $\lambda$, with $v^{(R)}=v/R^2$ and $\hat v^{(R)}_n=\hat v_n/\hat R_n^2$, and \eqref{eq:band_hazard} is \eqref{eq:band_width} for this estimator. The pointwise statement follows from   Theorem~\ref{theo:band}(\ref{theo:band:ii})   in the same way.
\end{proof}

The following lemma gives the uniform rate of the integrated estimator. It is used for Condition \eqref{cond:E} in the proof of Theorem~\ref{theo:cumhaz}, and to control the plug-in quantities in the cure models. The crude bound $\int|\hat\lambda_n-\lambda|$ is not sufficient, since it only gives the rate $w_n(a)$.

\begin{lemma}\label{lem:int_rate}
Let $0\leq\mathfrak s_0<\mathfrak s_1\leq\tau_0$. Assume that $(\hat f_n)$   admits a one-pass expansion   \eqref{eq:onepass} and satisfies \eqref{cond:I}, \eqref{cond:B} and \eqref{cond:E} on $[\mathfrak s_0,\mathfrak s_1]$ with exponents $a,\mathfrak b$, that $R(\cdot,\x)$ is Lipschitz and bounded away from zero on $[\mathfrak s_0,\mathfrak s_1]$, and that the integrated contributions $\psi^\Lambda_i(t;\mathfrak z)=\int_{\mathfrak s_0}^t\psi_i(s;\mathfrak z)R(s,\x)^{-1}\dd s$ satisfy the bounds of Assumption~\ref{ass:cumhaz}(\ref{ass:cumhaz:iii}) on $[\mathfrak s_0,\mathfrak s_1]$ with an exponent $a_\Lambda\leq a$. Assume moreover that, almost surely, $\sup_t|\hat R_n(t)-R(t,\x)|=O(r^R_n)$ and $\sup_t|\int_{\mathfrak s_0}^t\mathfrak r_n(s)R(s,\x)^{-1}\dd s|=O(r^R_n)$ on $[\mathfrak s_0,\mathfrak s_1]$, with $r^R_n\to0$. Then, with $\hat\Lambda_n(t)=\int_{\mathfrak s_0}^t\hat f_n/\hat R_n$ and $\Lambda(t)=\int_{\mathfrak s_0}^tf(s,1,\x)/R(s,\x)\dd s$,
\[
\sup_{t\in[\mathfrak s_0,\mathfrak s_1]}|\hat\Lambda_n(t)-\Lambda(t)|=O\big(n^{-\mathfrak b}+w_n(a_\Lambda)+r^R_n\big),\qquad \sup_{t\in[\mathfrak s_0,\mathfrak s_1]}\Big|\sum_{i\leq n}\xi^\Lambda_i(t)\Big|=O\big(nw_n(a_\Lambda)\big),
\]
almost surely.
\end{lemma}

\begin{proof}
Let $r=\inf_{[\mathfrak s_0,\mathfrak s_1]}R(\cdot,\x)$. Almost surely, for $n$ large, $\inf\hat R_n\geq r/2$ and $\sup|\hat f_n|\leq C$ on $[\mathfrak s_0,\mathfrak s_1]$, and
\begin{equation}\label{eq:Lambda_decomp}
\hat\Lambda_n(t)=\frac1n\sum_{i=1}^n\psi_i^\Lambda(t;Z_i)+\mathfrak r^\Lambda_n(t),\qquad \mathfrak r^\Lambda_n(t)=\int_{\mathfrak s_0}^t\frac{\mathfrak r_n(s)}{R(s,\x)}\dd s+\int_{\mathfrak s_0}^t\hat f_n(s)\,\frac{R(s,\x)-\hat R_n(s)}{R(s,\x)\hat R_n(s)}\dd s,
\end{equation}
so that $\sup_t|\mathfrak r^\Lambda_n(t)|=O(r^R_n)$. By \eqref{cond:B}, $\sup_t|\bar\Lambda_n(t)-\Lambda(t)|=O(n^{-\mathfrak b})$, with $\bar\Lambda_n=n^{-1}\sum_{i\leq n}\EE\psi^\Lambda_i(\cdot;Z_i)$. The centered integrated contributions are Lipschitz in $t$, by \eqref{cond:I}:
\[
|\xi^\Lambda_i(t)-\xi^\Lambda_i(s)|=\Big|\int_s^t\frac{\xi_i(u)}{R(u,\x)}\dd u\Big|\leq C_Mr^{-1}i^{a}|t-s|.
\]
The proof of Lemma~\ref{lemma:chaining}(\ref{lemma:chaining:a}) uses the Lipschitz property only through the centered contributions; since $a<1\leq(3+a_\Lambda)/2$, it applies with the growth exponent $a_\Lambda$ and gives $\sup_t|\sum_{i\leq n}\xi^\Lambda_i(t)|=O(nw_n(a_\Lambda))$ almost surely. The two bounds follow.
\end{proof}

\begin{proof}[Proof of Theorem~\ref{theo:cumhaz}]
Let $r=\inf_{[0,t_1]}R(\cdot,\x)$. As in the proof of Theorem~\ref{theo:hazard}, $\sup_{[0,t_1]}|\hat R_n-R|=o(w_n(a_\Lambda))$ almost surely, by Assumption~\ref{ass:cumhaz}(\ref{ass:cumhaz:iv}). Hence, by \eqref{eq:Lambda_decomp} with $\mathfrak s_0=0$ and the last part of Assumption~\ref{ass:cumhaz}(\ref{ass:cumhaz:iii}), $\sup_{[0,t_1]}|\mathfrak r^\Lambda_n|=o(w_n(a_\Lambda))$ almost surely. We check the conditions of   Theorem~\ref{theo:band}(\ref{theo:band:i})   for $\hat\Lambda_n$ as a one-pass estimator of $\Lambda$ on $\mathcal T^\Lambda$, with exponents $a_\Lambda$ and $\mathfrak b$. Condition \eqref{cond:I} is Assumption~\ref{ass:cumhaz}(\ref{ass:cumhaz:iii}) together with the bound on $\mathfrak r_n^\Lambda$; Condition \eqref{cond:V} is part of Assumption~\ref{ass:cumhaz}(\ref{ass:cumhaz:iii}). For Condition \eqref{cond:B}, $\bar\Lambda_n(t)-\Lambda(t)=\int_0^t\{\bar f_n(s)-f(s,1,\x)\}R(s)^{-1}\dd s$, which is $O(n^{-\mathfrak b})$ uniformly by Condition \eqref{cond:B} for $\hat f_n$ on $[0,t_1]$, and $2\mathfrak b\geq1-a_\Lambda$ by assumption. Condition \eqref{cond:E} follows from Lemma~\ref{lem:int_rate} on $[0,t_1]$, with $r^R_n=w_n(a_\Lambda)$. For Condition \eqref{cond:M}, with $s=\pi_n(t)$,
\[
|\mathfrak S^\Lambda_n(t)-\mathfrak S^\Lambda_n(s)|=\Big|\int_s^t\frac{\mathfrak S_n(u)}{R(u)}\dd u\Big|\leq\frac{|t-s|}{r}\sup_{u\in[0,t_1]}|\mathfrak S_n(u)|=O\big(\varsigma_nn\{w_n(a)+n^{-\mathfrak b}\}\big)=o(nw_n(a_\Lambda)),
\]
where $\mathfrak S^\Lambda_n=\sum_{i\leq n}\xi^\Lambda_i$ and $\sup_u|\mathfrak S_n(u)|=O(n\{w_n(a)+n^{-\mathfrak b}\})$ almost surely by \eqref{eq:onepass}, \eqref{cond:E}, \eqref{cond:B} and \eqref{cond:I} for $\hat f_n$ on $[0,t_1]$, and the last equality follows from Assumption~\ref{ass:cumhaz}(\ref{ass:cumhaz:v}), since $w_n(a)/w_n(a_\Lambda)=n^{(a-a_\Lambda)/2}$. Theorem~\ref{theo:band}(\ref{theo:band:i}) concludes, and Theorem~\ref{theo:band}(\ref{theo:band:ii}) gives the pointwise statement in the same way.
\end{proof}

\begin{proof}[Proof of Corollary~\ref{cor:survival}]
Let $[\hat\Lambda_n(t)-L^*_n(t),\hat\Lambda_n(t)+U^*_n(t)]$ be the unknown band sequence associated with Theorem~\ref{theo:cumhaz}, with $\sup_t|L^*_n/\hat L^\Lambda_n-1|\to0$ and $\sup_t|U^*_n/\hat L^\Lambda_n-1|\to0$ almost surely. Since $u\mapsto e^{-u}$ is decreasing, on the event of probability at least $1-\alpha$ where $\Lambda(t)\in[\hat\Lambda_n(t)-L^*_n(t),\hat\Lambda_n(t)+U^*_n(t)]$ for all $n$ and $t$, we have $S(t)=e^{-\Lambda(t)}\in[\hat S_n(t)e^{-U^*_n(t)},\hat S_n(t)e^{L^*_n(t)}]$ for all $n$ and $t$. This is a nonasymptotic band sequence centered at $\hat S_n$, with lower and upper half-widths $\hat S_n(1-e^{-U^*_n})$ and $\hat S_n(e^{L^*_n}-1)$. Let $\chi_-(u)=(1-e^{-u})/u$ and $\chi_+(u)=(e^u-1)/u$, extended by $\chi_\pm(0)=1$; both are continuous at $0$. Since $\sup_t\hat L^\Lambda_n(t)\to0$,
\[
\frac{1-e^{-U^*_n}}{1-e^{-\hat L^\Lambda_n}}=\frac{U^*_n}{\hat L^\Lambda_n}\,\frac{\chi_-(U^*_n)}{\chi_-(\hat L^\Lambda_n)}\to 1\qquad\text{and}\qquad\frac{e^{L^*_n}-1}{e^{\hat L^\Lambda_n}-1}=\frac{L^*_n}{\hat L^\Lambda_n}\,\frac{\chi_+(L^*_n)}{\chi_+(\hat L^\Lambda_n)}\to 1
\]
uniformly in $t$, almost surely, which is \eqref{eq:def_band_ratio}.
\end{proof}

\subsection{On the choice of the weights and the use of shape constraints}\label{app:weights}

Recall that $s$ denotes a grid point, that is, a point of $G$, and $t$ a location in $\mathcal T$ where the band is evaluated; at time $n$, the band at $t$ is calibrated at the grid point $\pi_n(t)\in G_{\ell_n}$.

\paragraph{The two roles of the weights.} In the proof of   Theorem~\ref{theo:band}(\ref{theo:band:i}),   the weights enter twice. In Step~1, they are used in the union bound over the couplings at the grid points $s\in G$, through $u_s(\varepsilon)=\log(1/\varepsilon)+\log(1/\varpi_s)$; this only affects the rate $O(n^a\log n)$ of the coupling, which is negligible with respect to the width, and any summable weights satisfying \eqref{eq:log_alpha_s} can be used there. In Step~2, they define the levels $\alpha_s$ of the Gaussian boundaries at the grid points $s$, and this is where they affect the width. The allocation \eqref{eq:alpha_s} factorizes as $\varpi_s=\varpi^{\mathrm{dep}}_{\ell(s)}\varpi^{\mathrm{loc}}_{\ell(s)}$, with $\varpi^{\mathrm{dep}}_\ell=\{(\ell-\ell_0+1)(\ell-\ell_0+2)\}^{-1}$ across depths and $\varpi^{\mathrm{loc}}_\ell=1/N_\ell$ across the grid points of depth $\ell$. Since the grid is refined as $n$ grows, the first factor is an allocation over resolution epochs, in the spirit of the summable weights over epochs used to stitch boundaries in time \citep{HRMS2021}; polynomial weights cost $O(\log\ell)$ in $\log(1/\alpha_s)$, geometric weights such as those of \citet{MH2023} cost $O(\ell)$. The second factor is an allocation over locations, which can be chosen freely.

\paragraph{A balanced allocation.} The irregularity of the half-width along $t$ described in item~\ref{rem:band:irreg} of Remark~\ref{rem:band} comes from counting each grid point $s$ once, at the depth $\ell(s)$ where it first appears. A balanced alternative, closer to \citet{MH2023}, indexes the union bound by the pairs $(\ell,s)$ with $\ell\geq\ell_0$ and $s\in G_\ell$, so that each grid point $s$ is counted at every depth $\ell\geq\ell(s)$, with
\[
\bar\alpha_\ell=\frac{\alpha\,\varpi^{\mathrm{dep}}_\ell}{|G_\ell|},\qquad \sum_{\ell\geq\ell_0}\sum_{s\in G_\ell}\bar\alpha_\ell=\alpha,
\]
and uses at time $n$ the same level $\bar\alpha_{\ell_n}$ at all locations $t\in\mathcal T$, in place of $\alpha_{\pi_n(t)}$ in \eqref{eq:band_width}. The proof of   Theorem~\ref{theo:band}(\ref{theo:band:i})    applies with one change: in Step~2, the event $\mathcal E$ is replaced by $\{\forall\ell\geq\ell_0,\ \forall s\in G_\ell,\ \forall n\geq1:\ |T_n(s)|<\mathcal B_\rho(V_n(s),\bar\alpha_\ell)\}$, of probability at least $1-\alpha$, and Step~3 uses, for each location $t$, the pair $(\ell_n,\pi_n(t))$. Steps~0 and~4 are unchanged, since $\log(1/\bar\alpha_{\ell_n})\leq\log(1/\alpha)+C\log(n+1)$. With this allocation, the half-width no longer varies irregularly along $t$. At the grid points $s$ of depth $\ell_n$, which form about half of $G_{\ell_n}$, the level $\bar\alpha_{\ell_n}$ is about half of the level $\alpha_s$ given by \eqref{eq:alpha_s}, a difference of $\log2$ in $\log(1/\alpha_s)$. At the grid points $s$ of smaller depth, $\bar\alpha_{\ell_n}$ is smaller than $\alpha_s$, so that the band is wider at the locations $t$ calibrated at these grid points. The constant level used for the fine grid in Section~\ref{sec:numerics} is of this type.

\paragraph{Target-oriented weights and their effect.} Within a depth, the weights can favor a region of interest, for instance a landmark time or the region where a decision is to be taken, by taking $\varpi^{\mathrm{loc}}_\ell(s)$ proportional to a fixed positive function of $s$, normalized over the grid points $s$ of depth $\ell$. The weights must be fixed before the data are observed. Their effect on the width is limited: by \eqref{eq:robbins_b}, the squared half-width at a location $t$ is proportional to $\log(\mathfrak V_n\hat v_n(\pi_n(t))\rho^2+1)+2\log(1/\alpha_{\pi_n(t)})$, and the first term grows like $(1+a)\log n$, so that changing the level $\alpha_{\pi_n(t)}$ by a fixed factor changes the width at $t$ by a factor which tends to one. The shape of the bands is driven by the variance profile, which is already adapted to each target. The half-width of the band for the hazard at $t$ is proportional to $\{\lambda(t)/R(t,\x)\}^{1/2}$: it is small where the hazard is small and large late in the follow-up, where few subjects are at risk. The variance profile $v^\Lambda(t)$ of the cumulative hazard vanishes as $t\downarrow0$, and the band for the survival function is multiplicative, hence narrow near the origin and shrinking with $S$.

\paragraph{Monotone targets and the grid.} For the monotone targets $\Lambda$ and $S$, the interpolation between grid points can also be obtained from monotonicity, as in \citet{MH2023}, instead of Condition~\eqref{cond:M}. For a location $t$ between two consecutive grid points $s_-\leq t\leq s_+$ of $G_{\ell_n}$, $\Lambda(s_-)\leq\Lambda(t)\leq\Lambda(s_+)$, so that the coverage at the grid points gives a band at every location $t$, with limits $\hat\Lambda_n(s_-)-\hat L^\Lambda_n(s_-)$ and $\hat\Lambda_n(s_+)+\hat L^\Lambda_n(s_+)$. This band contains the increment $\hat\Lambda_n(s_+)-\hat\Lambda_n(s_-)$, of order $\varsigma_n$, and keeps the order $w_n(a_\Lambda)$ of the width only if $\varsigma_n=o(w_n(a_\Lambda))$, that is, with a grid polynomially finer than required by Assumption~\ref{ass:cumhaz}(\ref{ass:cumhaz:v}), which only asks $\varsigma_n=o(n^{-\gamma/2})$ when $a-a_\Lambda=\gamma$, and even only $\varsigma_n\to0$ by Appendix~\ref{app:grid_cumhaz}.   By item~\ref{rem:band:width} of Remark~\ref{rem:band},   a finer grid increases the constant of the width. The smoothness of the centered process at the scale of the bandwidth, expressed by Condition~\eqref{cond:M}, is what allows a coarse grid, and monotonicity does not improve on it. It can be used here only to post-process the bands, as discussed in Appendix~\ref{subsec:app_range_shape}. The same comparison applies to a hazard which is known to be monotone,   as for increasing or decreasing failure rate distributions.

\subsection{A coarser grid for the cumulative hazard}\label{app:grid_cumhaz}

Assumption~\ref{ass:cumhaz}(\ref{ass:cumhaz:v}) requires $\varsigma_nn^{(a-a_\Lambda)/2}\to0$, that is $\varsigma_n=o(n^{-\gamma/2})$ for the learners of Section~\ref{sec:examples}. This comes from the Lipschitz bound used for Condition \eqref{cond:M} in the proof of Theorem~\ref{theo:cumhaz}. That bound ignores the size of the increments of the integrated process, which are of Brownian type. For the three learners,
\begin{equation}\label{eq:incr_Lambda}
\Var\{\xi^\Lambda_i(t)-\xi^\Lambda_i(s)\}\leq C\,i^{a_\Lambda}|t-s|,\qquad s,t\in[0,t_1].
\end{equation}
For the kernel learner, the increment of $\psi^\Lambda_i(\cdot;\mathfrak z)$ between $s$ and $t$ is $\1\{\iota=1\}\boldsymbol K_{h(i)}(\z_c-\x_c)\1\{\z_d=\x_d\}\int_s^tK_{b(i)}(y-u)R(u,\x)^{-1}\dd u$, with the boundary kernel near the origin. The integral is bounded by $C\min\{1,|t-s|/b(i)\}$ and vanishes unless $y\in[s-b(i),t+b(i)]$. The second moment is thus at most $Ch(i)^{-d_c}\min\{1,|t-s|/b(i)\}^2\{|t-s|+2b(i)\}\leq Ci^{a_\Lambda}|t-s|$. The histogram is handled in the same way. For the series learner, the integral in $u$ over $(s,t]$ is the projection of $\1_{(s,t]}$ on the first $J(i)+1$ basis functions, evaluated at $y$, whose squared $L^2$ norm is at most $|t-s|$.

The bound \eqref{eq:incr_Lambda}, together with $|\xi^\Lambda_i(t)-\xi^\Lambda_i(s)|\leq C\min\{i^{a_\Lambda},i^{a}|t-s|\}$, gives Condition \eqref{cond:M} for $\hat\Lambda_n$ as soon as $\varsigma_n\to0$. The argument is the one used for Condition \eqref{cond:M} in the proof of Proposition~\ref{prop:NA}: Bernstein's inequality for the increments over pairs of points of $G_{\ell'_n}$ at distance at most $2\varsigma_n$, with the threshold $\epsilon\{n^{1+a_\Lambda}\log n\}^{1/2}$, and the Lipschitz bound between $t$ and the nearest point of $G_{\ell'_n}$. We omit the details. For $d_c=0$, \eqref{eq:incr_Lambda} holds with $a_\Lambda=0$, and Proposition~\ref{prop:NA} is an instance of this argument.

Hence the second part of Assumption~\ref{ass:cumhaz}(\ref{ass:cumhaz:v}) can be replaced by $\varsigma_n\to0$, which is part of \eqref{eq:ell_n}. The bands for $\Lambda$ and $S$ can then be calibrated on a grid refined slowly, with $\ell_n=o(\log n)$, while the band for $\lambda$ keeps the grid \eqref{eq:rule_grid}. By item~\ref{rem:band:width} of Remark~\ref{rem:band}, the factor $1+a_\Lambda+2\gamma$ in the squared half-width becomes $1+a_\Lambda+o(1)$. The width is thus reduced by the factor $\{(1+a_\Lambda)/(1+a_\Lambda+2\gamma)\}^{1/2}$, about $0.89$ for $d_c=2$ and $\gamma=1/6$, and $0.88$ for $\gamma=1/5$. The interpolation error is then negligible only by a logarithmic factor, so that the gain may be smaller in finite samples. This grid is not used in Section~\ref{sec:numerics}.

\subsection{On the inference with range and shape constraints}\label{subsec:app_range_shape}

The bands proposed in Sections~\ref{subsec:hazard}--\ref{subsec:cumhaz}  can be adapted to the nature of the targets at no cost in coverage. Write $\underline S_n(t)=\hat S_n(t)e^{-\hat L^\Lambda_n(t)}$ and $\overline S_n(t)=\hat S_n(t)e^{\hat L^\Lambda_n(t)}$ for the limits of the band of Corollary~\ref{cor:survival}. Since $S_{T\mid\X}(\cdot\mid\x)$ is nonincreasing with values in $[0,1]$, these limits can be replaced by
\[
\overline S^{\,\mathrm m}_n(t)=\min\Big\{1,\ \inf_{u\in\mathcal T^\Lambda,\,u\leq t}\overline S_n(u)\Big\},\qquad \underline S^{\,\mathrm m}_n(t)=\sup_{u\in\mathcal T^\Lambda,\,u\geq t}\underline S_n(u),\qquad t\in\mathcal T^\Lambda,
\]
and, since $\Lambda_{T\mid\X}(\cdot\mid\x)$ is nonnegative and nondecreasing, the limits $\underline\Lambda_n=\hat\Lambda_n-\hat L^\Lambda_n$ and $\overline\Lambda_n=\hat\Lambda_n+\hat L^\Lambda_n$ can be replaced by $\max\{0,\sup_{u\leq t}\underline\Lambda_n(u)\}$ and $\inf_{u\geq t}\overline\Lambda_n(u)$, with $u\in\mathcal T^\Lambda$. The lower limit of the band for the hazard can be replaced by $\max\{0,\hat\lambda_n-\hat L^\lambda_n\}$, and, if the hazard is known to be monotone, as for increasing or decreasing failure rate distributions, the band for the hazard can be replaced by its monotone hull in the same way. The resulting bands are contained in the original ones. On the event considered in the proof of Proposition~\ref{prop:stopping}, the target lies in the original bands, inflated by the factor $1+\varepsilon$, for all $n\geq m$ and all locations $t$, hence also in their truncated monotone hulls. Proposition~\ref{prop:stopping} therefore holds for the post-processed bands, built from the inflated half-widths.

\section{A possible refinement: exact coverage via a joint Gaussian approximation}\label{app:refinement}

The bands of   Theorem~\ref{theo:band}(\ref{theo:band:i})    control the Gaussian partial sums $T_n(s)$ separately at each grid point and combine them by a union bound, with the allocation \eqref{eq:alpha_s}. This is what makes the construction simple: the only probabilistic inputs are a one-dimensional coupling at each point and Robbins' boundary. The price is a coverage which is not asymptotically exact. We indicate what an exact version would require, and why we did not pursue it.

The process $t\mapsto \mathfrak S_n(t)$ is, at each time $n$, a sum of independent processes $\xi_i(\cdot)$ with covariance $\Sigma_i(s,s')=\mathrm{Cov}\{\xi_i(s),\xi_i(s')\}$, which, for the learners of Section~\ref{subsec:learners}, is concentrated on $|s-s'|\lesssim i^{-\gamma}$ and has the form of a convolution kernel at scale $i^{-\gamma}$. A sharp band would be based on a Gaussian process $\mathbb T_n(t)=\sum_{i\leq n}\mathbb G_i(t)$, with $\mathbb G_i$ independent, centered, Gaussian, with covariance $\Sigma_i$, and on the two following ingredients.
\begin{enumerate}
\item A \emph{joint} strong approximation, uniform in $n$: $\sup_{t\in\mathcal T}|\mathfrak S_n(t)-\mathbb T_n(t)|=o(n^{(1+a)/2}(\log n)^{-1/2})$ almost surely, with a tail bound uniform in $n$. Restricted to the grid $G_{\ell}$, this is a multivariate version of Lemma~\ref{lemma:blocked} in dimension $|G_\ell|$, which grows with $n$; the multivariate results of \citet[Theorems~2--3]{Z2013} give an error depending polynomially on the dimension, so that a grid at the scale of the bandwidth would be needed, with the interpolation handled as in Condition \eqref{cond:M}. The blocks would again be dyadic, and the gluing across blocks would have to preserve the joint Gaussian structure, which the product construction does.
\item A \emph{time-uniform boundary for the supremum} of the approximating process: a function $\mathfrak B^{\sup}_n$ with $\PP(\exists n: \sup_{t}|\mathbb T_n(t)|/\sqrt{V_n(t)}\geq\mathfrak B^{\sup}_n)\leq\alpha$. Robbins' mixture argument does not apply directly to a supremum. One route is a time-uniform version of the Gaussian comparison and extreme value results for kernel-type processes, as in the fixed-$n$ theory of  \cite{BR1973}: the supremum of the normalized process over $\mathcal T$ behaves like $\sqrt{2\log(1/b)}$, with $b\asymp n^{-\gamma}$ the current bandwidth, which has the same leading order as $\{2\log(1/\alpha_s)\}^{1/2}\sim(2\gamma\log n)^{1/2}$ for the union bound over the grid \eqref{eq:rule_grid}. The term $\log V_n\asymp(1+a)\log n$ due to the normal mixture boundary would remain. With this boundary, an exact version would thus not change the leading constant of the width, and its gain would be a coverage which is asymptotically exact rather than conservative. This also requires $\rho$ tuned to the start of the monitoring; see Section~\ref{subsec:defs}.
\end{enumerate}
Both ingredients would require substantial work, and we are not aware of a result of the second kind even for the simplest kernel density estimator in the streaming setting. Since the union bound over a grid at the scale of the bandwidth already delivers the optimal rate, with the expected leading constant for the supremum part, we regard the exact version as a refinement rather than a necessity.


\section{Proofs for Section~\ref{sec:examples}}\label{app:examples}

\subsection{The augmented half-cosine basis}\label{app:basis}

Here we discuss the properties \eqref{cond:S1}--\eqref{cond:S3} of the augmented half-cosine basis considered in Section~\ref{subsubsec:series}. Property \eqref{cond:S1} is established in    \citet[Section~4.1 and Appendix~F]{PS2026}.   Property \eqref{cond:S2} is the purpose of the augmentation: for $J\geq2$, the span of $\varphi_0,\ldots,\varphi_J$ contains the quadratic polynomials, so that one may subtract from $q$ a quadratic polynomial $\mathsf p_q$ with $(q-\mathsf p_q)'(0)=(q-\mathsf p_q)'(\tau_0)=0$, and the cosine coefficients of $q-\mathsf p_q$ are then $O(j^{-3})$ after three integrations by parts, which gives $\sum_{j>J}O(j^{-3})=O(J^{-2})$ uniformly on $[0,\tau_0]$. Without the augmentation the cosine coefficients of $q$ decay only like $\{(-1)^jq'(\tau_0)-q'(0)\}j^{-2}$ and the uniform approximation error is of order $J^{-1}$. Property \eqref{cond:S3} is the localization of the squared projection kernel. For the pure cosine part, it follows from the representation of $D_J$ as a sum of two Dirichlet kernels, in the variables $y-t$ and $y+t$, whose squares are, up to normalization, Fej\'er kernels: the first one is an approximate identity at $y=t$, and the second one carries negligible mass near $y=t$ when $t$ is at a positive distance from $0$ and $\tau_0$. The augmentation adds to the kernel of the cosines a rank-two term $e_1(y)e_1(t)+e_2(y)e_2(t)$, where $e_1,e_2$ span the part of $t,t^2$ orthogonal to these cosines; $e_1,e_2$ have unit $L^2$-norm and are of order $J^{-1/2}$ uniformly on compacts of $(0,\tau_0)$, so that the additional terms in $\int D_J(y,t)^2q(y)\dd y$ are $O(1)$, negligible with respect to $J$. The details are given in   Lemma~\ref{lem:S3} below, 
together   with a bound on the partial sums of functions of bounded variation which is needed for the cumulative hazard. We note that  \citet{PS2026} establish an integrated form of \eqref{cond:S3}, against test functions in $t$, which is not sufficient here because the variance profile is needed uniformly in $t$. The restriction to compacts of $(0,\tau_0)$ in \eqref{cond:S3} is not an artifact: at $t=0$, the sum $\sum_{j\leq J}\varphi_j(0)^2$ is of order $J$ with a constant larger than the interior value $1/\tau_0$, so that the variance profile of the series learner is not of the form required in \eqref{cond:V} uniformly up to the origin. As for the kernel learner, this only matters for the hazard function, for which $t_0>0$.  

Throughout this subsection we take $\tau_0=1$; the general case follows by the change of variable $t\mapsto t/\tau_0$, which multiplies $\widetilde D_J$ by $\tau_0^{-1}$ and gives the constant $1/\tau_0$ in \eqref{cond:S3}. Let $c_0\equiv1$ and $c_j(t)=\sqrt2\cos(\pi jt)$, $j\geq1$, let $V_J=\mathrm{span}(c_0,\ldots,c_J)$, let $S^{\cos}_J$ be the orthogonal projection on $V_J$ in $L^2[0,1]$ and $D^{\cos}_J(y,t)=\sum_{j\leq J}c_j(y)c_j(t)$ its kernel. The augmented basis is an orthonormal basis of $W_J=V_J+\mathrm{span}(\mathsf p_1,\mathsf p_2)$, $\mathsf p_1(t)=t$, $\mathsf p_2(t)=t^2$, so that the projection $\widetilde S_J=S_{J+2}$ on the span of $\varphi_0,\ldots,\varphi_{J+2}$, with kernel $\widetilde D_J=D_{J+2}$, is the projection on $W_J$, and
\begin{equation}\label{eq:DJ_decomp}
\widetilde D_J=D^{\cos}_J+e_1\otimes e_1+e_2\otimes e_2,\qquad \widetilde S_J=S^{\cos}_J+P_{E_J},
\end{equation}
where $E_J=W_J\ominus V_J$ is two-dimensional, with orthonormal basis $e_1,e_2$, and $P_{E_J}$ is the projection on $E_J$. Since $S_J=\widetilde S_{J-2}$ and $D_J=\widetilde D_{J-2}$ for $J\geq3$, and $(J-2)/J\to1$, the statements below for $\widetilde S_J$ and $\widetilde D_J$ give \eqref{cond:S2}, \eqref{cond:S3} and the same bounds for $S_J$ and $D_J$. We use the Dirichlet kernel $\mathcal D_J(u)=1+2\sum_{j\leq J}\cos(ju)=\sin\{(J+\frac12)u\}/\sin(u/2)$, for which
\begin{equation}\label{eq:dirichlet}
D^{\cos}_J(y,t)=\tfrac12\{\mathcal D_J(\pi(y-t))+\mathcal D_J(\pi(y+t))\},\quad \int_{-\pi}^{\pi}\mathcal D_J(u)^2\dd u=2\pi(2J+1),\quad |\mathcal D_J(u)|\leq\frac{1}{|\sin(u/2)|}.
\end{equation}

\begin{lemma}\label{lem:S3}
\begin{enumerate}
\item[(i)]\condlabel{i}{lem:S3:i} (Augmentation) Every unit vector $e\in E_J$ satisfies $\|e\|_\infty\leq CJ^{1/2}$ and, for $\delta\in(0,1/2)$, $\sup_{t\in[\delta,1-\delta]}|e(t)|\leq C\delta^{-1}J^{-1/2}$, with an absolute constant $C$.
\item[(ii)]\condlabel{ii}{lem:S3:ii} (Partial sums of functions of bounded variation) There is an absolute constant $C$ such that, for every function $g$ of bounded variation on $[0,1]$ and every $J\geq1$,
\[
\|\widetilde S_Jg\|_\infty\leq C\{\|g\|_\infty+\mathrm{TV}(g)\},\qquad \|\widetilde S_Jg-g\|_{L^2[0,1]}\leq C\,\mathrm{TV}(g)\,J^{-1/2}.
\]
\item[(iii)]\condlabel{iii}{lem:S3:iii} (Approximation, \eqref{cond:S2}) If $q$ has two derivatives on $[0,1]$ and $q''$ is Lipschitz with constant $L$, then $\|\widetilde S_Jq-q\|_\infty\leq CLJ^{-2}$, with an absolute constant $C$.
\item[(iv)]\condlabel{iv}{lem:S3:iv} (Localization, \eqref{cond:S3}) Let $\mathcal Q$ be a uniformly bounded and uniformly equicontinuous family of functions on $[0,1]$ and $\delta\in(0,1/2)$. Then
\[
\sup_{q\in\mathcal Q}\ \sup_{t\in[\delta,1-\delta]}\Big|\frac1J\int_0^1\widetilde D_J(y,t)^2q(y)\dd y-q(t)\Big|\to0\qquad\text{as }J\to\infty.
\]
\end{enumerate}
\end{lemma}

\begin{proof}
(\ref{lem:S3:i}) Since $W_J=V_J\oplus E_J$ has dimension $J+3$, $E_J=\{\mathsf p-S^{\cos}_J\mathsf p: \mathsf p\in\mathrm{span}(\mathsf p_1,\mathsf p_2)\}$. For $\mathsf p=\alpha_1\mathsf p_1+\alpha_2\mathsf p_2$, two integrations by parts give, for $j\geq1$, since $\mathsf p''$ is constant,
\[
\langle \mathsf p,c_j\rangle=\sqrt2\int_0^1\mathsf p(y)\cos(\pi jy)\dd y=\frac{\sqrt2\{(-1)^j\mathsf p'(1)-\mathsf p'(0)\}}{(\pi j)^2}.
\]
Hence $\mathsf p-S^{\cos}_J\mathsf p=\sum_{j>J}\langle \mathsf p,c_j\rangle c_j$, the series converging uniformly, and
\begin{align*}
\|\mathsf p-S^{\cos}_J\mathsf p\|_\infty&\leq\frac{2(|\mathsf p'(0)|+|\mathsf p'(1)|)}{\pi^2}\sum_{j>J}\frac1{j^2}\leq\frac{2(|\mathsf p'(0)|+|\mathsf p'(1)|)}{\pi^2J},\\
\|\mathsf p-S^{\cos}_J\mathsf p\|_2^2&=\frac{2}{\pi^4}\sum_{j>J}\frac{\{(-1)^j\mathsf p'(1)-\mathsf p'(0)\}^2}{j^4}\geq\frac{c\{\mathsf p'(0)^2+\mathsf p'(1)^2\}}{J^3},
\end{align*}
the last inequality because two consecutive terms of the sum add up to at least $2\{\mathsf p'(0)^2+\mathsf p'(1)^2\}(j+1)^{-4}$. The first bound in (\ref{lem:S3:i}) follows. For the second one, write $(-1)^j\cos(\pi jt)=\cos(\pi j(1+t))$, so that
\[
(\mathsf p-S^{\cos}_J\mathsf p)(t)=\frac{2}{\pi^2}\Big\{\mathsf p'(1)\sum_{j>J}\frac{\cos(\pi j(1+t))}{j^2}-\mathsf p'(0)\sum_{j>J}\frac{\cos(\pi jt)}{j^2}\Big\}.
\]
By Abel summation, $|\sum_{j>J}\cos(j\theta)j^{-2}|\leq(J+1)^{-2}\sup_{N}|\sum_{j=J+1}^N\cos(j\theta)|\leq(J+1)^{-2}|\sin(\theta/2)|^{-1}$. For $t\in[\delta,1-\delta]$, $|\sin(\pi t/2)|$ and $|\sin(\pi(1+t)/2)|=\cos(\pi t/2)$ are both at least $\sin(\pi\delta/2)\geq\delta$. Hence $|(\mathsf p-S^{\cos}_J\mathsf p)(t)|\leq C(|\mathsf p'(0)|+|\mathsf p'(1)|)\delta^{-1}J^{-2}$, and dividing by $\|\mathsf p-S^{\cos}_J\mathsf p\|_2$ gives the claim.

(\ref{lem:S3:ii}) By \eqref{eq:dirichlet} and the substitution $s\mapsto-s$ in the second term, $S^{\cos}_Jg(y)=\int_0^1D^{\cos}_J(y,s)g(s)\dd s=(2\pi)^{-1}\int_{-\pi}^{\pi}\mathcal D_J(u)\tilde g(y+u/\pi)\dd u$, where $\tilde g$ is the even, $2$-periodic extension of $g$; this is the $J$-th partial sum of the Fourier series of $\tilde g$ at $\pi y$. Writing $\tilde g(y+u/\pi)+\tilde g(y-u/\pi)=2\tilde g(y)+\chi_y(u)$ for $u\in(0,\pi]$, we   obtain   $S^{\cos}_Jg(y)=g(y)+(2\pi)^{-1}\int_0^\pi\mathcal D_J(u)\chi_y(u)\dd u$, and an integration by parts with $\mathcal I_J(u)=\int_0^u\mathcal D_J$, which satisfies $\sup_{J,u}|\mathcal I_J(u)|<\infty$  \citep[Vol.~I, Chap.~II, \S9]{Z2002}, gives $|\int_0^\pi\mathcal D_J\chi_y|\leq\sup|\mathcal I_J|\{|\chi_y(\pi)|+\mathrm{TV}(\chi_y)\}\leq C\{\|g\|_\infty+\mathrm{TV}(g)\}$. Hence $\|S^{\cos}_Jg\|_\infty\leq C\{\|g\|_\infty+\mathrm{TV}(g)\}$. Next, an integration by parts gives $|\langle g,c_j\rangle|\leq\sqrt2\,\mathrm{TV}(g)/(\pi j)$, so that $\|g-S^{\cos}_Jg\|_2^2=\sum_{j>J}\langle g,c_j\rangle^2\leq C\,\mathrm{TV}(g)^2J^{-1}$. Since $e_k\perp V_J$, $|\langle g,e_k\rangle|=|\langle g-S^{\cos}_Jg,e_k\rangle|\leq\|g-S^{\cos}_Jg\|_2\leq C\,\mathrm{TV}(g)J^{-1/2}$, and by (\ref{lem:S3:i}), $\|P_{E_J}g\|_\infty\leq\sum_{k=1,2}|\langle g,e_k\rangle|\|e_k\|_\infty\leq C\,\mathrm{TV}(g)$. The first bound follows from \eqref{eq:DJ_decomp}. The second one follows from $\|\widetilde S_Jg-g\|_2\leq\|S^{\cos}_Jg-g\|_2$, since $W_J\supset V_J$.

(\ref{lem:S3:iii}) Let $\mathsf p_q(t)=q'(0)t+\frac12\{q'(1)-q'(0)\}t^2$ and $r=q-\mathsf p_q$, so that $r'(0)=r'(1)=0$, $r''$ is Lipschitz with constant $L$, and $\mathsf p_q\in W_J$, hence $\widetilde S_Jq-q=\widetilde S_Jr-r$. Three integrations by parts, the boundary terms vanishing because $r'(0)=r'(1)=0$ and $\sin(\pi j\cdot)$ vanishes at $0$ and $1$, give $\langle r,c_j\rangle=\sqrt2(\pi j)^{-3}\int_0^1r'''(y)\sin(\pi jy)\dd y$, where $r'''$ exists almost everywhere and $|r'''|\leq L$. Hence $|\langle r,c_j\rangle|\leq\sqrt2L(\pi j)^{-3}$, $\|S^{\cos}_Jr-r\|_\infty\leq\sum_{j>J}\sqrt2|\langle r,c_j\rangle|\leq CLJ^{-2}$ and $\|S^{\cos}_Jr-r\|_2\leq CLJ^{-5/2}$. As in (\ref{lem:S3:ii}), $\|P_{E_J}r\|_\infty\leq2\|r-S^{\cos}_Jr\|_2\max_k\|e_k\|_\infty\leq CLJ^{-2}$ by (\ref{lem:S3:i}), and \eqref{eq:DJ_decomp} concludes.

(\ref{lem:S3:iv}) Let $t\in[\delta,1-\delta]$ and $q\in\mathcal Q$, with $\|q\|_\infty\leq M$ and modulus of continuity bounded by $\omega$. By \eqref{eq:DJ_decomp},
\begin{multline*}
\int_0^1\widetilde D_J(y,t)^2q(y)\dd y=\int_0^1D^{\cos}_J(y,t)^2q(y)\dd y\\
+2\sum_{k=1,2}e_k(t)\int_0^1D^{\cos}_J(y,t)e_k(y)q(y)\dd y+\sum_{k,k'=1,2}e_k(t)e_{k'}(t)\int_0^1e_ke_{k'}q.
\end{multline*}
By the Cauchy--Schwarz inequality and $\|D^{\cos}_J(\cdot,t)\|_2^2=\sum_{j\leq J}c_j(t)^2\leq2J+1$, the second term is bounded by $4\max_k|e_k(t)|M\sqrt{2J+1}\leq CM\delta^{-1}$ by (\ref{lem:S3:i}), and the third one by $CM\delta^{-2}J^{-1}$. For the first term, by \eqref{eq:dirichlet}, $D^{\cos}_J(y,t)^2=\frac14\{\mathcal D_J(\pi(y-t))^2+2\mathcal D_J(\pi(y-t))\mathcal D_J(\pi(y+t))+\mathcal D_J(\pi(y+t))^2\}$. For $y\in[0,1]$, $\pi(y+t)\in[\pi\delta,\pi(2-\delta)]$, so that $|\mathcal D_J(\pi(y+t))|\leq1/\sin(\pi\delta/2)\leq\delta^{-1}$: the third part contributes at most $M\delta^{-2}/4$, and the cross part at most $\delta^{-1}M\int_0^1|\mathcal D_J(\pi(y-t))|\dd y\leq\delta^{-1}M\{2(2J+1)\}^{1/2}$, by the Cauchy--Schwarz inequality and \eqref{eq:dirichlet}. For the first part, substituting $u=\pi(y-t)$ and letting $\varepsilon\in(0,\pi\delta]$,
\[
\frac14\int_0^1\mathcal D_J(\pi(y-t))^2q(y)\dd y=\frac{1}{4\pi}\int_{|u|\leq\varepsilon}\mathcal D_J(u)^2q(t+u/\pi)\dd u+O\big(M\varepsilon^{-2}\big),
\]
because $[-\varepsilon,\varepsilon]\subset[-\pi t,\pi(1-t)]$ and $\int_{\varepsilon<|u|\leq\pi}\mathcal D_J^2\leq2\pi\sin^{-2}(\varepsilon/2)$. On $|u|\leq\varepsilon$, $q(t+u/\pi)=q(t)+O(\omega(\varepsilon))$, and $\int_{|u|\leq\varepsilon}\mathcal D_J^2=2\pi(2J+1)-O(\varepsilon^{-2})$ by \eqref{eq:dirichlet}. Collecting the bounds,
\[
\frac1J\int_0^1\widetilde D_J(y,t)^2q(y)\dd y=q(t)+O\big(\omega(\varepsilon)+M\varepsilon^{-2}J^{-1}+M\delta^{-2}J^{-1}+M\delta^{-1}J^{-1/2}\big),
\]
uniformly in $t\in[\delta,1-\delta]$ and $q\in\mathcal Q$. Taking $\varepsilon=J^{-1/4}$ concludes.
\end{proof}

\subsection{Proofs of the propositions}\label{app:props}

We write $f(y,\z_c)$ for $f(y,1,(\z_c,\x_d))$, $\bar f_h(y)=\int f(y,\z_c)\boldsymbol K_h(\z_c-\x_c)\dd\z_c$, and $\|f\|_\infty$ for the supremum of $f$ on $[0,\tau_0]\times\mathcal N$. In all the proofs, $i$ is large enough so that $h(i)$ is smaller than the radius of $\mathcal N$ and $b(i)<\min(t_0,\tau_0-t_1)$ when $t_0>0$; the finitely many first contributions are bounded and do not affect the conditions.

\begin{proof}[Proof of Proposition~\ref{prop:kernel}]
(\ref{prop:kernel:i}) Condition \eqref{cond:I}: $|\psi^{\mathrm K}_i|\leq\|K\|_\infty^{p}(c_bc_h^{d_c})^{-1}i^{\gamma p}$, and
\begin{align*}
\EE[\psi^{\mathrm K}_i(t;Z)^2]&=\int\!\!\int K_{b}^2(y-t)\boldsymbol K_h^2(\z_c-\x_c)f(y,\z_c)\dd y\dd\z_c\\
&=\frac{1}{bh^{d_c}}\int\!\!\int K^2(u)\boldsymbol K^2(\mathbf v)f(t+bu,\x_c+h\mathbf v)\dd u\dd\mathbf v
\end{align*}
(with $b=b(i)$, $h=h(i)$), which is at most $\|f\|_\infty c_K^{p}(bh^{d_c})^{-1}=Ci^{\gamma p}$; there is no remainder. Condition \eqref{cond:V}: by the uniform continuity of $f$, the last display equals $(bh^{d_c})^{-1}\{c_K^pf(t,\x_c)+o(1)\}$ uniformly in $t\in\mathcal T$, while $\EE[\psi^{\mathrm K}_i]^2\leq\|f\|_\infty^2$. Hence $\Var\{\xi_i(t)\}=\mathfrak c_Ki^{\gamma p}f(t,1,\x)\{1+o(1)\}$ uniformly in $t$, and summing, $V_n(t)=\mathfrak V_nv(t)\{1+o(1)\}$ uniformly, with $\mathfrak V_n=\sum_{i\leq n}i^{\gamma p}\sim n^{1+\gamma p}/(1+\gamma p)$ (the two normalizations differ by a factor $1+o(1)$ and either can be used). The function $v$ is continuous and positive on $\mathcal T$ by (\ref{ass:D2a}) and (\ref{ass:D3}), and $\sup_{\mathcal T}|\hat v_n/v-1|\to0$ almost surely follows from Condition \eqref{cond:E} and (\ref{ass:D3}). Condition \eqref{cond:B}: by the symmetry of $K$ and $\boldsymbol K$ and a second order Taylor expansion with the bounded second derivatives of (\ref{ass:D2a}), $\EE[\psi^{\mathrm K}_i(t;Z)]-f(t,\x_c)=\int\!\!\int K(u)\boldsymbol K(\mathbf v)\{f(t+bu,\x_c+h\mathbf v)-f(t,\x_c)\}\dd u\dd\mathbf v=O(b^2+h^2)=O(i^{-2\gamma})$ uniformly in $t$, and $n^{-1}\sum_{i\leq n}i^{-2\gamma}=O(n^{-2\gamma})$ since $2\gamma<1$. Condition \eqref{cond:E}: $K$ being Lipschitz, $|\psi^{\mathrm K}_i(t;\mathfrak z)-\psi^{\mathrm K}_i(s;\mathfrak z)|\leq L_Kb^{-2}\|\boldsymbol K\|_\infty h^{-d_c}|t-s|\leq Ci^{\gamma(p+1)}|t-s|$; Lemma~\ref{lemma:chaining}(\ref{lemma:chaining:a}) applies with $\kappa=\gamma(p+1)<(3+\gamma p)/2$, which holds since $\gamma<1/p$. Condition \eqref{cond:M}: the last bound is \eqref{eq:incr_size} with $a=\gamma p$ and $\vartheta_i=i^{-\gamma}$, for which $i^{a}/\vartheta_i=i^{\kappa}$. If $|t-s|\leq b$, the function $y\mapsto K_b(y-t)-K_b(y-s)$ is bounded by $L_Kb^{-2}|t-s|$ and supported on a set of length at most $4b$, so that $\EE[\{\psi^{\mathrm K}_i(t;Z)-\psi^{\mathrm K}_i(s;Z)\}^2]\leq Ch^{-d_c}b^{-3}|t-s|^2\leq Ci^{\gamma p}\vartheta_i^{-2}|t-s|^2$; if $|t-s|>b$, the bound on $\EE[\psi^{\mathrm K}_i(t;Z)^2]$ gives the same conclusion, that is \eqref{eq:incr_var}. Lemma~\ref{lemma:chaining}(\ref{lemma:chaining:d}) applies, since $\varsigma_n=o(n^{-\gamma})$.

(\ref{prop:kernel:ii}) The matrix of the linear system defining $(\mathsf k_0(\theta),\mathsf k_1(\theta))$ has determinant $\int_{-1}^\theta K\int_{-1}^\theta u^2K-(\int_{-1}^\theta uK)^2>0$ for $\theta\in[0,1]$, by the Cauchy--Schwarz inequality, the equality case being excluded since $K\dd u$ is not concentrated at one point of $[-1,0]$. Hence $\mathsf k_0,\mathsf k_1$ are continuous on $[0,1]$, Lipschitz if $K$ is, with $\mathsf k_0(1)=1$, $\mathsf k_1(1)=0$, so that $t\mapsto K^t_{B,b}(y-t)$ is continuous at $t=b$. It follows that $|K^t_{B,b}(s)|\leq C_BK_b(s)$ with $C_B=\|\mathsf k_0\|_\infty+\|\mathsf k_1\|_\infty$, and that $|K^t_{B,b}(y-t)-K^s_{B,b}(y-s)|\leq Cb^{-2}|t-s|$. Conditions \eqref{cond:I} and \eqref{cond:E} on $[0,t_1]$ then follow exactly as in (\ref{prop:kernel:i}). For Condition \eqref{cond:B}, for $t\geq b$ the computation of (\ref{prop:kernel:i}) applies; for $t\in[0,b)$, with $\theta=t/b$ and $u=(y-t)/b\in[-\theta,1]$,
\[
\int_0^{\tau_0}K^t_{B,b}(y-t)f(y,\z_c)\dd y=\int_{-\theta}^1\{\mathsf k_0(\theta)-\mathsf k_1(\theta)u\}K(u)f(t+bu,\z_c)\dd u=f(t,\z_c)+O(b^2),
\]
by the defining moment conditions and a second order expansion of $f(t+bu,\z_c)$ in $u$, uniformly in $t$ and $\z_c$; smoothing in $\z_c$ adds $O(h^2)$.

Integrated contributions. Let $r=\inf_{[0,t_1]}R(\cdot,\x)$ and $\Phi_b(y,t)=\int_0^tK^s_{B,b}(y-s)R(s,\x)^{-1}\dd s$, so that $\psi^\Lambda_i(t;\mathfrak z)=\1\{\iota=1\}\Phi_{b(i)}(y,t)\boldsymbol K_{h(i)}(\z_c-\x_c)\1\{\z_d=\x_d\}$. Since $\int_0^t|K^s_{B,b}(y-s)|\dd s\leq C_B\int K_b(y-s)\dd s\leq C_B$, we have $|\Phi_b|\leq C_B/r$, hence $|\psi^\Lambda_i|\leq Ch(i)^{-d_c}=Ci^{\gamma d_c}$ and $\Var\{\xi^\Lambda_i(t)\}\leq\EE[\psi^\Lambda_i(t;Z)^2]\leq Ci^{\gamma d_c}$ uniformly in $t\in[0,t_1]$. For the variance profile, we claim that
\begin{equation}\label{eq:Phi_b}
\sup_{t\in[0,t_1]}\sup_{\z_c\in\mathcal N}\Big|\int_0^{\tau_0}\Phi_b(y,t)^2f(y,\z_c)\dd y-\int_0^t\frac{f(y,\z_c)}{R(y,\x)^2}\dd y\Big|=O(b).
\end{equation}
Indeed, for $y\in[2b,t-b]$, the window $[y-b,y+b]$ is contained in $[b,t]$, where the boundary kernel is the standard kernel, so that $\Phi_b(y,t)=\int_{-1}^1K(u)R(y-bu,\x)^{-1}\dd u=R(y,\x)^{-1}+O(b)$, because $1/R(\cdot,\x)$ is Lipschitz; for $y>t+b$, $\Phi_b(y,t)=0=\1\{y\leq t\}$; and the remaining values of $y$ form a set of Lebesgue measure at most $4b$, on which both integrands are bounded. This proves \eqref{eq:Phi_b}. Consequently,
\begin{align*}
\EE[\psi^\Lambda_i(t;Z)^2]&=\frac{1}{h^{d_c}}\int\boldsymbol K(\mathbf v)^2\Big\{\int_0^t\frac{f(y,\x_c+h\mathbf v)}{R(y,\x)^2}\dd y+O(b)\Big\}\dd\mathbf v\\
&=\frac{1}{h^{d_c}}\Big\{c_K^{d_c}\int_0^t\frac{f(y,1,\x)}{R(y,\x)^2}\dd y+O(b+h)\Big\},
\end{align*}
uniformly in $t\in[0,t_1]$, using that $f$ is Lipschitz in $\z_c$, while $\EE[\psi^\Lambda_i]^2=O(1)$. Hence $\Var\{\xi^\Lambda_i(t)\}=\mathfrak c_\Lambda i^{\gamma d_c}\{v^\Lambda(t)/\mathfrak c_\Lambda+O(i^{-\gamma})\}$, and $V^\Lambda_n(t)=\mathfrak V^\Lambda_nv^\Lambda(t)\{1+o(1)\}$ uniformly in $t\in\mathcal T^\Lambda$, where $v^\Lambda$ is continuous and, by (\ref{ass:D3}), bounded below by $\mathfrak c_\Lambda\int_0^{t_0}f(y,1,\x)R(y,\x)^{-2}\dd y>0$ on $\mathcal T^\Lambda$. Finally, $\sup_{t\in\mathcal T^\Lambda}|\hat v^\Lambda_n(t)/v^\Lambda(t)-1|\to0$ almost surely because $\int_0^t\hat f_n/\hat R_n^2-\int_0^tf/R^2\to0$ uniformly in $t$, by Condition \eqref{cond:E} on $[0,t_1]$ and the uniform consistency of $\hat R_n$, and $v^\Lambda$ is bounded below.
\end{proof}

\begin{proof}[Proof of Proposition~\ref{prop:series}]
We use the properties \eqref{cond:S1}--\eqref{cond:S3} of the augmented basis, established in Lemma~\ref{lem:S3} and  \cite{PS2026}, and write $J=J(i)$, $h=h(i)$, $D_J=D_{J(i)}$. Recall that $\int_0^{\tau_0}D_J(y,t)^2\dd y=\sum_{j\leq J}\varphi_j(t)^2\leq CJ$ by orthonormality and \eqref{cond:S1}, and that $\int_0^{\tau_0}D_J(y,t)q(y)\dd y=S_Jq(t)$.

(\ref{prop:series:i}) Condition \eqref{cond:I}: $|D_J(y,t)|\leq CJ$ by \eqref{cond:S1}, so $|\psi^{\mathrm S}_i|\leq CJh^{-d_c}\leq Ci^{\gamma p}$, and
\[
\EE[\psi^{\mathrm S}_i(t;Z)^2]=\int\boldsymbol K_h(\z_c-\x_c)^2\int_0^{\tau_0}D_J(y,t)^2f(y,\z_c)\dd y\dd\z_c\leq\|f\|_\infty\frac{c_K^{d_c}}{h^{d_c}}\sum_{j\leq J}\varphi_j(t)^2\leq\frac{CJ}{h^{d_c}}.
\]
Condition \eqref{cond:V}: the family $\{f(\cdot,\z_c):\z_c\in\mathcal N\}$ is uniformly bounded and uniformly equicontinuous on $[0,\tau_0]$ by (\ref{ass:D2b}), so that \eqref{cond:S3}, that is Lemma~\ref{lem:S3}(\ref{lem:S3:iv}) after scaling, applied with $\mathcal T=[t_0,t_1]\subset(0,\tau_0)$, gives $\int_0^{\tau_0}D_J(y,t)^2f(y,\x_c+h\mathbf v)\dd y=J\{\tau_0^{-1}f(t,\x_c+h\mathbf v)+o(1)\}$ uniformly in $t\in\mathcal T$ and $\mathbf v\in[-1,1]^{d_c}$. Since $f$ is Lipschitz in $\z_c$, $\EE[\psi^{\mathrm S}_i(t;Z)^2]=Jh^{-d_c}\{c_K^{d_c}\tau_0^{-1}f(t,1,\x)+o(1)\}$ uniformly in $t$, while $\EE[\psi^{\mathrm S}_i(t;Z)]=\int\boldsymbol K_h(\z_c-\x_c)S_J\{f(\cdot,\z_c)\}(t)\dd\z_c$ is bounded by \eqref{cond:S2} and (\ref{ass:D2b}). Hence $\Var\{\xi_i(t)\}=J(i)i^{\gamma d_c}c_K^{d_c}(\tau_0c_h^{d_c})^{-1}f(t,1,\x)\{1+o(1)\}$ and $V_n(t)=\mathfrak V_nv(t)\{1+o(1)\}$ uniformly in $t\in\mathcal T$, with $\mathfrak V_n=\sum_{i\leq n}J(i)i^{\gamma d_c}$, which is of exact order $n^{1+\gamma p}$ since $J(i)\geq c_Ji^{\gamma}/2$ for $i$ large; $v$ is continuous and positive by (\ref{ass:D3}), and $\hat v_n/v\to1$ uniformly by Condition \eqref{cond:E}. Condition \eqref{cond:B}: by (\ref{cond:S2}), applied to $q=f(\cdot,\z_c)$, whose second derivative in $y$ is Lipschitz uniformly in $\z_c\in\mathcal N$ by (\ref{ass:D2b}), $S_J\{f(\cdot,\z_c)\}(t)=f(t,\z_c)+O(J^{-2})$ uniformly in $t\in[0,\tau_0]$ and $\z_c$; smoothing in $\z_c$ adds $O(h^2)$, as in the proof of Proposition~\ref{prop:kernel}. Hence $\EE[\psi^{\mathrm S}_i(t;Z)]-f(t,1,\x)=O(J(i)^{-2}+h(i)^2)=O(i^{-2\gamma})$ uniformly in $t\in[0,\tau_0]$, which gives \eqref{cond:B} with $\mathfrak b=2\gamma$. Condition \eqref{cond:E}: by (\ref{cond:S1}), $|\partial_tD_J(y,t)|\leq\sum_{j\leq J}\|\varphi_j\|_\infty\|\varphi_j'\|_\infty\leq CJ^2$, so that $|\psi^{\mathrm S}_i(t;\mathfrak z)-\psi^{\mathrm S}_i(s;\mathfrak z)|\leq CJ^2h^{-d_c}|t-s|\leq Ci^{\gamma(p+1)}|t-s|$ and Lemma~\ref{lemma:chaining}(\ref{lemma:chaining:a}) applies as in Proposition~\ref{prop:kernel}. Condition (\ref{cond:M}): the last bound is \eqref{eq:incr_size} with $a=\gamma p$ and $\vartheta_i=1/J(i)$, for which $i^{a}/\vartheta_i\leq Ci^{\kappa}$ with $\kappa=\gamma(p+1)$, and, by orthonormality and (\ref{cond:S1}), $\int_0^{\tau_0}\{D_J(y,t)-D_J(y,s)\}^2\dd y=\sum_{j\leq J}\{\varphi_j(t)-\varphi_j(s)\}^2\leq CJ^3|t-s|^2$, so that $\EE[\{\psi^{\mathrm S}_i(t;Z)-\psi^{\mathrm S}_i(s;Z)\}^2]\leq CJ^3h^{-d_c}|t-s|^2\leq Ci^{\gamma p}\vartheta_i^{-2}|t-s|^2$, that is \eqref{eq:incr_var}; Lemma~\ref{lemma:chaining}(\ref{lemma:chaining:d}) applies, since $\varsigma_n=o(n^{-\gamma})$.

(\ref{prop:series:ii}) Conditions \eqref{cond:I}, \eqref{cond:B} and \eqref{cond:E} on $[0,t_1]$ have just been proved, since none of the arguments above except \eqref{cond:V} used $t_0>0$. For the integrated contributions, let $g_t=\1_{[0,t]}R(\cdot,\x)^{-1}$, a function of bounded variation on $[0,\tau_0]$ with $\|g_t\|_\infty+\mathrm{TV}(g_t)\leq C$ uniformly in $t\in[0,t_1]$, because $1/R(\cdot,\x)$ is Lipschitz and bounded. Then
\[
\psi^\Lambda_i(t;\mathfrak z)=\1\{\iota=1,\,y\leq\tau_0\}\,S_J g_t(y)\,\boldsymbol K_h(\z_c-\x_c)\1\{\z_d=\x_d\},\qquad S_Jg_t(y)=\int_0^tD_J(y,s)R(s,\x)^{-1}\dd s.
\]
By Lemma~\ref{lem:S3}(\ref{lem:S3:ii}), $|S_Jg_t(y)|\leq C$ uniformly in $J$, $y$ and $t$, so that $|\psi^\Lambda_i|\leq Ch^{-d_c}=Ci^{\gamma d_c}$ and $\Var\{\xi^\Lambda_i(t)\}\leq Ci^{\gamma d_c}$. Moreover, for $q=f(\cdot,\z_c)$,
\[
\Big|\int_0^{\tau_0}\{(S_Jg_t)^2-g_t^2\}q\Big|\leq\|q\|_\infty\|S_Jg_t-g_t\|_2\{\|S_Jg_t\|_2+\|g_t\|_2\}\leq CJ^{-1/2},
\]
uniformly in $t$ and $\z_c$, by Lemma~\ref{lem:S3}(\ref{lem:S3:ii}), and $\int_0^{\tau_0}g_t^2q=\int_0^tq(y)R(y,\x)^{-2}\dd y$. The rest of the argument is identical to the one for the kernel learner, with $O(J^{-1/2})$ in place of $O(b)$ in \eqref{eq:Phi_b}: $\EE[\psi^\Lambda_i(t;Z)^2]=h^{-d_c}\{c_K^{d_c}\int_0^tf(y,1,\x)R(y,\x)^{-2}\dd y+O(h+J^{-1/2})\}$ uniformly in $t\in[0,t_1]$, which gives Assumption~\ref{ass:cumhaz}(\ref{ass:cumhaz:iii}) with the same $\mathfrak V^\Lambda_n$, $v^\Lambda$ and $\hat v^\Lambda_n$ as for the kernel learner.
\end{proof}

\begin{proof}[Proof of Proposition~\ref{prop:histogram}]
Write $k=k_i$, let $A=A_k(t)$ be the cell of $\mathcal C_k$ containing $(t,\x_c)$, of volume $|A|=\tau_02^{-kp}$ and diameter at most $C2^{-k}$, and note that $2^{-k}\leq2i^{-\gamma}$ and $2^{kp}\leq i^{\gamma p}$.

(\ref{prop:histogram:i}) Condition \eqref{cond:I}: $|\psi^{\mathrm H}_i|\leq|A|^{-1}\leq\tau_0^{-1}i^{\gamma p}$ and $\EE[\psi^{\mathrm H}_i(t;Z)^2]=|A|^{-2}\int_Af\leq\|f\|_\infty|A|^{-1}$. Condition \eqref{cond:B}: since $f$ is Lipschitz, $\EE[\psi^{\mathrm H}_i(t;Z)]=|A|^{-1}\int_Af(y,\z_c)\dd y\dd\z_c=f(t,1,\x)+O(2^{-k})=f(t,1,\x)+O(i^{-\gamma})$ uniformly in $t\in[0,t_1]$, and $n^{-1}\sum_{i\leq n}i^{-\gamma}=O(n^{-\gamma})$. Condition \eqref{cond:V}: $\EE[\psi^{\mathrm H}_i(t;Z)^2]=|A|^{-1}\{f(t,1,\x)+O(2^{-k})\}$ and $\EE[\psi^{\mathrm H}_i]^2=O(1)$, so that $\Var\{\xi_i(t)\}=\tau_0^{-1}2^{k_ip}f(t,1,\x)\{1+O(i^{-\gamma})\}$ uniformly in $t$; summing over $i\leq n$ gives $V_n(t)=\mathfrak V_nv(t)\{1+o(1)\}$ with $\mathfrak V_n=\tau_0^{-1}\sum_{i\leq n}2^{k_ip}$, which satisfies $\tau_0^{-1}2^{-p}\sum_{i\le n}i^{\gamma p}\leq\mathfrak V_n\leq\tau_0^{-1}\sum_{i\leq n}i^{\gamma p}$; the function $v=f(\cdot,1,\x)$ is continuous and positive on $\mathcal T$ by (\ref{ass:D3}), and $\hat v_n/v\to1$ uniformly by Condition \eqref{cond:E}. Conditions \eqref{cond:E} and (\ref{cond:M}): for $i\leq n$, $k_i\leq k_n$ and the partitions are nested, so that $\hat f^{\mathrm H}_n$ and $\bar f^{\mathrm H}_n$ are constant on the time-intervals of $\mathcal C_{k_n}$ intersected with $\mathcal T$, of which there are at most $2^{k_n}+1\leq n^{\gamma}+1$, with endpoints in $G_{k_n}\subset G_{\ell_n}$ since $\ell_n\geq k_n$; Lemma~\ref{lemma:chaining}(\ref{lemma:chaining:c}) applies.

(\ref{prop:histogram:ii}) Let $B$ be the cell of the partition of $[0,1)^{d_c}$ into cubes of side $2^{-k}$ containing $\x_c$, of volume $2^{-kd_c}$, and $I(y)$ the dyadic interval of length $\tau_02^{-k}$ containing $y$. Then
\[
\psi^\Lambda_i(t;\mathfrak z)=\1\{\iota=1\}\,2^{kd_c}\1\{\z_c\in B\}\1\{\z_d=\x_d\}\,\Theta_k(y,t),\qquad \Theta_k(y,t)=\frac{2^k}{\tau_0}\int_{I(y)\cap[0,t]}\frac{\dd s}{R(s,\x)}.
\]
We have $0\leq\Theta_k\leq1/r$, $\Theta_k(y,t)=R(y,\x)^{-1}+O(2^{-k})$ if $I(y)\subset[0,t]$, and $\Theta_k(y,t)=0$ if $I(y)\cap[0,t]=\emptyset$; hence $\Theta_k(y,t)=\1\{y\leq t\}R(y,\x)^{-1}+O(2^{-k})$ for all $y$ outside the interval $I(t)$, of length $\tau_02^{-k}$. It follows that $|\psi^\Lambda_i|\leq r^{-1}2^{k_id_c}\leq r^{-1}i^{\gamma d_c}$ and, since $f$ is Lipschitz,
\[
\EE[\psi^\Lambda_i(t;Z)^2]=2^{2kd_c}\int_B\int_0^{\tau_0}\Theta_k(y,t)^2f(y,\z_c)\dd y\dd\z_c=2^{kd_c}\Big\{\int_0^t\frac{f(y,1,\x)}{R(y,\x)^2}\dd y+O(2^{-k})\Big\},
\]
uniformly in $t\in[0,t_1]$, while $\EE[\psi^\Lambda_i]^2=O(1)$. The conclusion follows as in the proof of Proposition~\ref{prop:kernel}(\ref{prop:kernel:ii}), with $\mathfrak V^\Lambda_n=\sum_{i\leq n}2^{k_id_c}$, of exact order $n^{1+\gamma d_c}$.
\end{proof}

\begin{proof}[Proof of Proposition~\ref{prop:denominator}] 
Let $u_m=u_m^+-u_m^-$ be as in (\ref{ass:D4}), and let $C_u$ bound $|u_m^+|+|u_m^-|$ on $[0,\tau_0]\times\mathcal Z$. Write $g=g(i)$ and $\mathfrak l_i=\log(i+1)$. Size and variance: $|\psi^R_{m,i}|\leq C_u\|\boldsymbol K\|_\infty g^{-d_c}\leq Ci^{\beta d_c}\mathfrak l_i^{-\nu d_c}$ and $\EE[\psi^R_{m,i}(t;Z)^2]\leq C_u^2\int\boldsymbol K_g(\z_c-\x_c)^2f_\X(\z_c,\x_d)\dd\z_c\leq C_u^2c_K^{d_c}g^{-d_c}\sup_{\mathcal N}f_\X(\cdot,\x_d)$, which is $O(i^{\beta d_c}\mathfrak l_i^{-\nu d_c})$ since $f_\X(\cdot,\x_d)$ is bounded on $\mathcal N$ by (\ref{ass:D4}). Bias: $\EE[\psi^R_{m,i}(t;Z)]=\int R_m(t,(\z_c,\x_d))\boldsymbol K_g(\z_c-\x_c)\dd\z_c=R_m(t,\x)+O(g^2)$ uniformly in $t\in[0,\tau_0]$, by the symmetry of $\boldsymbol K$ and the bounded second derivatives in $\z_c$ of (\ref{ass:D4}). Hence the bias of $\hat R_{m,n}$ is of order $n^{-1}\sum_{i\leq n}g(i)^2=O(n^{-2\beta}\mathfrak l_n^{2\nu}+n^{-1}\mathfrak l_n^{1+2\nu})$, the second term being needed only for $\beta\geq1/2$ and being of smaller order than the second term of $r^R_n$. Stochastic term: with $\psi^{(1)}_i=u^+_m\boldsymbol K_g(\z_c-\x_c)\1\{\z_d=\x_d\}$ and $\psi^{(2)}_i=u^-_m\boldsymbol K_g(\z_c-\x_c)\1\{\z_d=\x_d\}$, both non-decreasing in $t$ since $\boldsymbol K\geq0$, with the size and variance bounds above, and with $t\mapsto\EE[\psi^{(j)}_i(t;Z)]=\int\EE[u_m^{\pm}(t;Z)\mid\X=(\z_c,\x_d)]f_\X(\z_c,\x_d)\boldsymbol K_g(\z_c-\x_c)\dd\z_c$ Lipschitz on $[0,\tau_0]$ with a constant not depending on $i$ by (\ref{ass:D4}), Lemma~\ref{lemma:chaining}(\ref{lemma:chaining:b}) applies for $\nu=0$, with $\kappa=0$ and $a=\beta d_c<1$. For $\nu>0$, the same proof applies,  with the bound $Cn^{\beta d_c}\mathfrak l_n^{-\nu d_c}$ on the size, $Cn^{1+\beta d_c}\mathfrak l_n^{-\nu d_c}$ on the sum of the variances and $v_n=C_0\{n^{1+\beta d_c}\mathfrak l_n^{1-\nu d_c}\}^{1/2}$ in \eqref{eq:bernstein}, and gives the second term of $r^R_n$.
\end{proof}

\begin{proof}[Proof of Corollary~\ref{cor:exponents}]
Item (\ref{cor:exponents:i}). By Proposition~\ref{prop:kernel}(\ref{prop:kernel:i}) or \ref{prop:series}(\ref{prop:series:i}), Assumption~\ref{ass:hazard}(\ref{ass:hazard:i}) holds with $a=\gamma p$, $\mathfrak b=2\gamma$, and $2\mathfrak b\geq1-a$ is equivalent to $\gamma\geq1/(p+4)=1/(d_c+5)$. Assumption~\ref{ass:hazard}(\ref{ass:hazard:ii}) is (\ref{ass:D4}). For Assumption~\ref{ass:hazard}(\ref{ass:hazard:iii}), Proposition~\ref{prop:denominator} gives the rate $r^R_n$, which is $o(w_n(\gamma p))$ as soon as $\beta\geq(1-\gamma p)/4$ and $\beta d_c<\gamma p$; the case $2\beta=(1-\gamma p)/2$ is admissible because $2\nu<1/2$. The constraint $\beta<1/d_c$ of Section~\ref{subsec:denominator} is implied by $\beta<\gamma p/d_c<1/d_c$. Theorem~\ref{theo:hazard} applies.

Item (\ref{cor:exponents:ii}). By Proposition~\ref{prop:kernel}(\ref{prop:kernel:ii}) or \ref{prop:series}(\ref{prop:series:ii}), Assumption~\ref{ass:cumhaz}(\ref{ass:cumhaz:i})--(\ref{ass:cumhaz:iii}) hold with $a_\Lambda=\gamma d_c$, the remainder $\mathfrak r_n$ being zero. In Assumption~\ref{ass:cumhaz}(\ref{ass:cumhaz:v}), $2\mathfrak b\geq1-a_\Lambda$ is equivalent to $\gamma\geq1/(d_c+4)$, and the condition on the grid holds since $a-a_\Lambda=\gamma$. Assumption~\ref{ass:cumhaz}(\ref{ass:cumhaz:iv}) requires $r^R_n=o(w_n(\gamma d_c))$. The first term of $r^R_n$ is $o(w_n(\gamma d_c))$ if and only if $\beta\geq(1-\gamma d_c)/4$, the equality being admissible because $2\nu<1/2$. The second term is $o(w_n(\gamma d_c))$ if and only if $\beta<\gamma$, or $\beta=\gamma$ and $\nu>0$. Theorem~\ref{theo:cumhaz} and Corollary~\ref{cor:survival} apply.

Item (\ref{cor:exponents:iii}). For the histogram, Assumption~\ref{ass:hazard}(\ref{ass:hazard:i}) holds by Proposition~\ref{prop:histogram}(\ref{prop:histogram:i}), Assumption~\ref{ass:cumhaz}(\ref{ass:cumhaz:i}) and (\ref{ass:cumhaz:iii}) hold by Proposition~\ref{prop:histogram}(\ref{prop:histogram:ii}), and Assumption~\ref{ass:cumhaz}(\ref{ass:cumhaz:ii}) follows from (\ref{ass:D4});  $\mathfrak b=\gamma$ and the conditions $2\mathfrak b\geq1-\gamma p$ and $2\mathfrak b\geq1-\gamma d_c$ read $\gamma\geq1/(d_c+3)$ and $\gamma\geq1/(d_c+2)$; the condition on the grid holds since $\ell_n\geq k_n$ gives $\varsigma_n\leq2\tau_0n^{-\gamma}$.
\end{proof}

\subsection{Other one-pass averages for the at-risk components}\label{app:denominator_alt}

The kernel in \eqref{eq:psi_R} can be replaced by other weights in $\x_c$. The covariate value $\x$ is fixed, and the denominators only enter the results through their uniform rate, without any condition of the type \eqref{cond:M}; the weights therefore need not come from nested partitions. Let $\delta(i)=c_\delta i^{-\beta}\{\log(i+1)\}^{\nu}$, with $\beta\in(0,1/d_c)$, $\nu\geq0$ and $c_\delta>0$, be the resolution of the $i$-th contribution, and let
\begin{equation}\label{eq:psi_R_W}
\psi^{W}_{m,i}(t;\mathfrak z)=u_m(t;\mathfrak z)\,W_i(\z_c)\,\1\{\z_d=\x_d\},\qquad \hat R^{W}_{m,n}(t)=\frac1n\sum_{i=1}^n\psi^{W}_{m,i}(t;Z_i),
\end{equation}
where the $W_i:\mathcal X_c\to\RR$ are bounded measurable functions, fixed in advance and depending on $\x_c$, so that $\hat R^W_{m,n}$ is a one-pass average. We write $|\cdot|_\infty$ for the maximum norm on $\RR^{d_c}$. The weights are said to be \emph{local of order $\mathfrak q\in\{1,2\}$} if there exists a constant $C_W$ such that, for all $i$ large enough,
\begin{enumerate}
\item[(W1)]\condlabel{W1}{cond:W1} $W_i(\z_c)=0$ if $|\z_c-\x_c|_\infty>C_W\delta(i)$, and $\sup_{\z_c}|W_i(\z_c)|\leq C_W\delta(i)^{-d_c}$;
\item[(W2)]\condlabel{W2}{cond:W2} $\int W_i(\z_c)\dd\z_c=1$ and, if $\mathfrak q=2$, $\int(\z_c-\x_c)W_i(\z_c)\dd\z_c=0$.
\end{enumerate}
No sign condition is imposed. The kernel weights $\boldsymbol K_{g(i)}(\cdot-\x_c)$ of \eqref{eq:psi_R} are local of order $2$, with $\delta(i)=g(i)$. For the next three examples, $\mathcal X_c\subset[0,1)^{d_c}$, as in Section~\ref{subsubsec:histogram}, and $k=k^R_i=\max\{0,\lceil\log_2(1/\delta(i))\rceil\}$, so that $2^{-k}\in(\delta(i)/2,\delta(i)]$ for $i$ large.

\emph{Histogram.} Let $W^{\mathrm H}_i=2^{kd_c}\1_{B_i}$, where $B_i$ is the cell containing $\x_c$ of the partition of $[0,1)^{d_c}$ into dyadic cubes of side $2^{-k}$. These weights are local of order $1$. Their first moment $\int(\z_c-\x_c)W^{\mathrm H}_i(\z_c)\dd\z_c$ is the offset between the center of $B_i$ and $\x_c$. It is of order $\delta(i)$, and it cannot vanish at two consecutive resolutions, since the center of a dyadic cube is a vertex of its sub-cubes.

\emph{Centered cell.} Let $W^{\mathrm C}_i=\{2\delta(i)\}^{-d_c}\1\{|\z_c-\x_c|_\infty\leq\delta(i)\}$. These weights are local of order $2$; this is \eqref{eq:psi_R} with the uniform kernel and $g(i)=\delta(i)$, and the proof of Proposition~\ref{prop:denominator} does not use the Lipschitz continuity of $K$.

\emph{Frequency polygon.} Let $\mathfrak c$ be the center of a dyadic cube of side $2^{-k}$ such that $\boldsymbol\theta=2^k(\x_c-\mathfrak c)\in[0,1)^{d_c}$. For $\mathbf j\in\{0,1\}^{d_c}$, let $B^{\mathbf j}_i$ be the cube with center $\mathfrak c+2^{-k}\mathbf j$ and $\lambda_{\mathbf j}=\prod_{r=1}^{d_c}\theta_r^{j_r}(1-\theta_r)^{1-j_r}$, and let $W^{\mathrm P}_i=2^{kd_c}\sum_{\mathbf j}\lambda_{\mathbf j}\1_{B^{\mathbf j}_i}$. The corresponding contribution is the value at $\x_c$ of the multilinear interpolation, between the centers of the cells, of the histogram contribution of resolution $k$. 
These weights are nonnegative, with $\sum_{\mathbf j}\lambda_{\mathbf j}=1$ and $\sum_{\mathbf j}\lambda_{\mathbf j}(\mathfrak c+2^{-k}\mathbf j)=\x_c$, so that they are local of order $2$. They require the $2^{d_c}$ cell sums around $\x_c$ at each resolution.

\emph{Cosine series.} Let $\mathcal X_c=[0,1]^{d_c}$, after an affine transformation, $J_R(i)=\lfloor1/\delta(i)\rfloor\vee1$ and
\[
W^{\mathrm S}_i(\z_c)=\prod_{r=1}^{d_c}D^{\cos}_{J_R(i)}(z_{c,r},x_{c,r}),
\]
where $D^{\cos}_J$ is the kernel of the cosine projection $S^{\cos}_J$ of Appendix~\ref{app:basis}, with $\tau_0=1$. The corresponding estimator $\hat R^{\mathrm S}_{m,n}(t)$ is the value at $\x_c$ of a recursive tensor-product cosine series estimator, the analogue in $\x_c$ of the learner of Section~\ref{subsubsec:series}. These weights are signed and they are not local: every observation with $\X_{i,d}=\x_d$ contributes. The requirements of Assumption~\ref{ass:data}(\ref{ass:D4}), which concern the neighborhood $\mathcal N$ of $\x_c$, then have to hold on the whole of $\mathcal X_c$:
\begin{enumerate}
\item[(D4$'$)]\condlabel{D4$'$}{ass:D4prime} $f_\X(\cdot,\x_d)$ is bounded on $\mathcal X_c$; $t\mapsto\EE[u_m^{\pm}(t;Z)\mid\X=(\z_c,\x_d)]f_\X(\z_c,\x_d)$, with $u_m^\pm$ as in (\ref{ass:D4}), is Lipschitz on $[0,\tau_0]$, uniformly in $\z_c\in\mathcal X_c$; and, for each coordinate $r$, the function $z_{c,r}\mapsto R_m(t,(\z_c,\x_d))$ has on $[0,1]$ a bounded first derivative and a Lipschitz second derivative, with bounds which do not depend on $t\in[0,\tau_0]$ and on the other coordinates of $\z_c$.
\end{enumerate}

\begin{proposition}[At-risk components, other weights]\label{prop:denominator_alt}
Let Assumption~\ref{ass:data}(\ref{ass:D1}) hold, and let $m\in\{1,2\}$.
\begin{enumerate}
\item[(i)]\condlabel{i}{prop:denominator_alt:i} (Local weights) Under Assumption~\ref{ass:data}(\ref{ass:D4}), if the weights are local of order $\mathfrak q\in\{1,2\}$, then
\begin{gather*}
\sup_{t\in[0,\tau_0]}|\hat R^{W}_{m,n}(t)-R_m(t,\x)|=O\big(r^{W}_n(\mathfrak q)\big)\quad\text{almost surely},\\
r^{W}_n(\mathfrak q)=\frac{\{\log(n+1)\}^{\mathfrak q\nu}}{n^{\mathfrak q\beta}}+\frac{\{\log(n+1)\}^{(1-\nu d_c)/2}}{n^{(1-\beta d_c)/2}}.
\end{gather*}
For $\mathfrak q=1$, the requirement on the second order partial derivatives in (\ref{ass:D4}) can be replaced by the Lipschitz continuity of $\z_c\mapsto R_m(t,(\z_c,\x_d))$ on $\mathcal N$, uniformly in $t$.
\item[(ii)]\condlabel{ii}{prop:denominator_alt:ii} (Cosine series) Under (\ref{ass:D4prime}), with $\x_c$ in the interior of $\mathcal X_c=[0,1]^{d_c}$,
\begin{gather*}
\sup_{t\in[0,\tau_0]}|\hat R^{\mathrm S}_{m,n}(t)-R_m(t,\x)|=O\big(r^{\mathrm S}_n\big)\quad\text{almost surely},\\
r^{\mathrm S}_n=\frac{\{\log(n+1)\}^{2\nu+d_c-1}}{n^{2\beta}}+\frac{\{\log(n+1)\}^{(1-\nu d_c)/2}}{n^{(1-\beta d_c)/2}}.
\end{gather*}
\end{enumerate}
\end{proposition}

The rate $r^W_n(2)$ is the rate $r^R_n$ of Proposition~\ref{prop:denominator}, with $\delta(i)$ in place of $g(i)$, and $r^{\mathrm S}_n=r^R_n$ when $d_c=1$.

\begin{proof}
(\ref{prop:denominator_alt:i}) Let $i$ be large enough for (\ref{cond:W1}) and (\ref{cond:W2}) to hold and for the set $\{\z_c:|\z_c-\x_c|_\infty\leq C_W\delta(i)\}$ to be contained in $\mathcal N$; the finitely many other contributions are bounded and change the estimator by $O(1/n)$. Write $\delta=\delta(i)$, let $W_i=W_i^+-W_i^-$ be the decomposition of $W_i$ into its positive and negative parts, and let $u_m=u_m^+-u_m^-$ and $C_u$ be as in the proof of Proposition~\ref{prop:denominator}. Then
\[
u_mW_i=\big(u_m^+W_i^++u_m^-W_i^-\big)-\big(u_m^+W_i^-+u_m^-W_i^+\big),
\]
where both terms are non-decreasing in $t$. By (\ref{cond:W1}), each of them is bounded by $2C_uC_W\delta^{-d_c}$, the second moment of the corresponding contribution is at most $4C_u^2C_W^2\delta^{-2d_c}(2C_W\delta)^{d_c}\sup_{\mathcal N}f_\X(\cdot,\x_d)\leq C\delta^{-d_c}$ and, by (\ref{ass:D4}), its expectation is Lipschitz in $t$ with a constant proportional to $\int|W_i|\leq C_W(2C_W)^{d_c}$. These are the bounds used for the stochastic term in the proof of Proposition~\ref{prop:denominator}, with $\delta(i)$ in place of $g(i)$, and the argument given there, which is valid for every $\nu\geq0$, yields the second term of $r^W_n(\mathfrak q)$. For the bias, let $\rho_t(\z_c)=R_m(t,(\z_c,\x_d))$. By (\ref{cond:W2}),
\[
\EE[\psi^W_{m,i}(t;Z)]-R_m(t,\x)=\int\{\rho_t(\z_c)-\rho_t(\x_c)\}W_i(\z_c)\dd\z_c,
\]
and, if $\mathfrak q=2$, the same identity holds with $\rho_t(\z_c)-\rho_t(\x_c)-\nabla\rho_t(\x_c)^\top(\z_c-\x_c)$ in the integral. On the support of $W_i$, these differences are $O(\delta)$ if $\rho_t$ is Lipschitz and $O(\delta^2)$ if its second order partial derivatives are bounded, uniformly in $t$; under (\ref{ass:D4}) both properties hold, since the functions $\rho_t$ are uniformly bounded on the ball $\mathcal N$, with uniformly bounded second order partial derivatives. Hence the bias of the $i$-th contribution is $O(\delta(i)^{\mathfrak q})$ uniformly in $t\in[0,\tau_0]$, and $n^{-1}\sum_{i\leq n}\delta(i)^{\mathfrak q}=O(n^{-\mathfrak q\beta}\{\log(n+1)\}^{\mathfrak q\nu})$ if $\mathfrak q\beta<1$. Otherwise $\mathfrak q=2$ and $\beta\geq1/2$, the average is $O(n^{-1}\{\log(n+1)\}^{1+2\nu})$, of smaller order than the second term of $r^W_n(2)$.

(\ref{prop:denominator_alt:ii}) Write $J=J_R(i)$. By \eqref{eq:dirichlet}, $|D^{\cos}_J|\leq2J+1$ and $\int_0^1D^{\cos}_J(z,x)^2\dd z\leq2J+1$, so that $|W^{\mathrm S}_i|\leq(2J+1)^{d_c}\leq C\delta(i)^{-d_c}$, $\int W^{\mathrm S}_i(\z_c)^2f_\X(\z_c,\x_d)\dd\z_c\leq C(2J+1)^{d_c}$ by (\ref{ass:D4prime}), and $\int|W^{\mathrm S}_i|\leq(2J+1)^{d_c/2}$ by the Cauchy--Schwarz inequality. The decomposition used for (\ref{prop:denominator_alt:i}) applies, and the expectations of its two terms are Lipschitz in $t$ with a constant $Ci^{\beta d_c/2}$, which is allowed in Lemma~\ref{lemma:chaining}(\ref{lemma:chaining:b}). This gives the second term of $r^{\mathrm S}_n$. For the bias, we use a bound at an interior point: if $q$ has on $[0,1]$ a first derivative bounded by $M_1$ and a second derivative which is Lipschitz with constant $L$, then, for $J\geq1$, $\mathsf d\in(0,1/2)$ and $x\in[\mathsf d,1-\mathsf d]$,
\begin{equation}\label{eq:cos_interior}
|S^{\cos}_Jq(x)-q(x)|\leq C\big(M_1\mathsf d^{-1}+L\big)J^{-2},
\end{equation}
with an absolute constant $C$. Indeed, as in the proofs of Lemma~\ref{lem:S3}(\ref{lem:S3:i}) and (\ref{lem:S3:iii}), two integrations by parts give $\langle q,c_j\rangle=\sqrt2\{(-1)^jq'(1)-q'(0)\}(\pi j)^{-2}+\epsilon_j$, with $|\epsilon_j|\leq\sqrt2L(\pi j)^{-3}$. The cosine series of $q$ thus converges absolutely, and $S^{\cos}_Jq(x)-q(x)=-\sum_{j>J}\langle q,c_j\rangle c_j(x)$. The terms $\epsilon_j$ contribute at most $CLJ^{-2}$ to this sum, and the other ones at most $CM_1\mathsf d^{-1}J^{-2}$, by the Abel summation bound of the proof of Lemma~\ref{lem:S3}(\ref{lem:S3:i}). Now let $\rho_t$ be as above and let $S^{(r)}$ denote the projection $S^{\cos}_J$ acting on the $r$-th coordinate. Then
\[
\EE[\psi^{\mathrm S}_{m,i}(t;Z)]-R_m(t,\x)=\big(S^{(1)}\cdots S^{(d_c)}\rho_t-\rho_t\big)(\x_c)=\sum_{r=1}^{d_c}\big(S^{(1)}\cdots S^{(r-1)}(S^{(r)}-\mathrm{Id})\rho_t\big)(\x_c).
\]
In the $r$-th term, the coordinates $r,\ldots,d_c$ are evaluated at $\x_c$. By \eqref{eq:cos_interior} and (\ref{ass:D4prime}), the function $(S^{(r)}-\mathrm{Id})\rho_t$, with its coordinates $r,\ldots,d_c$ fixed at those of $\x_c$, is bounded by $CJ^{-2}$, uniformly in the first $r-1$ coordinates and in $t$. Each of the remaining $r-1$ projections multiplies the supremum norm by at most the Lebesgue constant $\sup_{x}\int_0^1|D^{\cos}_J(z,x)|\dd z\leq C\log(J+2)$  \citep[Vol.~I, Chap.~II]{Z2002}. Hence the bias of the $i$-th contribution is $O(J^{-2}\{\log(J+2)\}^{d_c-1})=O(\delta(i)^2\{\log(i+1)\}^{d_c-1})$, uniformly in $t\in[0,\tau_0]$, and the first term of $r^{\mathrm S}_n$ follows by averaging, as in (\ref{prop:denominator_alt:i}).
\end{proof}

\begin{corollary}[Exponents, other weights]\label{cor:exponents_alt}
In Corollary~\ref{cor:exponents}, let $\hat R_{m,n}$ be given by \eqref{eq:psi_R_W} in place of \eqref{eq:psi_R}.
\begin{enumerate}
\item[(i)]\condlabel{i}{cor:exponents_alt:i} (Order two) If the weights are local of order $2$, as for the centered cell and the frequency polygon, Corollary~\ref{cor:exponents} holds without change, with $\nu\in[0,1/4)$.
\item[(ii)]\condlabel{ii}{cor:exponents_alt:ii} (Order one) If the weights are local of order $1$, as for the histogram, Corollary~\ref{cor:exponents} holds with the conditions on $(\beta,\nu)$ replaced by $(1-\gamma p)/2\leq\beta\leq\gamma p/d_c$ for the hazard, and by $(1-\gamma d_c)/2\leq\beta\leq\gamma$ for the cumulative hazard, with $\nu<1/2$ if the lower bound is attained and $\nu>0$ if the upper bound is attained. These conditions can be met if and only if $\gamma\geq d_c/\{(d_c+1)(d_c+2)\}$ for the hazard, and $\gamma\geq1/(d_c+2)$ for the cumulative hazard.
\item[(iii)]\condlabel{iii}{cor:exponents_alt:iii} (Cosine series) With the weights $W^{\mathrm S}_i$ and (\ref{ass:D4prime}) in place of Assumption~\ref{ass:data}(\ref{ass:D4}), Corollary~\ref{cor:exponents} holds without change if $d_c=1$. If $d_c\geq2$, it holds with strict lower bounds on $\beta$, that is $\beta>(1-\gamma p)/4$ for the hazard and $\beta>(1-\gamma d_c)/4$ for the cumulative hazard; the latter requires $\gamma>1/(d_c+4)$.
\end{enumerate}
\end{corollary}

\begin{proof}
The rates of Proposition~\ref{prop:denominator_alt} are compared with $w_n(a)$, for $a=\gamma p$ (hazard) and $a=\gamma d_c$ (cumulative hazard), as in the proof of Corollary~\ref{cor:exponents}, the rest of which is unchanged. Item (\ref{cor:exponents_alt:i}) holds because $r^W_n(2)$ is the rate $r^R_n$. For (\ref{cor:exponents_alt:ii}), the first term of $r^W_n(1)$ is $o(w_n(a))$ if and only if $\beta>(1-a)/2$, or $\beta=(1-a)/2$ and $\nu<1/2$; the second one if and only if $\beta d_c<a$, or $\beta d_c=a$ and $\nu>0$. The two bounds on $\beta$ are compatible if and only if $d_c(1-a)\leq2a$, that is $\gamma\geq d_c/\{(d_c+1)(d_c+2)\}$ when $a=\gamma p$, and $\gamma\geq1/(d_c+2)$ when $a=\gamma d_c$. For (\ref{cor:exponents_alt:iii}), the first term of $r^{\mathrm S}_n$ is $o(w_n(a))$ if and only if $2\beta>(1-a)/2$, or $2\beta=(1-a)/2$ and $2\nu+d_c-1<1/2$; the equality case is thus possible only if $d_c=1$, with $\nu<1/4$. For the cumulative hazard, $(1-\gamma d_c)/4<\beta\leq\gamma$ forces $\gamma>1/(d_c+4)$.
\end{proof}

\paragraph{Comparison with Corollary~\ref{cor:exponents}.} Table~\ref{tab:denominator_alt} summarizes the conditions. The weights of order two with a compact support are interchangeable with the kernel. The histogram is a learner of order one. With the histogram numerator of Corollary~\ref{cor:exponents}(\ref{cor:exponents:iii}), nothing is lost: the lower bounds $\gamma\geq1/(d_c+3)$ and $\gamma\geq1/(d_c+2)$ required there imply those of Corollary~\ref{cor:exponents_alt}(\ref{cor:exponents_alt:ii}), and $R_m$ only has to be Lipschitz in $\x_c$. With the kernel or series numerators, the rate-optimal exponent of the cumulative hazard, $\gamma=1/(d_c+4)$, is excluded in every dimension, and that of the hazard, $\gamma=1/(d_c+5)$, is excluded as soon as $d_c\geq2$; when $d_c=1$, it is admissible, since $d_c/\{(d_c+1)(d_c+2)\}=1/6=1/(d_c+5)$, with $\beta=1/3$ and $\nu\in(0,1/2)$. The smallest admissible $\gamma$ then gives the width $n^{-1/(d_c+2)}\sqrt{\log n}$ for both targets, instead of $n^{-2/(d_c+5)}\sqrt{\log n}$ and $n^{-2/(d_c+4)}\sqrt{\log n}$. This loss is due to the position of $\x_c$ in its cell, and not to binning: a cell centered at $\x_c$, or the interpolation between the cells, restores the second-order bias. The cosine series has a second-order bias at an interior point, although no adaptation to the boundary of $\mathcal X_c$ is made; the augmentation of Appendix~\ref{app:basis} would only be needed for a bias bound which is uniform up to the boundary. The series requires, however, the global condition (\ref{ass:D4prime}), and, when $d_c\geq2$, our bound for its bias carries the factor $\{\log n\}^{d_c-1}$ of the Lebesgue constants, which excludes the end of the range of exponents. We expect this factor to be removable under conditions on the mixed partial derivatives of $R_m$. For the modified current status data of Section~\ref{subsec:mcs}, the estimator \eqref{eq:hat_p} can be computed with the same weights $W_i$ in place of $\boldsymbol K_{g(i)}(\cdot-\x_c)$; it then has the rate of Proposition~\ref{prop:denominator_alt}, and the statements of Section~\ref{sec:mechanisms} hold under the conditions of Corollary~\ref{cor:exponents_alt}.

\begin{table}[ht]
\centering
\footnotesize
\setlength{\tabcolsep}{4pt}
\begin{tabular}{llll}
\hline
Weights in $\x_c$ & Hazard & Cumulative hazard & Smallest $\gamma$\\
\hline
Kernel, centered cell, polygon & $(1-\gamma p)/4\leq\beta<\gamma p/d_c$ & $(1-\gamma d_c)/4\leq\beta\leq\gamma$ & $1/(d_c+5)$;\ $1/(d_c+4)$\\
Histogram & $(1-\gamma p)/2\leq\beta\leq\gamma p/d_c$ & $(1-\gamma d_c)/2\leq\beta\leq\gamma$ & $d_c/\{(d_c+1)(d_c+2)\}$;\ $1/(d_c+2)$\\
Cosine series, $d_c=1$ & $(1-\gamma p)/4\leq\beta<\gamma p/d_c$ & $(1-\gamma d_c)/4\leq\beta\leq\gamma$ & $1/(d_c+5)$;\ $1/(d_c+4)$\\
Cosine series, $d_c\geq2$ & $(1-\gamma p)/4<\beta<\gamma p/d_c$ & $(1-\gamma d_c)/4<\beta\leq\gamma$ & $1/(d_c+5)$;\ $>1/(d_c+4)$\\
\hline
\end{tabular}
\caption{Conditions on $\beta$ for the one-pass   averages\ of the at-risk components, and smallest admissible $\gamma$ for the hazard and for the cumulative hazard with the kernel or series numerators. The equality cases carry the conditions on $\nu$ of Corollaries~\ref{cor:exponents} and \ref{cor:exponents_alt}.}\label{tab:denominator_alt}
\end{table}


\section{Complements and proofs for Section~\ref{sec:mechanisms}}\label{app:mechanisms}

\subsection{Left truncation and right censoring}\label{app:ltrc}

There is a latent vector $(T^*,C^*,L^*,\X^*)$ with $T^*\perp L^*\perp C^*\mid\X^*$, and a realization is observed if and only if $L^*\leq Y^*=T^*\wedge C^*$; the observed vector is $Z=(Y,I,L,\X)$, the conditional distribution of $(Y^*,I^*,L^*,\X^*)$ given $L^*\leq Y^*$, with $I^*=\1\{T^*\leq C^*\}$. Let $\mathfrak q(\x)=\PP(L^*\leq Y^*\mid\X^*=\x)$, assumed positive. Then, as in \citet[Section~6.4]{E2018},
\begin{align*}
f_{Y,I,\X}(t,1,\x)&=\frac{F_{L^*\mid\X^*}(t\mid\x)S_{C^*\mid\X^*}(t-\mid\x)}{\mathfrak q(\x)}f_{T^*\mid\X^*}(t\mid\x)f_\X(\x),\\
\PP(L\leq t\leq Y\mid\X=\x)&=\frac{F_{L^*\mid\X^*}(t\mid\x)S_{C^*\mid\X^*}(t-\mid\x)}{\mathfrak q(\x)}S_{T^*\mid\X^*}(t-\mid\x),
\end{align*}
where $f_\X$ is the density of the observed covariates, so that the hazard of $T^*$ given $\X^*=\x$ has the representation \eqref{eq:R_affine} with
\begin{equation}\label{eq:R_ltrc}
R(t,\x)=R_1(t,\x)=\PP(L\leq t\leq Y\mid\X=\x)f_\X(\x),\quad u_1(t;\mathfrak z)=\1\{\ell\leq t\}-\1\{y<t\},\quad \eta\equiv0,
\end{equation}
where $u_1(t;\mathfrak z)=\1\{\ell\leq t\leq y\}$, because $\ell\leq y$ for every observation. Both pieces are monotone in $t$, and (\ref{ass:mech:A1}) for $Y$ and $L$ gives the Lipschitz property in Assumption~\ref{ass:data}(\ref{ass:D4}). The at-risk function vanishes where $F_{L^*\mid\X^*}(t\mid\x)=0$, that is below the left end of the support of the truncation variable: the hazard is identified, and $R_1$ is bounded away from zero, on an interval $[s_0,t_1]$ with $F_{L^*\mid\X^*}(s_0\mid\x)>0$ and $S_{T^*\mid\X^*}(t_1\mid\x)S_{C^*\mid\X^*}(t_1\mid\x)>0$. With $s_0$ known and Assumption~\ref{ass:data} on $[s_0,\tau_0]$, the bands of Theorem~\ref{theo:hazard} apply on any $[t_0,t_1]$ with $t_0>s_0$. For the cumulative hazard, Theorem~\ref{theo:cumhaz} applies with the integral started at $s_0$ instead of $0$, see Remark~\ref{rem:s0} below: the targets are then $\Lambda_{T^*\mid\X^*}(t\mid\x)-\Lambda_{T^*\mid\X^*}(s_0\mid\x)$ and the conditional survival $S_{T^*\mid\X^*}(t\mid\x)/S_{T^*\mid\X^*}(s_0\mid\x)$ given $T^*>s_0$, which are the quantities identified under truncation  \citep{W1985}; when $L^*$ has an atom at $0$, or more generally when $R_1$ is bounded away from zero on $[0,t_1]$, one can take $s_0=0$.

\begin{remark}[Integration from a point other than the origin]\label{rem:s0}
Under truncation, the at-risk function may vanish at the origin and be small near it. We then fix a known $s_0$ with $\inf_{[s_0,t_1]}R(\cdot,\x)>0$. Under (\ref{ass:mech:A1}), this requires $F_{L^*\mid\X^*}(s_0\mid\x)>0$, so that observations with $Y<s_0$ occur. Theorem~\ref{theo:cumhaz} and Corollary~\ref{cor:survival} hold with the contributions $\psi_i(t;\mathfrak z)$ multiplied by $\1\{y\geq s_0\}$, with $\int_0^t$ replaced by $\int_{s_0}^t$ in the definitions of $\hat\Lambda_n$, $\psi^\Lambda_i$ and $v^\Lambda$, with $[0,t_1]$ replaced by $[s_0,t_1]$ in Assumption~\ref{ass:cumhaz} and with $t_0>s_0$; the targets are then $\Lambda_{T\mid\X}(t\mid\x)-\Lambda_{T\mid\X}(s_0\mid\x)$ and the conditional survival function given $T>s_0$. The truncated contributions estimate $f(t,1,\x)\1\{t\geq s_0\}$, which vanishes below $s_0$ as $f(\cdot,1,\x)$ does below $0$. The proofs are thus the same after a translation of the time origin to $s_0$: Assumption~\ref{ass:data} of Section~\ref{subsec:data_ass} is required on $[s_0,\tau_0]$, the boundary kernel of Proposition~\ref{prop:kernel} is used at $s_0$, that is $K^{t-s_0}_{B,b}(y-t)$ for $t\in[s_0,s_0+b)$, the basis of the series learner is built on $[s_0,\tau_0]$, and $s_0$ is an endpoint of the time-intervals of the histogram. Without the factor $\1\{y\geq s_0\}$, the boundary kernel would also weight the observations below $s_0$, with a bias of order $b$.
\end{remark}

\subsection{Complements and proofs for cure models}\label{app:cure_more}

\begin{remark}[Smoothness at $\tau$]\label{rem:cure_smooth}
Assumption~\ref{ass:data}(\ref{ass:D2}) is required on $[0,\tau_0]$, while under \eqref{eq:suff_followup} the sub-density $f(\cdot,1,\x)$ vanishes on $(\tau,\tau_0]$. Read literally, (\ref{ass:D2a}) and (\ref{ass:D2b}) require the density of the susceptibles and its first two derivatives to vanish at $\tau$. This excludes, for instance, a density with a jump at the end of its support. The restriction can be avoided. For the kernel and histogram learners, take $\tau_1>\tau$ as in Remark~\ref{rem:cure_tau}. The integrals $\hat\Lambda_n(\tau_1)$ and $\hat\varrho_n$ then smooth across $\tau$, and the kernel integrates to one, so that their bias keeps the order $b(i)^2$, or $2^{-k_i}$ for the histogram. On $\mathcal T$, the jump does not affect $\hat f_n$ as soon as $b(i)<\tau-t_1$. The proofs then use the bias of these integrals instead of the uniform bias of $\hat f_n$ near $\tau$. For the series learner, a jump inside $(0,\tau_0)$   degrades   the approximation on the whole interval. The basis is then built on $[0,\tau]$, that is $\tau_0=\tau$, in place of $\tau<\tau_0$, and (\ref{ass:D2b}) is only required on $[0,\tau]$, with one-sided derivatives at $\tau$. The details are omitted.
\end{remark}

\begin{remark}[Choice of $\tau$]\label{rem:cure_tau}
Any known $\tau_1\geq\tau$ with $\inf_{[0,\tau_1]}R_1(\cdot,\x)>0$ (or the infimum on $[s_0,\tau_1]$ under left truncation) can be used in place of $\tau$, since $\Lambda_{T\mid\X}(\cdot\mid\x)$ is constant on $[\tau,\tau_1]$; in practice $\tau_1$ is taken close to the end of the follow-up. Under left truncation, $\lambda_{T\mid\X}$ is identified only on $[s_0,\infty)$, with $F_{L^*\mid\X^*}(s_0\mid\x)>0$ as in Remark~\ref{rem:s0}. Representation \eqref{eq:cure_hazard} only involves $\lambda_{T\mid\X}$ on $[t,\tau]$. Theorem~\ref{theo:cure}(\ref{theo:cure:i}) thus holds with $[0,\tau]$ replaced by $[s_0,\tau]$ and $t_0>s_0$, since its proof only uses $[t_0,\tau]$. Representation \eqref{eq:cure_cumhaz} requires $s_0=0$, that is $F_{L^*\mid\X^*}(0\mid\x)>0$, as assumed in \eqref{eq:suff_followup}. Otherwise, $\mathfrak G$ applied to the integrals started at $s_0$, as in Remark~\ref{rem:s0}, gives $\Lambda_0(t\mid\x)-\Lambda_0(s_0\mid\x)$ and $S_0(t\mid\x)/S_0(s_0\mid\x)$, and Theorem~\ref{theo:cure}(\ref{theo:cure:ii}) holds for these targets, with the same proof.
\end{remark}

\begin{remark}[Cure threshold depending on the covariates]\label{rem:cure_tau_x}
In Section~\ref{subsec:cure} the threshold $\tau$ of \eqref{eq:suff_followup} does not depend on the covariates. A threshold $\tau(\x)$, as in \cite{PVK2020}, can be accommodated in the representation \eqref{eq:cure_hazard}, since only $\Lambda_{T\mid\X}(\tau(\x)\mid\x)$ enters. For the inference, $\tau(\x)$ has to be replaced by an estimator, and the additional error becomes part of the remainder; a sufficient condition is an estimator $\hat\tau_n(\x)$ such that $\Lambda_{T\mid\X}(\hat\tau_n(\x)\mid\x)-\Lambda_{T\mid\X}(\tau(\x)\mid\x)$ is $o(w_n(\gamma p))$ for the hazard and $o(w_n(\gamma d_c))$ for the cumulative hazard, which holds for instance when $\tau(\x)$ lies beyond the right end of the support of the susceptibles' distribution, so that the hazard of $T$ vanishes on a neighborhood of $\tau(\x)$, and $\hat\tau_n(\x)$ falls in this neighborhood eventually. We leave the detailed investigation for future work.
\end{remark}

\medskip

\begin{proof}[Proof of Theorem~\ref{theo:cure}]
Throughout, $\lambda=\lambda_{T\mid\X}(\cdot\mid\x)$, $\Lambda=\Lambda_{T\mid\X}(\cdot\mid\x)$ and $S=S_{T\mid\X}(\cdot\mid\x)$ refer to the improper variable $T$, whose hazard has the representation $f(\cdot,1,\x)/R_1(\cdot,\x)$ of the mechanism without cure. Let $r=\inf_{[0,\tau]}R_1(\cdot,\x)>0$ by \eqref{eq:suff_followup}. We write $\mathfrak b$ for the bias exponent of the learner, $\mathfrak b=2\gamma$ for the kernel and series learners and $\mathfrak b=\gamma$ for the histogram; these learners have $\mathfrak r_n\equiv0$. Since $S_0(t_1\mid\x)>0$, we have $S(t_1)>S(\tau)$, hence $\varrho(t)\geq\varrho(t_1)>0$ on $\mathcal T$, and $\varrho$ is Lipschitz on $[0,\tau]$ because $\lambda$ is bounded.

(\ref{theo:cure:i}) Write $\hat\lambda_{0,n}=\hat f_n/\hat R^{(0)}_n$ with $\hat R^{(0)}_n=\hat R_{1,n}(1-e^{-\hat\varrho_n})$ and $R^{(0)}=R_1(1-e^{-\varrho})$, which is Lipschitz and bounded below by $r(1-e^{-\varrho(t_1)})>0$ on $\mathcal T$. Since $\hat\varrho_n(t)=\{\hat\Lambda_n(\tau)-\hat\Lambda_n(t_0)\}-\{\hat\Lambda_n(t)-\hat\Lambda_n(t_0)\}$ only involves $\hat\lambda_n$ on $[t_0,\tau]$, Lemma~\ref{lem:int_rate} applies on $[\mathfrak s_0,\mathfrak s_1]=[t_0,\tau]$, without boundary correction, with $a=\gamma p$, $a_\Lambda=\gamma d_c$ (Propositions~\ref{prop:kernel}--\ref{prop:histogram}) and $r^R_n$ from \eqref{eq:rate_R}, and gives
\[
\sup_{t\in\mathcal T}|\hat\varrho_n(t)-\varrho(t)|=O\big(n^{-\mathfrak b}+w_n(\gamma d_c)+r^R_n\big)=o\big(w_n(\gamma p)\big)
\]
almost surely, under the conditions of Corollary~\ref{cor:exponents}(\ref{cor:exponents:i}) or (\ref{cor:exponents:iii}): $n^{-\mathfrak b}\leq n^{-(1-\gamma p)/2}$, $w_n(\gamma d_c)=o(w_n(\gamma p))$ since $d_c<p$, and $r^R_n=o(w_n(\gamma p))$ as in the proof of Corollary~\ref{cor:exponents}. Hence, for $n$ large, $\sup_{\mathcal T}|\hat R^{(0)}_n-R^{(0)}|\leq\sup|\hat R_{1,n}-R_1|+\|R_1\|_\infty\sup|e^{-\hat\varrho_n}-e^{-\varrho}|=o(w_n(\gamma p))$, and Theorem~\ref{theo:hazard}, through Proposition~\ref{prop:ratio}, applies to $\hat f_n/\hat R^{(0)}_n$ with the target $f(\cdot,1,\x)/R^{(0)}=\lambda_0(\cdot\mid\x)$ by \eqref{eq:cure_hazard}.

(\ref{theo:cure:ii}) Here $\hat\Lambda_n$ is computed on $[0,\tau]$, with the boundary correction for the kernel learner, and Lemma~\ref{lem:int_rate} on $[0,\tau]$ gives, under the conditions of Corollary~\ref{cor:exponents}(\ref{cor:exponents:ii}), $\delta_n:=\sup_{[0,\tau]}|\hat\Lambda_n-\Lambda|=O(n^{-\mathfrak b}+w_n(\gamma d_c)+r^R_n)=O(w_n(\gamma d_c))$ and $\sup_{[0,\tau]}|\mathfrak S^\Lambda_n|=O(nw_n(\gamma d_c))$ almost surely, where $\mathfrak S^\Lambda_n=\sum_{i\leq n}\xi^\Lambda_i$. Let $\mathcal U=\{(u,w): -\epsilon\leq u\leq\Lambda(t_1)+\epsilon,\ |w-\Lambda(\tau)|\leq\epsilon\}$ for a small $\epsilon>0$. On $\mathcal U$, $e^{-u}-e^{-w}\geq S(t_1)e^{-\epsilon}-S(\tau)e^{\epsilon}>0$ and $1-e^{-w}\geq1-S(\tau)e^{\epsilon}>0$, so that $\mathfrak G$ is twice continuously differentiable on a neighborhood of $\mathcal U$, with $\partial_u\mathfrak G(u,w)=e^{-u}/(e^{-u}-e^{-w})$ and $\partial_w\mathfrak G(u,w)=e^{-w}/(1-e^{-w})-e^{-w}/(e^{-u}-e^{-w})$; at $(u,w)=(\Lambda(t),\Lambda(\tau))$ these are $c_1(t)$ and $c_2(t)$ of \eqref{eq:cure_c1c2}, and $c_1(t)+c_2(t)=1/\{1-S(\tau)\}=1/\phi(\x)$. For $n$ large, $(\hat\Lambda_n(t),\hat\Lambda_n(\tau))\in\mathcal U$ for all $t\in\mathcal T$, and a second order Taylor expansion gives, uniformly in $t\in\mathcal T$,
\[
\hat\Lambda_{0,n}(t)-\Lambda_0(t)=c_1(t)\{\hat\Lambda_n(t)-\Lambda(t)\}+c_2(t)\{\hat\Lambda_n(\tau)-\Lambda(\tau)\}+O(\delta_n^2).
\]
Let $\mathfrak d(t)=\Lambda_0(t)-c_1(t)\Lambda(t)-c_2(t)\Lambda(\tau)$, a deterministic Lipschitz function on $\mathcal T$, which only enters the analysis. Inserting the decomposition $\hat\Lambda_n(s)-\Lambda(s)=n^{-1}\mathfrak S^\Lambda_n(s)+\{\bar\Lambda_n(s)-\Lambda(s)\}+\mathfrak r^\Lambda_n(s)$ for $s=t$ and $s=\tau$, we obtain the one-pass representation
\begin{gather*}
\hat\Lambda_{0,n}(t)=\frac1n\sum_{i=1}^n\psi^0_i(t;Z_i)+\mathfrak r^0_n(t),\\
\psi^0_i(t;\mathfrak z)=\mathfrak d(t)+\int_0^\tau w_t(s)\frac{\psi_i(s;\mathfrak z)}{R_1(s,\x)}\dd s=\mathfrak d(t)+c_1(t)\psi^\Lambda_i(t;\mathfrak z)+c_2(t)\psi^\Lambda_i(\tau;\mathfrak z),
\end{gather*}
with $w_t$ from \eqref{eq:cure_c1c2}, $\bar g^0_n(t)=\mathfrak d(t)+c_1(t)\bar\Lambda_n(t)+c_2(t)\bar\Lambda_n(\tau)$ and $\mathfrak r^0_n(t)=c_1(t)\mathfrak r^\Lambda_n(t)+c_2(t)\mathfrak r^\Lambda_n(\tau)+O(\delta_n^2)$. We check the conditions of   Theorem~\ref{theo:band}(\ref{theo:band:i})    on $\mathcal T$ for the target $\Lambda_0$, with the exponents $a_\Lambda=\gamma d_c$ and $\mathfrak b$.
Condition \eqref{cond:I}: $c_1,c_2$ are bounded on $\mathcal T$, so the bounds of Assumption~\ref{ass:cumhaz}(\ref{ass:cumhaz:iii}) transfer to $\xi^0_i=c_1\xi^\Lambda_i(t)+c_2\xi^\Lambda_i(\tau)$, and $\sup_t|\mathfrak r^0_n|=O(r^R_n)+O(\delta_n^2)=o(w_n(\gamma d_c))$, by the proof of Corollary~\ref{cor:exponents}. Condition \eqref{cond:B}: $\bar g^0_n(t)-\Lambda_0(t)=c_1(t)\{\bar\Lambda_n(t)-\Lambda(t)\}+c_2(t)\{\bar\Lambda_n(\tau)-\Lambda(\tau)\}=O(n^{-\mathfrak b})$ uniformly in $t$. Condition \eqref{cond:V}: for $s\in[0,\tau]$, $s\mapsto w_t(s)/R_1(s,\x)$ is of bounded variation, with total variation and supremum bounded uniformly in $t\in\mathcal T$, and the computations of the variance of the integrated contributions in the proofs of Propositions~\ref{prop:kernel}--\ref{prop:histogram}(\ref{prop:histogram:ii}) apply with this weight in place of $\1\{s\leq t\}R_1(s,\x)^{-1}$: for the kernel learner, \eqref{eq:Phi_b} holds with $\Phi_b$ replaced by $\int_0^\tau w_t(s)K^s_{B,b}(y-s)R_1(s,\x)^{-1}\dd s$ and $\1\{y\leq t\}R_1(y,\x)^{-1}$ by $w_t(y)R_1(y,\x)^{-1}$, the jump of $w_t$ at $t$ being handled exactly as the jump of the indicator; for the series learner, Lemma~\ref{lem:S3}(\ref{lem:S3:ii}) applies to $g=w_tR_1^{-1}$; for the histogram, the computation is unchanged. Hence,   $\Var\{\xi^0_i(t)\}=\mathfrak v_iv^0(t)\{1+o(1)\}$, with $\mathfrak v_i=i^{\gamma d_c}$ or $2^{k_id_c}$,   uniformly in $t\in\mathcal T$, with $v^0(t)=\mathfrak c_\Lambda\int_0^\tau w_t(s)^2f(s,1,\x)R_1(s,\x)^{-2}\dd s\geq\phi(\x)^{-2}\mathfrak c_\Lambda\int_0^{t_0}f(s,1,\x)R_1(s,\x)^{-2}\dd s>0$ by Assumption~\ref{ass:data}(\ref{ass:D3}), and $\mathfrak V^\Lambda_n$ as in Section~\ref{sec:examples}; $v^0$ is continuous since $c_1,c_2$ are, and $\hat v^0_n/v^0\to1$ uniformly because $\hat c_1,\hat c_2$ are uniformly consistent on $\mathcal T$ (as $\delta_n\to0$ and $S(t_1)>S(\tau)$) and $\int_0^\tau\hat w_t^2\hat f_n/\hat R_{1,n}^2\to\int_0^\tau w_t^2f/R_1^2$ uniformly. Condition (\ref{cond:M}): $c_1,c_2$ are Lipschitz on $\mathcal T$, and
\[
\mathfrak S^0_n(t)-\mathfrak S^0_n(s)=c_1(t)\{\mathfrak S^\Lambda_n(t)-\mathfrak S^\Lambda_n(s)\}+\{c_1(t)-c_1(s)\}\mathfrak S^\Lambda_n(s)+\{c_2(t)-c_2(s)\}\mathfrak S^\Lambda_n(\tau),
\]
where, for $s=\pi_n(t)$, the first term is $o(nw_n(\gamma d_c))$ as in the proof of Theorem~\ref{theo:cumhaz} and the other two are $O(\varsigma_n\sup_{[0,\tau]}|\mathfrak S^\Lambda_n|)=O(\varsigma_nnw_n(\gamma d_c))=o(nw_n(\gamma d_c))$. Condition \eqref{cond:E} follows from the expansion and $\delta_n\to0$.   Theorem~\ref{theo:band}(\ref{theo:band:i})    applies and gives the band for $\Lambda_0$; the band for $S_0$ follows by monotonicity, exactly as in the proof of Corollary~\ref{cor:survival}. The pointwise statements follow from   Theorem~\ref{theo:band}(\ref{theo:band:ii}).  
\end{proof}

\subsection{Further mechanisms}\label{app:further_mech}

\paragraph{Competing risks.} With $K$ competing causes, $T=\min_kT_k$, and right censoring by $C$, one observes $Y=T\wedge C$, with $I=k$ if $T=T_k\leq C$ and $I=0$ if $T>C$. It is assumed that $(T_1,\ldots,T_K)\perp C\mid\X$. The cause-specific hazard $\lambda_k(t\mid\x)=\lim_{h\downarrow0}h^{-1}\PP(t\leq T<t+h,I=k\mid T\geq t,\X=\x)$ has the representation $f_{Y,I,\X}(t,k,\x)/R_1(t,\x)$, with $R_1$ as in \eqref{eq:R_rc}. The numerator is estimated by any of the three learners applied with the indicator $\1\{I=k\}$, and Theorem~\ref{theo:hazard} gives anytime-valid bands for $\lambda_k$ without further work. The cumulative incidence function $F_k(t\mid\x)=\int_0^tS_{T\mid\X}(s\mid\x)\lambda_k(s\mid\x)\dd s$ and its derivative \citep{AP2026} can be handled in a similar way, with two changes. First, $S_{T\mid\X}$ is estimated by $\exp(-\hat\Lambda_n)$, integrated from the origin, which requires the boundary-corrected learner. Second, the cumulative incidence function requires a delta method combining the numerators of cause $k$ and of all causes. We do not detail these extensions. See  \cite{FG1999} and  \cite{G2024} for the background.

\vspace{-.3cm}

\paragraph{Twice censoring.} For data where $T$ is observed exactly, or only known to be larger, or only known to be smaller than the observed duration, the two twice-censoring models of \cite{PR2006a}, with $Y=\max\{C_{\mathrm l},\min(T,C_{\mathrm r})\}$ or $Y=\min\{\max(T,C_{\mathrm l}),C_{\mathrm r}\}$, admit explicit inverses as well; the at-risk functions are more involved functionals of the sub-distributions, and we do not detail them.

\vspace{-.3cm}

\paragraph{Reverse hazards.}   The approach extends to   the reverse hazard function $r_{T\mid\X}=f_{T\mid\X}/F_{T\mid\X}$, which characterizes the conditional distribution of $T$ through $F_{T\mid\X}(t\mid\x)=\exp\{-\int_t^\infty r_{T\mid\X}(s\mid\x)\dd s\}$. The reverse hazard is the natural object when the maximum of $T$ and a censoring variable is observed, as with left censoring, $Y=T\vee C$, for which $r_{T\mid\X}(t\mid\x)=f_{Y,I,\X}(t,1,\x)/\{\PP(Y\leq t\mid\X=\x)f_\X(\x)\}$ \citep{AP2026}. The second latent model of \cite{PR2006b} is the counterpart of modified current status data for this situation, and is obtained by a time reversal. Its at-risk function is affine in the nuisance, as in \eqref{eq:R_mcs}: with $I=1$ for the exact observations and $I=0$ for the left-censored ones, $R_1(t,\x)=\PP(Y\leq t,I=1\mid\X=\x)f_\X(\x)$, $R_2(t,\x)=\PP(Y\leq t,I=0\mid\X=\x)f_\X(\x)$ and $\eta=\mathfrak p(\x)$, estimated by \eqref{eq:hat_p}. The left-censoring denominator $\PP(Y\leq t\mid\X=\x)f_\X(\x)$ applies here only when $\mathfrak p\equiv1$. Model II of \cite{PR2006a}, the counterpart of twice censoring, also admits an explicit inverse, but its at-risk function is not of regression type. For left censoring and the second model of \cite{PR2006b}, the at-risk components are of the form $\EE[u(t;Z)\mid\X=\x]f_\X(\x)$ with non-decreasing indicators such as $u(t;\mathfrak z)=\1\{y\leq t,\iota=1\}$, and Assumption~\ref{ass:data}(\ref{ass:D4}) is unaffected. The cumulative reverse hazard integrates from $t$ to a known right end of the support of $T$ given $\X=\x$, an identification assumption analogous to sufficient follow-up, and the boundary correction of the kernel learner is needed at this point instead of the origin. With these changes, the results of Sections~\ref{sec:inference} and \ref{sec:examples} apply to these two mechanisms.

\section{The case without continuous covariates}\label{app:dc0}

Let $d_c=0$: the covariate vector is discrete, $\X=\X_d$, or absent, and $\PP(\X=\x)>0$. We consider right censoring, with $f(t)=f_{Y,I,\X}(t,1,\x)$ and $R(t)=\PP(Y\geq t,\X=\x)$, so that $\lambda=f/R$, where $\lambda$ and $\Lambda$ stand for $\lambda_{T\mid\X}(\cdot\mid\x)$ and $\Lambda_{T\mid\X}(\cdot\mid\x)$; the other mechanisms are handled in a similar way.

\vspace{-.3cm}

\paragraph{Bounded contributions.} The integrated contributions now have the growth exponent $a_\Lambda=\gamma d_c=0$, and the results of Section~\ref{sec:inference} extend to $a=0$, with $w_n(0)=n^{-1/2}\{\log(n+1)\}^{1/2}$. In the proof of Lemma~\ref{lemma:blocked}, $\lambda_k=(2C_M)^{-1}$ does not depend on $k$, and the sum over the blocks gives, on $E_u$, $|\sum_{j\leq n}(\xi_j-\eta_j)|\leq(2C_M/c)(K+1)\{2(K+1)+u\}\leq C\{\log^2(n+1)+u\log(n+1)\}$. In Step~1 of the proof of Theorem~\ref{theo:band}(\ref{theo:band:i}), the bound \eqref{eq:unif_coupling} becomes $O(\log^2n)$, which is $o(nw_n(0))$; the other steps, as well as Lemma~\ref{lemma:chaining}, only use $a<1$. Hence Theorem~\ref{theo:band} holds with $a=0$.   The same holds for Theorem~\ref{theo:cumhaz} with $a_\Lambda=0$ in Assumption~\ref{ass:cumhaz}(\ref{ass:cumhaz:iii}): besides Theorem~\ref{theo:band}, its proof and that of Lemma~\ref{lem:int_rate} only use $a_\Lambda\leq a<1$. 

\vspace{-.3cm}

\paragraph{The hazard function.} The numerator is smoothed in time only: Proposition~\ref{prop:kernel}(\ref{prop:kernel:i}) holds with $p=1$, $a=\gamma$, $\mathfrak b=2\gamma$ and $\mathfrak c_K=c_K/c_b$, and similarly for the series and histogram learners. The at-risk function is estimated without smoothing, by $\hat R_n(t)=n^{-1}\sum_{i\leq n}\1\{Y_i\geq t,\X_i=\x\}$, and $\sup_t|\hat R_n-R|=O(w_n(0))=o(w_n(\gamma))$ almost surely, by Lemma~\ref{lemma:chaining}(\ref{lemma:chaining:b}).   Theorem~\ref{theo:hazard} applies for $1/5\leq\gamma<1$, with the width $n^{-2/5}\sqrt{\log n}$ at $\gamma=1/5$,   for the kernel and series learners, and for $1/3\leq\gamma<1$, with the width $n^{-1/3}\sqrt{\log n}$ at $\gamma=1/3$, for the histogram learner.

\vspace{-.3cm}

\paragraph{The cumulative hazard and the denominator.} For $d_c\geq1$, the error of $\hat R_n$ is of smaller order than the integrated contributions $\psi^\Lambda_i$, which only involve the numerator. For $d_c=0$, both are of order $n^{-1/2}$, up to logarithmic factors. 
Let $W_i=Y_i$ if $I_i=1$ and $\X_i=\x$, and $W_i=+\infty$ otherwise; let $W'_i=Y_i$ if $\X_i=\x$, and $W'_i=+\infty$ otherwise. Then $H_1(t)=\PP(W\leq t)$ and $R(t)=\PP(t\leq W'<\infty)$ for $t\geq0$. The law of the iterated logarithm form of the time-uniform Dvoretzky--Kiefer--Wolfowitz inequality, \citet[Theorem~2]{HR2022}, which holds for any distribution on a totally ordered set, applied to $(W_i)$ and $(W'_i)$, gives, almost surely,
\begin{equation}\label{eq:dkw_lil}
	\delta_n:=\sup_{t\geq0}|\hat H_{1,n}(t)-H_1(t)|+\sup_{t\geq0}|\hat R_n(t)-R(t)|=O\big(\{n^{-1}\log\log n\}^{1/2}\big),
\end{equation}
with $H_1(t)=\PP(Y\leq t,I=1,\X=\x)$ and $\hat H_{1,n}(t)=n^{-1}\sum_{i\leq n}\1\{Y_i\leq t,I_i=1,\X_i=\x\}$. Hence $\sup_t|\hat R_n-R|=o(w_n(0))$, Assumption~\ref{ass:cumhaz}(\ref{ass:cumhaz:iv}) holds, and Theorem~\ref{theo:cumhaz} applies to the integrated learners, with $a_\Lambda=0$, and $\gamma\geq1/4$ for the kernel and series learners or $\gamma\geq1/2$ for the histogram learner. However, the error of $\hat R_n$ is then smaller than the width only by the factor $\{\log\log n/\log n\}^{1/2}$. Moreover, $\EE[\psi^\Lambda_i(t;Z)]^2$ is no longer negligible with respect to $\EE[\psi^\Lambda_i(t;Z)^2]$, and the variance profile of Proposition~\ref{prop:kernel}(\ref{prop:kernel:ii}) becomes $\int_0^tf/R^2-\Lambda(t)^2$.   It is therefore preferable to include the linear part of the error of the denominator in the contributions,   that is, to add $-\1\{\z=\x\}\int_0^{t\wedge y}f/R^2+\Lambda(t)$ to $\psi^\Lambda_i(t;\mathfrak z)$. These terms are bounded and centered, and are non-increasing in $t$ up to the Lipschitz function $\Lambda$,   so that Conditions \eqref{cond:I}, \eqref{cond:B}, \eqref{cond:E} and \eqref{cond:M} hold for them, with $a=0$, as in the proof of Proposition~\ref{prop:NA}.   The remainder then becomes
\[
-\int_0^t\frac{(\hat f_n-f)(\hat R_n-R)}{R^2}+\int_0^t\frac{\hat f_n(\hat R_n-R)^2}{R^2\hat R_n}=O\Big(\delta_n\sup_{[0,t_1]}|\hat f_n-f|+\delta_n^2\Big)=o(w_n(0))
\]
almost surely, by \eqref{eq:dkw_lil} and Conditions \eqref{cond:E} and \eqref{cond:B} for $\hat f_n$ on $[0,t_1]$. The variance profile tends to the classical one,  $\int_0^tf/R^2$, and smoothing becomes optional.  

\begin{proposition}[Nelson--Aalen estimator]\label{prop:NA}
Let $d_c=0$, let $f$ be bounded on $[0,t_1]$, with $\int_0^{t_0}f>0$, let $R(t_1)>0$, and let the grid be such that $\varsigma_n\to0$ and \eqref{eq:ell_n} holds. Then the Nelson--Aalen estimator
\[
\hat\Lambda^{\mathrm{NA}}_n(t)=\int_0^t\frac{\hat H_{1,n}(\dd s)}{\hat R_n(s)},\qquad t\in[0,t_1],
\]
satisfies Conditions \eqref{cond:I}, \eqref{cond:V}, \eqref{cond:B}, \eqref{cond:E} and \eqref{cond:M} on $\mathcal T^\Lambda=[t_0,t_1]$ for $\Lambda$, with $a=0$, any $\mathfrak b>0$, $\mathfrak V_n=n$,
\[
v^{\mathrm{NA}}(t)=\int_0^t\frac{f(s)}{R(s)^2}\dd s\qquad\text{and}\qquad\hat v^{\mathrm{NA}}_n(t)=\int_0^t\frac{\hat H_{1,n}(\dd s)}{\hat R_n(s)^2}.
\]
Consequently, the bands $[\hat\Lambda^{\mathrm{NA}}_n(t)\pm\hat L^{\mathrm{NA}}_n(t)]$, with $\hat L^{\mathrm{NA}}_n(t)=n^{-1}\mathcal B_\rho(n\hat v^{\mathrm{NA}}_n(\pi_n(t)),\alpha_{\pi_n(t)})$, form a $(1-\alpha)$-asymptotic confidence band sequence for $\Lambda$ on $\mathcal T^\Lambda$, and the bands $[e^{-\hat\Lambda^{\mathrm{NA}}_n(t)-\hat L^{\mathrm{NA}}_n(t)},e^{-\hat\Lambda^{\mathrm{NA}}_n(t)+\hat L^{\mathrm{NA}}_n(t)}]$ one for $S_{T\mid\X}(\cdot\mid\x)$. If $\ell_n=o(\log n)$, then $n\hat L^{\mathrm{NA}}_n(t)^2\sim\hat v^{\mathrm{NA}}_n(\pi_n(t))\log n$.
\end{proposition}

\begin{proof}
By \eqref{eq:dkw_lil}, $\inf_{[0,t_1]}\hat R_n\geq R(t_1)/2$ for $n$ large, almost surely. Expanding $1/\hat R_n=1/R-(\hat R_n-R)/R^2+(\hat R_n-R)^2/(R^2\hat R_n)$ gives $\hat\Lambda^{\mathrm{NA}}_n(t)=n^{-1}\sum_{i\leq n}\psi^{\mathrm{NA}}_i(t;Z_i)+\mathfrak r_n(t)$, with
\begin{gather*}
\psi^{\mathrm{NA}}_i(t;\mathfrak z)=\Lambda(t)+\1\{\z=\x\}\Big\{\frac{\1\{\iota=1,\,y\leq t\}}{R(y)}-\int_0^{t\wedge y}\frac{f(s)}{R(s)^2}\dd s\Big\},\\
\mathfrak r_n(t)=\int_0^t\frac{R-\hat R_n}{R\hat R_n}\,\dd(\hat H_{1,n}-H_1)+\int_0^t\frac{(\hat R_n-R)^2}{R^2\hat R_n}\,\dd H_1.
\end{gather*}
The second term of $\mathfrak r_n$ is $O(\delta_n^2)$. In the first one, $\varphi_n=(R-\hat R_n)/(R\hat R_n)$ satisfies $\sup_{[0,t_1]}|\varphi_n|\leq C\delta_n$ and has a total variation on $[0,t_1]$ bounded by a constant, since $R$ and $\hat R_n$ are monotone and bounded away from zero; an integration by parts gives a bound $C\delta_n$. Hence $\sup_t|\mathfrak r_n(t)|=O(\delta_n)=o(w_n(0))$ almost surely.

The centered contributions are $\xi_i(t)=\int_0^tR(s)^{-1}\dd M_i(s)$, with $M_i(s)=\1\{\X_i=\x\}\{\1\{Y_i\leq s,I_i=1\}-\int_0^s\1\{Y_i\geq u\}\lambda(u)\dd u\}$, which is a martingale in $s$ under $T\perp C\mid\X$. Hence $\EE\xi_i(t)=0$, so that Condition \eqref{cond:B} holds without bias, and
\[
\Var\{\xi_i(t)\}=\int_0^t\frac{\EE[\1\{Y\geq s,\X=\x\}]\lambda(s)}{R(s)^2}\dd s=\int_0^t\frac{f(s)}{R(s)^2}\dd s=v^{\mathrm{NA}}(t).
\]
The contributions are bounded by $\Lambda(t_1)+R(t_1)^{-1}+t_1\sup f/R(t_1)^2$, which gives Condition \eqref{cond:I} with $a=0$, and $V_n=nv^{\mathrm{NA}}$ exactly. The function $v^{\mathrm{NA}}$ is continuous and bounded below by $\int_0^{t_0}f/R^2>0$ on $\mathcal T^\Lambda$, and $\sup_t|\hat v^{\mathrm{NA}}_n-v^{\mathrm{NA}}|=O(\delta_n)$ by the same integration by parts, which gives Condition \eqref{cond:V}.
	
For Conditions \eqref{cond:E} and \eqref{cond:M}, write $\psi^{\mathrm{NA}}_i-\Lambda=\psi^{(1)}-\psi^{(2)}$, with $\psi^{(1)}(t;\mathfrak z)=\1\{\z=\x,\iota=1,y\leq t\}/R(y)$ and $\psi^{(2)}(t;\mathfrak z)=\1\{\z=\x\}\int_0^{t\wedge y}f/R^2$. Both are non-decreasing in $t$ and bounded, with Lipschitz expectations, and Lemma~\ref{lemma:chaining}(\ref{lemma:chaining:b}) gives Condition \eqref{cond:E}. Their increments satisfy $|\psi^{(j)}(t';\mathfrak z)-\psi^{(j)}(t;\mathfrak z)|\leq C$ and $\EE[\{\psi^{(j)}(t';Z)-\psi^{(j)}(t;Z)\}^2]\leq C|t'-t|$, since $f$ is bounded. Let $\varepsilon>0$, let $\ell'_n$ be as in the proof of Lemma~\ref{lemma:chaining}(\ref{lemma:chaining:d}), and let $\mathfrak S^{(j)}_n$ be the partial sums of the centered $\psi^{(j)}(\cdot;Z_i)$. For $t,t'\in G_{\ell'_n}$ with $|t'-t|\leq2\varsigma_n$, Bernstein's inequality gives
\[
\PP\big\{|\mathfrak S^{(j)}_n(t')-\mathfrak S^{(j)}_n(t)|\geq\varepsilon\sqrt{n\log n}\big\}\leq2\exp\Big(-\frac{\varepsilon^2n\log n}{2C\{2n\varsigma_n+\varepsilon\sqrt{n\log n}\}}\Big)\leq2(n+1)^{-6}
\]
for $n$ large, since $\varsigma_n\to0$. The union bound over the at most $\{2(n+1)^2+1\}^2$ pairs and the Borel--Cantelli lemma, applied for each rational $\varepsilon$, show that the maximum of these increments is $o(\sqrt{n\log n})$ almost surely. For $t\in\mathcal T^\Lambda$, let $s_-\leq t\leq s_+$ be consecutive points of $G_{\ell'_n}$. By monotonicity, as in the proof of Lemma~\ref{lemma:chaining}(\ref{lemma:chaining:b}), $\mathfrak S^{(j)}_n(s_-)-Cn(n+1)^{-2}\leq\mathfrak S^{(j)}_n(t)\leq\mathfrak S^{(j)}_n(s_+)+Cn(n+1)^{-2}$, and $\pi_n(t)$, $s_-$ and $s_+$ are within $2\varsigma_n$ of each other. Hence $\sup_t|\mathfrak S_n(t)-\mathfrak S_n(\pi_n(t))|=o(\sqrt{n\log n})=o(nw_n(0))$, which is Condition \eqref{cond:M}.
	
Theorem~\ref{theo:band}(\ref{theo:band:i}), with $a=0$, and the proof of Corollary~\ref{cor:survival} give the bands. Finally, if $\ell_n=o(\log n)$, then $\log(1/\alpha_s)=o(\log n)$ for $s\in G_{\ell_n}$, and $n^2\hat L^{\mathrm{NA}}_n(t)^2\sim n\hat v^{\mathrm{NA}}_n(\pi_n(t))\log n$ by \eqref{eq:robbins_b}.
\end{proof}

\paragraph{  Two routes for the cumulative hazard.}   When $d_c=0$, bands of width of order $n^{-1/2}\sqrt{\log n}$ for $\Lambda$ and $S$ can thus be obtained in two ways: from the integrated learners of Section~\ref{sec:examples}, with the linear correction of the contributions described above, or from the Nelson--Aalen estimator. The first route has the price of smoothing: a tuning sequence, a smoothness condition on $f$, a boundary correction at the origin for the kernel learner, and a bias, controlled for $\gamma\geq1/4$ (kernel and series learners) or $\gamma\geq1/2$ (histogram learner). Its only advantage is that the same learner also estimates $\lambda$, but only under these constraints, which exclude the exponents giving the narrowest bands for $\lambda$, namely $\gamma=1/5$ and $\gamma=1/3$. The Nelson--Aalen estimator requires no tuning, no smoothness of $f$ beyond boundedness and no boundary correction; it has no bias and the classical variance, and Condition~\eqref{cond:M} holds for any grid with $\varsigma_n\to0$. It admits a one-pass expansion, which is all the inference requires (Section~\ref{subsec:onepass}), but its exact computation is not one-pass: $\hat R_n$ enters at all the observed event times, and each new observation may change it at all the earlier ones, so that the observed times have to be stored. The gap is small, however. The integrals of the first route are computed on a fixed grid of time points, chosen for a planned horizon (Section~\ref{subsubsec:memory}); on such a grid, the Nelson--Aalen estimator is computed in one pass as well, from the numbers of events and of individuals at risk between consecutive grid points, with a quadrature error of the same type. For $\Lambda$ and $S$, we therefore recommend the Nelson--Aalen estimator, exact or binned, and a smoothed learner for $\lambda$ only. 

\vspace{-.3cm}

\paragraph{  The Kaplan--Meier estimator.}   The observed event times are almost surely distinct, since $H_1$ has a density. The Kaplan--Meier estimator $\hat S^{\mathrm{KM}}_n(t)$ is the product of the factors $1-x_j$, with $x_j=1/\{n\hat R_n(Y_j)\}$, over the observed event times $Y_j\leq t$ with $\X_j=\x$, while $\hat\Lambda^{\mathrm{NA}}_n(t)=\sum_jx_j$. Since $0\leq-\log(1-x)-x\leq x^2/(1-x)$ for $x\in[0,1)$, and $n\hat R_n(Y_j)\geq n\hat R_n(t)$, we have, when $n\hat R_n(t)\geq2$,
\[
0\leq-\log\hat S^{\mathrm{KM}}_n(t)-\hat\Lambda^{\mathrm{NA}}_n(t)\leq\varepsilon_n(t):=\frac{\hat\Lambda^{\mathrm{NA}}_n(t)}{n\hat R_n(t)-1},
\]
and we set $\varepsilon_n(t)=+\infty$ otherwise. By \eqref{eq:dkw_lil}, $\sup_{[0,t_1]}\varepsilon_n=O(n^{-1})$ almost surely. Hence, when $\Lambda(t)$ lies in the band of Proposition~\ref{prop:NA}, $S_{T\mid\X}(t\mid\x)$ lies in $[\hat S^{\mathrm{KM}}_n(t)e^{-\hat L^{\mathrm{NA}}_n(t)},\hat S^{\mathrm{KM}}_n(t)e^{\hat L^{\mathrm{NA}}_n(t)+\varepsilon_n(t)}]$, which is therefore a $(1-\alpha)$-asymptotic confidence band sequence centered at the Kaplan--Meier estimator, with a correction of smaller order than the width.

\vspace{-.3cm}

\paragraph{Comparison with time-uniform bands for distribution functions.} Without censoring and without covariates, the empirical distribution function is a one-pass   average, with the contributions $\1\{y\leq t\}$, $a=0$ and no bias, and the arguments of the proof of Proposition~\ref{prop:NA} give Conditions \eqref{cond:I}--\eqref{cond:M} for any grid with $\varsigma_n\to0$. Theorem~\ref{theo:band} then gives bands for $F$ whose half-width is of order $\{F(t)(1-F(t))\log n/n\}^{1/2}$, adapted to the variance, while \cite{HR2022} obtain bands of constant width shrinking at the rate $\{n^{-1}\log\log n\}^{1/2}$, which is the best possible by the law of the iterated logarithm. The additional factor $\{\log n/\log\log n\}^{1/2}$ comes from the normal mixture boundary, whose square grows like $\log V_n\asymp\log n$; with a grid growing slowly, $\ell_n=o(\log n)$, the allocation \eqref{eq:alpha_s} only adds $o(\log n)$. We expect that the stitched boundaries of \cite{HRMS2021}, applied to the Gaussian partial sums of Step~2 of the proof of Theorem~\ref{theo:band}(\ref{theo:band:i}), would reduce this factor when $d_c=0$, at the price of sharper bounds on the remainders; for $d_c\geq1$, the factor $\sqrt{\log n}$ is due to the union bound over a grid at the scale of the bandwidth in any case. Under right censoring, \cite{HZR2026} assume that the censoring distribution is known; Proposition~\ref{prop:NA} covers an unknown censoring distribution, estimated from the same data, and discrete covariates.   Bands of the order $\{n^{-1}\log\log n\}^{1/2}$ for $\Lambda$ and $S$ are also available, although neither adapted to the variance nor sharp in the constants: on the event where the bands of \cite{HR2022} for $H_1$ and $R$, each at level $1-\alpha/2$, hold for all $n$, the integration by parts in the proof of Proposition~\ref{prop:NA} bounds $\sup_{[0,t]}|\hat\Lambda^{\mathrm{NA}}_n-\Lambda|$ by an explicit multiple of their widths, of order $1/\hat R_n(t)$. At a fixed sample size,  \cite{BLM1999} prove a Dvoretzky--Kiefer--Wolfowitz inequality for the Kaplan--Meier estimator, uniform on the whole line thanks to the weight $1-G(t-)$, with $G$ the censoring distribution.

\section{Additional simulation results: a stress design}\label{app:stress}

The stress design differs from the design of Section~\ref{sec:numerics} by a stronger dependence between the continuous covariates, a stronger curvature in $x_1$, heavier censoring and a narrow   peak: the correlation of the copula is $0.8$, $g(x_1,x_2)=-0.9+0.8\sin(2\pi x_1)-0.8x_2^2$, the rate of $E$ is $0.35+0.5X_1$, which gives about $64\%$ of observations with $I=0$, and $h_1(t)=0.25+2.2\exp\{-(t-1)^2/(2\cdot0.18^2)\}$, with low values of the hazard outside the peak. The fourth covariate value is $(0.3,0.7,0)$; with $q=0.15$, about $14\%$ of the event times are not recalled, that is $5\%$ of the observations. The procedures are those of Section~\ref{sec:numerics}. Table~\ref{tab:sim_stress} gives the results. 

\begin{table}[ht]
\centering
\footnotesize
\setlength{\tabcolsep}{4pt}
\begin{tabular}{lll|ccccc|cc}
		\hline
		& & & \multicolumn{5}{c|}{Largest miscoverage over the four $\x$} & \multicolumn{2}{c}{Ratio to the oracle width}\\
		$\gamma$ & $q$ & Target & Fixed-$n$ & Fine & Fine, l.m. & Band. & Band., l.m. & Fine, l.m. & Band., l.m.\\
		\hline
		$1/6$ & $0$    & $\lambda$ & 0.507 & 0.009 & 0.000 & 0.640 & 0.132 & 2.37--3.03 & 1.24--1.62\\
		$1/6$ & $0$    & $S$       & 0.429 & 0.000 & 0.000 & 0.500 & 0.266 & 1.99--3.18 & 1.03--1.65\\
		$1/6$ & $0.15$ & $\lambda$ & 0.521 & 0.014 & 0.000 & 0.617 & 0.127 & 2.44--3.13 & 1.28--1.67\\
		$1/6$ & $0.15$ & $S$       & 0.416 & 0.000 & 0.000 & 0.374 & 0.175 & 2.14--3.29 & 1.11--1.71\\
		$1/5$ & $0$    & $\lambda$ & 0.554 & 0.132 & 0.011 & 0.447 & 0.034 & 2.90--3.20 & 1.58--1.75\\
		$1/5$ & $0$    & $S$       & 0.484 & 0.000 & 0.000 & 0.051 & 0.009 & 3.00--3.74 & 1.62--2.01\\
		$1/5$ & $0.15$ & $\lambda$ & 0.556 & 0.161 & 0.013 & 0.509 & 0.054 & 2.95--3.20 & 1.60--1.76\\
		$1/5$ & $0.15$ & $S$       & 0.474 & 0.000 & 0.000 & 0.043 & 0.005 & 3.08--3.86 & 1.66--2.07\\
		\hline
\end{tabular}
\caption{Stress design. Same quantities as in Table~\ref{tab:sim}.}\label{tab:sim_stress}
\vspace{-.5cm}
\end{table}

\paragraph{The tuning and the denominator.} With $\gamma=1/6$, the band with the local-maximum variance and the grid at the scale of the bandwidth reaches a miscoverage of $0.132$ for the hazard and $0.266$ for the survival function under right censoring. Both are attained at the fourth covariate value, which lies against the correlation of $(X_1,X_2)$. There, the estimates are below the target over the whole of $\mathcal T$, for both targets, and the main cause is the bias of the denominator. The denominator is smoothed in $\x_c$ with $g(i)$, about twice $h(i)$ in this tuning, so that the curvature of $R(t,\cdot)$ in $\x_c$ no longer cancels between the numerator and the denominator. A computation of the expectations by quadrature, for the smoothing in $\x_c$, gives at $n=10^4$ a relative bias of $\hat R_n$ between $8\%$ and $63\%$ on $\mathcal T$, and a relative bias of the hazard estimator between $-34\%$ and $-11\%$. With $g=h$, the latter would be between $-13\%$ and $-7\%$. At the third covariate value, the failures for the hazard, $0.104$, occur at the peak, $t\approx0.94$, with the estimate below the target. There, the bias of the denominator partly offsets the flattening of the   peak   by the smoothing in time. The bias of the ratio thus depends on both bandwidths and on the covariate value. For the survival function, $\gamma=1/6$ is the rate-optimal exponent, so that the bias of the numerator is already of the order of the standard deviation, and that of the denominator adds to it. With $q=0.15$, these miscoverages decrease to $0.127$ and $0.175$: recall increases the variance, hence the width, by the factor $\mathfrak p(\x)^{-1/2}$, but not the bias. With $\gamma=1/5$, the bias is of smaller order, and the extra width of the fine grid provides a safety margin: the band with the local-maximum variance and the fine grid holds in all the configurations.

\vspace{-.3cm}

\paragraph{Recall and the right end of $\mathcal T$.} With $\gamma=1/5$ and $q=0.15$, the miscoverage of the recommended band for the hazard reaches $0.054$, against $0.034$ under right censoring. The largest miscoverage occurs at the second covariate value, where all the failures are early in the stream, with a median first exit near $n=1.8\times10^4$, and at the right end of $\mathcal T$, where few subjects are at risk under $64\%$ of censoring, so that the error of the at-risk estimator is not yet negligible; recall makes it larger, through the smaller number of exactly observed durations and the estimation of $\mathfrak p$. Two checks confirm this explanation. With $\mathcal T=[0.25,1.5]$, the miscoverage drops to $0.024$, and to $0$ at the second covariate value. With $\mathcal T=[0.25,1.75]$ and the bands monitored from $n=2\times10^4$ instead of $10^4$, as allowed by Proposition~\ref{prop:stopping}, it is $0.035$, and $0.020$ from $n=5\times10^4$. The local-maximum window, one-sided at $t_1$, may also contribute, but it is one-sided in the first check as well. A level $q=0.30$, not shown,   leads to the same conclusions: the width is multiplied by $\mathfrak p(\x)^{-1/2}$ within $1\%$, the failures with $\gamma=1/6$ decrease further, and those with $\gamma=1/5$ reach $0.076$, and $0.020$ when the bands are monitored from $n=5\times10^4$.

\end{document}